\documentclass[12pt]{article}

\newtheorem{definition}{Definition}[section]
\newtheorem{theorem}[definition]{Theorem}
\newtheorem{lemma}[definition]{Lemma}
\newtheorem{proposition}[definition]{Proposition}
\newtheorem{corollary}[definition]{Corollary}
\newtheorem{remark}[definition]{Remark}
\newtheorem{problem}{Problem}

\usepackage{amssymb}
\usepackage{amsmath}
\usepackage{latexsym}

\newcommand{\C}{{\mathbb C}}

\newcommand{\Z}{{\mathbb Z}}
\newcommand{\N}{{\mathbb N}}
\newcommand{\T}{\mathcal{T}}
\newcommand{\A}{\mathcal{A}}
\newcommand{\B}{\mathcal{B}}
\newcommand{\LL}{\mathcal{L}}
\newcommand{\I}{{\mathcal I}}
\newcommand{\J}{{\mathcal J}}

\newcommand{\W}{{\mathcal W}}

\begin{document}

\title{The augmented tridiagonal algebra}

\author{
Tatsuro Ito{\footnote{
Division of Mathematical and Physical Sciences,
Kanazawa University,
Kakuma-machi,
Kanazawa 920-1192, Japan
}}
{\footnote{Supported in part by JSPS grant 18340022}} 
$\;$ and
Paul Terwilliger{\footnote{
Department of Mathematics, University of
Wisconsin, 480 Lincoln Drive, Madison WI 53706-1388 USA}
}
\\
\\
{\small 
Dedicated to Professor Eiichi Bannai on the occasion of his retirement}
}

\date{}
\maketitle

\begin{abstract} 
Motivated by investigations of the tridiagonal pairs of linear 
transformations, we introduce the augmented tridiagonal algebra 
$\T_q$. This is an infinite-dimensional associative $\C$-algebra 
with $1$. 
We classify the finite-dimensional 
irreducible representations of $\T_q$. 
All such representations are explicitly constructed via embeddings 
of $\T_q$ into the $U_q(sl_2)$-loop algebra.
As an application, tridiagonal pairs over $\C$ are classified 
in the case where $q$ is not a root of unity. 

\medskip
\noindent 
{\bf Keywords}. P- and Q-polynomial association scheme, 
Terwilliger algebra, tridiagonal pair, Leonard pair, 
tridiagonal algebra, $q$-Onsager algebra. 
\hfil\break
\noindent {\bf 2000 Mathematics Subject Classification}. 
Primary: 17B37. 
Secondary: 05E30, 
33D45. 

\end{abstract}

\section{Introduction}
The purpose of this paper is to introduce the augmented tridiagonal 
algebra $\T_q$ and classify its finite-dimensional irreducible 
representations. We explain our motivations in Sections~1.1,~1.2 and 
summarize our results in Sections 1.3, 1.4.
Throughout this paper, we choose the complex number field 
$\C$ as the ground field. 
An algebra means an associative $\C$-algebra with $1$.

\subsection{Tridiagonal pairs: a background in combinatorics}
The standard generators for the subconstituent algebra 
(Terwilliger algebra) of a P- and Q-polynomial association scheme
\cite{BI}
give rise to a tridiagonal pair
when they are restricted to an irreducible submodule of 
the standard module \cite[Example 1.4]{ITT},
 \cite[Lemmas 3.9, 3.12]{Ta}. 
This fact motivates our ongoing investigation of
the tridiagonal pairs
\cite{ITT},
\cite{ITa},
\cite{ITb},
\cite{ITc},
\cite{ITd}, 
\cite{ITe}.

\medskip 
Let $V$ denote a finite-dimensional nonzero vector space over $\C$. 
Let ${\rm{End}}(V)$ denote the $\C$-algebra of all
$\C$-linear transformations of $V$.
By a \textit{tridiagonal pair} (\textit{TD-pair}) on $V$
we mean an ordered pair $A$, $A^*$ of elements in 
${\rm{End}}(V)$ that satisfy (i)--(iv) below:

\begin{enumerate}
\item
$A$ and $A^*$ are diagonalizable.
\item
There exists an ordering $V_0, V_1,\ldots, V_d$ of the  
eigenspaces of $A$ such that 
\begin{eqnarray*}
A^* V_i \subseteq V_{i-1} + V_i+ V_{i+1} \qquad \qquad (0 \leq i \leq d),
\end{eqnarray*}
where $V_{-1} = 0$, $V_{d+1}= 0$.
\item 
There exists an ordering $V^*_0, V^*_1,\ldots, V^*_{d^*}$ of the 
eigenspaces of $A^*$ such that 
\begin{eqnarray*}
A V^*_i \subseteq V^*_{i-1} + V^*_i+ V^*_{i+1} \qquad \qquad 
(0 \leq i \leq d^*), 
\end{eqnarray*}
where $V^*_{-1} = 0$, $V^*_{d^*+1}= 0$.
\item
$V$ is irreducible as an $\langle A, A^* \rangle $-module, 
where $\langle A, A^* \rangle $ is the subalgebra of 
${\rm{End}}(V)$
generated by $A$, $A^*$. 
\end{enumerate}

\noindent
A TD-pair $A$, $A^* \in {\rm{End}}(V)$ is \textit{isomorphic} 
to a TD-pair 
$B, B^* \in {\rm{End}}(V')$ whenever there exists an isomorphism 
$\psi: V \rightarrow V'$ of vector spaces such that 
$B \psi = \psi A$ and $B^* \psi = \psi A^*$.

\medskip 
In this subsection, we summarize the basic properties of a TD-pair 
$A$, $A^*$; they will be used to introduce the augmented tridiagonal algebra 
$\mathcal{T}_q$ in the next subsection. 
First we remark that $A$  and $A^*$ have the same number of eigenvalues, 
i.e. $d=d^*$ [2, Lemma~4.5]: we call this common 
integer the \textit{diameter} of the pair. 
A TD-pair with $d=0$ is called \textit {trivial}. 
We usually assume  $d \geq 1$ to avoid the trivial TD-pairs. 
Under this assumption, 
there exist exactly two orderings of the eigenspaces 
of $A$ (resp. $A^*$) that satisfy the condition (ii) (resp. (iii)): 
if 
$V_0, V_1,\ldots, V_d$ (resp. $V^*_0, V^*_1,\ldots, V^*_d$) 
is one of these, then the other is the reversed ordering 
$V_d, V_{d-1},\ldots, V_0$ (resp. $V^*_d, V^*_{d-1}, \ldots , V^*_0 $). 
We understand that one of such orderings is chosen and fixed 
unless otherwise stated. 

\medskip 
By \cite[Theorem 10.1]{ITT} there exist scalers 
$ \beta,\, \gamma,\, \gamma^*, \delta,\, \delta^*$ in $\C$ such that 
\begin{eqnarray*}
\lbrack A, ~A^2A^*-\beta A\,A^*A + A^*A^2\rbrack 
&=&\gamma\, \lbrack A, ~A\,A^*+A^*A\rbrack 
+ \delta\, \lbrack A, ~A^*\rbrack,\\
\lbrack A^*, ~A^{*2}A-\beta A^*A\,A^* + A\,A^{*2}\rbrack 
&=&\gamma^* \lbrack A^*, ~A^*A+A\,A^*\rbrack  
+ \delta^* \lbrack A^*, ~A\rbrack, 
\end{eqnarray*}
where $\lbrack X,Y\rbrack $ means $X\,Y-Y\,X$. 
The sequence of scalars  $ \beta,\, \gamma,\, \gamma^*, \delta,\, \delta^*$ 
is unique if $d \geq 3$.  
The above relations
are called the \textit{tridiagonal relations} (\textit{TD-relations}) [8]. 
We fix a nonzero $q \in \C$ such that 
\[\beta = q^2 + q^{-2}.\] 
Let $\theta_i$ (resp. $\theta^*_i$ ) denote the eigenvalue of 
\textit{A} for $V_i$ (resp. $A^*$ for $V_i^*$) ($0 \leq i \leq d$). 
They are expressed as follows \cite[Theorem 11.2]{ITT}.

\begin{quote}
{\bf Type I}  
($q^2 \neq \pm 1$): there exist scalars $a,\, a^*, b,\, b^*, c,\,c^*$ 
such that
\begin{eqnarray*}
\theta_i &=& a + b\,q^{2i} + c\,q^{-2i}  
\qquad ~~\,(0 \leq i \leq d),\\
\theta_i^* &=& a + b^* q^{2i} + c^* q^{-2i} 
\qquad ~(0 \leq i \leq d).
\end{eqnarray*}
In this case, $\gamma = -{(q- q^{-1})}^2 a$, 
$\gamma^* = -{(q- q^{-1})}^2 a^*$, 
$\delta = {(q - q^{-1})}^2 a^2 - {(q^2 -q^{-2})}^2 b\,c$, 
$\delta^* = {(q - q^{-1})}^2 {a^*}^2 - {(q^2 -q^{-2})}^2 b^* c^*$.
\\

\noindent
{\bf Type II}  
$(q^2=1)$ : there exist scalars $a, a^*, b, b^*, c,\,c^*$ such that 
\begin{eqnarray*}
\theta_i &=& a + b\,i + c\,i^2  
\qquad ~~~\:(0 \leq i \leq d),\\
\theta_i^*  &=& a^* + b^*i + c^*i^2
\qquad ~(0 \leq i \leq d).
\end{eqnarray*}
In this case, 
$\gamma=2\,c, ~\gamma^*=2\,c^*, ~\delta=b^2-c^2-4\,ac$, 
$\delta^* = {b^*}^2 -{c^*}^2-4\,a^*c^*.$ 
\\

\noindent
{\bf Type III}  
$(q^2=-1)$ : there exist scalars  $a, a^*, b, b^*, c,\,c^*$ such that 
\begin{eqnarray*}
\theta_i &=& a + b\,{(-1)}^i + c\,{(-1)}^i\,i 
~~\:\qquad (0 \leq i \leq d),\\
\theta_i^* &=& 
a^* + b^*{(-1)}^i + c^* {(-1)}^i\,i \qquad (0 \leq i \leq d).
\end{eqnarray*}
In this case, 
$\gamma=4\,a,~ \gamma^*=4\,a^*,~ \delta=-4\,a^2 +c^2,~ 
\delta^*=-4\,{a^*}^2 + {c^*}^2$ . 
\end{quote}

\medskip
In this paper, we treat TD-pairs of Type I. 
If a TD-pair of Type I comes from a P- and Q-polynomial association scheme 
with sufficiently large diameter, 
then \textit{q} is not a root of unity, i.e., $q^n \neq 1$ for any 
nonzero integer $n$ \cite[Proposition 7.7]{BI}. 
\textit{From now on, we fix a nonzero scalar $q \in \C$ and assume that 
$q$ is not a root of unity.} 
One of the effects of this assumption is as follows. 
Let us call the conditions (ii), (iii) for a TD-pair 
the \textit{TD-structures}.
Then under the diagonalizability condition (i) and the irreducibility 
condition (iv), the TD-relations imply the TD-structures 
\cite[Theorem 3.10]{Tb}. 
This allows us to work with the TD-relations 
instead of the TD-structures. 
We first establish the representation theory of 
the augmented tridiagonal algebra $\mathcal{T}_q$. 
The classification of TD-pairs of Type I 
will be given as an application of the representation theory.

\medskip
If $A, A^*$ are a TD-pair on $V$, then 
$\lambda A + \mu I$, $\lambda^* A^* + \mu^* I$ are also a TD-pair on $V$ 
with the same eigenspaces. Here  
$\lambda, \lambda^*, \mu, \mu^* \in \C$, 
$\lambda \neq 0, ~\lambda^* \neq 0$ 
and $I$ is the identity map.  
The parameter $\beta$ and hence \textit{q} are 
invariant under the affine transformations 
$A~\mapsto~\lambda A + \mu I,~
A^* \mapsto \lambda^* A^* + \mu^* I$.    
Also the diameter $d$ is invariant under the affine transformations. 
For fixed $d$ and $q$, consider a TD-pair $A, A^*$ of Type I 
with diameter \textit{d} and 
the parameter $\beta=q^2+q^{-2}$. The TD-pair $A, A^*$ 
can be standardized to have the following eigenvalues 
by applying appropriate affine transformations and, if necessary, reversing 
the ordering of the eigenspaces $V_i$  of $A$ or 
of the eigenspaces $V_i^{*}$ of $A^*$: 
\begin{eqnarray}
\label{theta}
\theta_i &=& b\,q^{2i-d} + \varepsilon\,b^{-1} q^{d-2i} 
\qquad ~~\,(0 \leq i \leq d),
\\
\label{theta*}
\theta_i^* &=& \varepsilon^* b^* q^{2i-d}+ {b^*}^{-1} q^{d-2i} 
\qquad (0 \leq i \leq d) 
\end{eqnarray}
for some constants $b, b^*$ ($b \neq 0, ~b^* \neq 0$) and 
$\varepsilon, \varepsilon^* \in$ \{1,~0\}. 
A TD-pair $A, A^*$ is called a 
\textit{standardized} TD-pair of Type I, if 
$A$, $A^*$ have 
eigenvalues $\{\theta_i\}_{i=0}^d$, 
$\{\theta_i^*\}_{i=0}^d$ 
as in (\ref{theta}), (\ref{theta*}) respectively 
for some integer $d \geq 1$ and nonzero $b$, $b^*$ $\in \C$ 
under suitable orderings of the eigenspaces $\{V_i\}_{i=0}^d$, 
$\{V_i^*\}_{i=0}^d$. 

\medskip
If $d=1$, then $\theta_0$, $\theta_1$ (resp. $\theta_0^*$, $\theta_1^*$) 
can be any pair of distinct scalars by applying a suitable affine 
transformation to $A$ (resp. $A^*$), in particular for stadardization, 
$(\varepsilon, b)$ and $(\varepsilon^*, b^*)$ can be chosen arbitrarily 
from $\{0, 1\} \times \C^{\times}~ \backslash ~\{(1, \pm 1)\}$, 
where $\C^{\times} = \C \backslash \{0\}$. 
Assume $d \geq 2$. Then the pair 
$\varepsilon, \varepsilon^*$ 
is uniquely determined by $A, A^*$ regardless of 
standardization, but the scalars $b, b^*$ are not. 
If $\varepsilon = 1$ (resp. $\varepsilon^* = 1$), then 
$b$ (resp. $b^*$) is determined up to the $\pm$ sign 
by $A$ (resp. $A^*$) and by 
the ordering of the eigenspaces of $A$ (resp. $A^*$). 
In this case,  
$b$ (resp. $b^*$) is changed to $b^{-1}$ (resp. ${b^*}^{-1}$) 
when we reverse the ordering of the eigenspaces of $A$ (resp. $A^*$).  
If $\varepsilon = 0$ (resp. $\varepsilon^* = 0$), then 
$b$ (resp. $b^*$) can be an arbitrary nonzero scalar.  
In this case, the ordering of the eigenspaces of $A$ (resp. $A^*$) 
is uniquely determined when standardized. 

\medskip
If ($\varepsilon, \varepsilon^*) = (0,1)$, we further standardize 
the TD-pair $A, A^*$ 
by interchanging $A, A^*$ and then reversing the ordering of the 
eigenspaces $V_i^*$ so that the standardized TD-pair has 
($\varepsilon, \varepsilon^*) = (1,0)$. Thus a standardized TD-pair 
of Type~I has 
\begin{eqnarray*}
(\varepsilon, \varepsilon^*) = (1,1),~(1,0)~ \mbox{or}~ (0,0)
\end{eqnarray*}
and is called of the \textit{1st, 2nd, 3rd kind}, accordingly. 

\medskip
The TD-relations for a standardized TD-pair $A, A^*$  of Type I are 

\[
({\rm TD})\left\{
\begin{array}{lll}
[A,~A^2A^*-\beta AA^*A+A^*A^2] 
&=& \varepsilon \,\delta [A,~A^*],
\\
\lbrack A^*,~A^{*^2} A -\beta A^*AA^* +A A^{*^2} \rbrack
&=&\varepsilon^* \delta [A^*,~A], 
\end{array}
\right.
\]
where $\beta = q^2 + q^{-2}$ and $\delta = -{(q^2 - q^{-2})}^2$. 
Conversely, 
if a TD-pair $A, A^*$ satisfies the above TD-relations (TD), 
then we have $a=a^*=0$, $b\,c=\varepsilon$, $b^*c^*=\varepsilon^*$ 
in the general expression of the eigenvalues for Type I, and so 
with suitable orderings of the eigenspaces, 
$A, A^*$ have eigenvalues in the form of (\ref{theta}), (\ref{theta*}) 
for some integer $d \geq 1$ and some nonzero $b$, $b^*$, 
i.e., $A, A^*$ are a standardized TD-pair of Type I. 
Thus 
given 
$(\varepsilon, \varepsilon^*) \in \{(1,1), (1,0), (0,0)\}$ 
and a nonzero scalar $q$ that is not a root of unity, 
a TD-pair $A, A^*$ is a standardized TD-pair if and only if 
it satisfies the above TD-relations (TD).  
In view of this fact, we call (TD) the 
\textit {stadardized} TD-relations of Type I.

\medskip
Given TD-pair $A, A^* \in {\rm{End}}(V)$ with eigenspaces 
${\{V_i\}}^d_{i=0},~{\{V_i^*\}}^d_{i=0}$, 
the underlying vector space $V$ has the \textit{split decomposition} 
\cite[Theorem 4.6]{ITT}:

\[ V= \bigoplus^d_{i=0} U_i , \]
where

\[ U_i=(V_0^* + \cdots + V_i^*) \cap (V_i + \cdots + V_d). \]
For a TD-pair $A, A^* \in {\rm{End}}(V)$ of Type I 
with eigenvalues (\ref{theta}), (\ref{theta*}), 
let $K \in {\rm{End}}(V)$ denote the diagonalizable transformation for which 
$U_i$ is the eigenspace belonging to the eigenvalue 
$q^{2i-d} ~(0 \leq i \leq d)$. We define the  \textit{raising map} $R$ 
and the  \textit{lowering map} $L$ by 
\begin{eqnarray*}
R&=&A-b\,K - \varepsilon\, b^{-1} K^{-1},\\
L&=&A^* - \varepsilon^* b^* K - {b^*}^{-1} K^{-1} . 
\end{eqnarray*} 
Then indeed $R$ (resp. $L$ ) has the \textit{raising} 
(resp. \textit{lowering}) \textit{property} 
\cite[Theorem~4.6, Corollary~6.3]{ITT}: 
\begin{eqnarray*} 
RU_i&\subseteq& U_{i+1} \qquad (0 \leq i \leq d), \\
LU_i&\subseteq& U_{i-1} \qquad (0 \leq i \leq d),  
\end{eqnarray*} 
where $U_{-1}=U_{d+1}=0$. By the raising, lowering properties of $R$, $L$,
we get
\[
\hspace{-5.5cm}
({\rm TD})_0'\left\{
\begin{array}{lll}
KRK^{-1}&=&q^2 R,\\
KLK^{-1}&=&q^{-2} L,
\end{array} 
\right.
\]
and conversely the relations $(TD)'_0$ imply the raising, lowering properties  
of $R$, $L$. Writing $(TD)'_0$ in terms of $A, A^*, K$,  we get the 
\textit{generalized q-Weyl relations}: 
\[
({\rm TD})_0 \left\{
\begin{array}{lll}
(qAK-q^{-1} KA)/(q - q^{-1})&=&bK^2 + \varepsilon b^{-1} I,\\
(qKA^* -q^{-1} A^* K)/(q - q^{-1})&=&\varepsilon^* b^* K^2 + {b^*}^{-1} I,
\end{array} 
\right.
\]
where $I$ is the identity map. Writing the tridiagonal relations 
(TD) for $A, A^*$ in terms of $R$, $L$, $K$, we get 

\[
({\rm TD})'\left\{
\begin{array}{lll}
\lbrack R, ~R^2 L - \beta RLR + LR^2 \rbrack
&=&\delta'(\varepsilon^* s^2 R^2 K^2 
- \varepsilon s^{-2} K^{-2} R^2),\\
\lbrack L, ~L^2 R - \beta LRL + R L^2 \rbrack 
&=&\delta'(-\varepsilon^* s^2 K^2 L^2 
+ \varepsilon s^{-2} L^2 K^{-2}),
\end{array} 
\right.
\]
where 
$\beta = q^2 + q^{-2}$, $\delta' = -(q-q^{-1})(q^2-q^{-2})(q^3-q^{-3})q^4$, 
$s^2 = b\,b^*$.

\subsection{The TD-algebra $\A$ and the augmented TD-algebra $\T$}
Fix a nonzero scalar $q \in \C$ 
which is not a root of unity. We also fix 
$(\varepsilon, \varepsilon^*) \in \{(1,1), (1,0), (0,0)\}$.
Let $\A=\A_q^{(\varepsilon, \varepsilon^*)}$ 
denote the associative $\C$-algebra with $1$ defined by 
genarators $z, z^*$ subject to the relations 
\[({\rm TD})\left\{
\begin{array}{lll}
[z,~z^2z^*-\beta z\,z^*z+z^*z^2] 
&=& \varepsilon \,\delta \,[z,~z^*],
\\
\lbrack z^*,~{z^*}^2 z -\beta z^*z\,z^* +z \,{z^*}^2 \rbrack
&=&\varepsilon^* \delta \,[z^*,~z], 
\end{array}
\right.
\]
where $\beta = q^2 + q^{-2}$ and $\delta = -{(q^2 - q^{-2})}^2$. 
When we need to specify $(\varepsilon,\varepsilon^*)$, we write 
$({\rm TD})_{\rm I}$, $({\rm TD})_{\rm II}$, $({\rm TD})_{\rm III}$ 
for the relations (TD) and 
$\A_{\rm I}, \:\A_{\rm II}, \:\A_{\rm III}$ 
for the algebra $\A$ according to $(\varepsilon,\varepsilon^*)
=(1,1),(1,0),(0,0)$. The algebra $\A$ is called the 
\textit{tridiagonal algebra} (\textit{TD-algebra}) \cite{Tb} 
of the 1st, 2nd, 3rd kind, 
accordingly. 
$({\rm TD})_{\rm III}$ is the $q$-Serre relations and $\A_{\rm III}$ is 
isomorphic 
to the positive part of the quantum affine algebra $U_q({\widehat{sl}}_2)$. 
$({\rm TD})_{\rm I}$ can be regarded as a $q$-analogue of 
the Dolan- Grady relations 
and we call $\A_{\rm I}$ the \textit{q-Onsager algebra}. 

\medskip
Let $\T=\T_q^{(\varepsilon, \varepsilon^*)}$ denote the associative 
$\C$-algebra with $1$ defined by generators 
$x, y, k, k^{-1}$ subject to the relations 

\[
\hspace{-5.2cm}({\rm TD})_0'\left\{
\begin{array}{l}
k\,k^{-1}=k^{-1}k=1,\\
k\,x\,k^{-1}=q^2x,\\
k\,y\,k^{-1}=q^{-2}y, 
\end{array} 
\right.
\]
and
\[
({\rm TD})'\left\{
\begin{array}{l}
\lbrack
x,~x^2 y -\beta\, x\,y\,x+y\,x^2 
\rbrack
= \delta'(\varepsilon^* x^2 k^2 - \varepsilon\, k^{-2} x^2),\\
\lbrack
y,~y^2 x- \beta\, y\,x\,y + x\,y^2
\rbrack
=\delta' (- \varepsilon^* k^2 y^2 + \varepsilon\, y^2 k^{-2}),
\end{array} 
\right.
\]
where $\beta=q^2 + q^{-2}$, 
$\delta'=-(q-q^{-1})(q^2 - q^{-2})(q^3 - q^{-3})q^4$. 
When we need to specify $(\varepsilon, \varepsilon^*)$, 
we write $({\rm TD})_{\rm I}'$, 
$({\rm TD})_{\rm II}'$, 
$({\rm TD})_{\rm III}'$ for the relations $({\rm TD})'$ 
and $\T_{\rm I}, \T_{\rm II}, \T_{\rm III}$ for 
the algebra $\T$ according to 
$(\varepsilon, \varepsilon^*)=(1,1), (1,0), (0,0)$. 
The algebra $\T$ is called the \textit{augmented tridiagonal algebra} 
(\textit{augmented TD-algebra}) of the 1st, 2nd, 3rd, kind, accordingly.  
$\T_{\rm III}$ is isomorphic to the Borel subalgebra of 
the quantum affine algebra 
$U_q({\widehat{sl}}_2)$.

\medskip
The augmented TD-algebra $\T$ has another presentation. 
Fix a nonzero scalar $t \in \C$. 
Define the elements $z_t, z^*_t \in \T$ to be 
\begin{eqnarray}
\label{zt}
z_t=x+t\,k+\varepsilon\, t^{-1} k^{-1},\\
\label{z*t}
z^*_t=y+\varepsilon^* t^{-1} k + t\, k^{-1}. 
\end{eqnarray}
Then $\T$ is generated by $z_t, z^*_t, k, k^{-1}$  and the 
following relations hold : 

\[
({\rm TD})_0 \left\{
\begin{array}{lll}
k\, k^{-1} = k^{-1} k = 1 & & \\
(q \,z_t k - q^{-1} k\, z_t)/ (q - q^{-1}) &=&t\, k^2 + \varepsilon \,t^{-1},\\
(q\, k\, z_t^*  - q^{-1} z_t^* k)/ (q - q^{-1}) 
&=& \varepsilon^* t^{-1} \,k^2 + t,
\end{array} 
\right.
\]
and
\[
\hspace{-2mm}
({\rm TD}) \left\{
\begin{array}{lll}
\lbrack
z_t,~z^2_t z^*_t - \beta z_t z^*_t z_t + z_t^* z_t^2
\rbrack
&=& \varepsilon \, \delta\, 
\lbrack
z_t,~z_t^*
\rbrack,\\
\lbrack
z^*_t, {z^*_t}^2 z_t - \beta z^*_t z_t z^*_t + z_t {z^*_t}^2 
\rbrack
&=&  
\varepsilon^* \delta\, 
\lbrack
z^*_t, ~z_t
\rbrack,
\end{array} 
\right.
\]
where $\beta=q^2+q^{-2}$, $\delta=-{(q^2 - q^{-2})}^2$. 
One routinely verifies that $\T$ is isomorphic to the algebra 
generated by symbols $z_t, \,z^*_t, \,k, \,k^{-1}$ with $(TD)_0, \;(TD)$ 
the defining relations.

\begin{proposition}
\label{prop: iota}
There exists an algebra homomorphism $\iota_t$ from $\A$ to $\T$ that 
sends $z,\,z^*$ to $z_t,\,z_t^*$, respectively :
$$ \iota_t : \A \longrightarrow \T 
\qquad ( z, \,z^* \mapsto z_t, \,z^*_t) .$$
Moreover $\iota_t$ is injective.
\end{proposition}

It is obvious that the correspondence $z, \,z^* \mapsto z_t, \,z^*_t$ 
can be extended to 
an algebra homomorphism from $\A$ to $\T$. 
The injectivity of $\iota_t$ will be proved in Section 2.

\begin{lemma}
\label{lemma: weight-space decomposition}
Let 
$\rho: \T \longrightarrow {\rm{End}}(V)$ be a finite-dimensional 
irreducible representation of $\T$. 
Then $\rho (k)$ is diagonalizable 
with eigenvalues $\{s \,q^{2i-d}~| ~0 \leq i \leq d \}$ for some 
nonzero $s \in \C$ and an integer $d \geq 0$. 
Let $V = \bigoplus_{i=0}^d U_i$ 
denote the eigenspace decomposition of $\rho (k)$, 
where $U_i$ is the eigenspace belonging to $s \,q^{2i-d}$. 
Then regarding $V$ as an irreducible $\T$-module via $\rho$, 
we have 
$$ x \,U_i \subseteq U_{i+1},~~ 
y\, U_i \subseteq U_{i-1} ~~~~ (0 \leq i \leq d), $$
where $U_{-1}=U_{d+1}=0$. 
In particular $\rho(x),\,\rho(y)$  are nilpotent. 
\end{lemma}

The scalar $s$ (resp. the integer $d$) is 
called the \textit{type} (resp. \textit{diameter}) 
of the representation $\rho$ 
and the $\T$-module $V$. 
We call the direct sum $V=\bigoplus_{i=0}^d U_i$ 
the \textit{weight-space decomposition} and 
$U_0$ the \textit{highest weight space}.

\medskip
\noindent
Proof. For $\theta \in \C$, 
set $U(\theta) = \{v \in V ~|~ k v = \theta v \}$.
Note that $\theta$ is an eigenvalue of $\rho(k)$ if and only if 
$U(\theta) \neq 0 $ , and in this case $U(\theta)$ is the corresponding 
eigenspace. Using the relations $kx=q^2 xk$ and $ky=q^{-2} yk$ , 
we find $x\,U(\theta) \subseteq U(q^2 \theta)$ and 
$y\, U(\theta) \subseteq U(q^{-2} \theta)$.
Now assume that $\theta$ is an eigenvalue of $\rho(k)$. 
Observe that $\theta \neq 0$ since $k^{-1}$ exists, and that 
$\sum_{i \in \Z} U(q^{2i} \theta)$ is invariant 
under each of $x,\,y, \,k^{\pm 1}$ 
and the sum is a finite sum by dim $V < \infty$.
These elements $x,\,y,\,k^{\pm}$ generate $\T$ and the $\T$-module $V$ 
is irreducible, so we have $V =\sum_{i \in \Z} U(q^{2i} \theta)$.
This yields the lemma. 
\hfill $\Box $

\begin{proposition}
\label{prop: C1}
Let $\rho: \T \longrightarrow {\rm{End}}(V)$ be a finite-dimensional 
irreducible representation of $\T$ with type $s$, diameter $d$, 
and $V = \bigoplus_{i=0}^d U_i$ the weight-space decomposition. 
Let $z_t$, $z^*_t$ be as in (\ref{zt}), (\ref{z*t}). 
\begin{enumerate}
\item[$(i)$]
$\rho(z_t)$ is diagonalizable if and only if the scalars 
$$ \theta_i = st q^{2i-d} + \varepsilon s^{-1} t^{-1} q^{d-2i} 
~~~~~ (0 \leq i \leq d)$$
are mutually distinct. In this case, $\{\theta_i \}_{i=0}^d$ is the set of 
eigenvalues of $\rho(z_t)$ and it holds that 
$$V_i+V_{i+1}+\cdots+V_d=U_i+U_{i+1}+\cdots+U_d~~~(0 \leq i \leq d),$$
where $V_i$ is the eigenspace of $\rho(z_t)$ belonging to $\theta_i$. 
\item[$(ii)$]
$\rho(z^*_t)$ is diagonalizable 
if and only if the scalars 
$$ \theta^*_i = \varepsilon^*s t^{-1} q^{2i-d} + s^{-1} t q^{d-2i} 
~~~~~ (0 \leq i \leq d)$$ 
are mutually distinct. In this case, $\{\theta^*_i \}_{i=0}^d$ is the set of 
eigenvalues of $\rho(z^*_t)$ 
and it holds that 
$$V_0^*+V_1^*+\cdots+V_i^*=U_0+U_1+\cdots+U_i~~~(0 \leq i \leq d),$$
where $V_i^*$ is the eigenspace of $\rho(z_t^*)$ belonging to $\theta_i^*$.
\end{enumerate}
\end{proposition}

Proposition \ref{prop: C1} will be proved in Section 2. 

\medskip
Recall we are given in advance 
$(\varepsilon, \varepsilon^*) \in \{(1,1), (1,0), (0,0)\}$ 
and a nonzero scalar $q$ that is not a root of unity. 
Let $\rho: \A \longrightarrow {\rm{End}}(V)$ be a finite-dimensional 
irreducible representation of the TD-algebra 
$\A=\A_q^{(\varepsilon, \varepsilon^*)}$. 
We assume that $\rho$ 
satisfies the following property $({\rm C}_1)$:  
\begin{quote}
$({\rm C}_1)$: ~$\rho(z), \,\rho(z^*)$ are both diagonalizable.  
\end{quote} 
Set $A=\rho(z),\, A^*=\rho(z^*)$. 
Then $A,\, A^*$ satisfy the TD-relations. 
The TD-relations for $A, \,A^*$ imply the TD-structures, 
i.e., the conditions (ii), (iii) for a TD-pair hold for $A,\,A^*$,  
while the conditions (i), (iv) hold for $A,\,A^*$ 
by the property $({\rm C}_1)$ and the irreducibility of $\rho$.  
So $A, \,A^* \in {\rm{End}}(V)$ are a TD-pair on $V$. 
Moreover since the TD-relations (TD) 
for $A,\, A^*$ is the standardized TD-relations of Type I,  
the TD-pair $A,\, A^*$ is a standardized TD-pair of Type~I on $V$. 

\medskip
Conversely, 
let us start with a standardized TD-pair $A,\, A^*$ of Type I on $V$, 
where we understand $q$ and $(\varepsilon,\varepsilon^*)$ 
are chosen in advance and fixed.  
Consider the TD-algebra $\A=\A_q^{(\varepsilon,\varepsilon^*)}$.
Then by the TD-relations (TD) for $A,\, A^*$, 
we obtain a finite-dimensional 
representation $\rho$ of $\A$ that sends $z, \,z^*$ to 
$A,\, A^*$,  respectively: 
$$\rho: \A \longrightarrow {\rm{End}}(V) ~~~
(z, z^* \mapsto A, A^*). $$ 
This representation $\rho$ is irreducible and 
satisfies the property $({\rm C}_1)$ by the conditions (iv), (i)  
for the TD-pair $A,\, A^*$. 

\medskip
We restate what we saw in the previous two paragraphs 
as a proposition below. 
We are given in advance 
$(\varepsilon, \varepsilon^*) \in \{(1,1), (1,0), (0,0)\}$ 
and 
a nonzero scalar $q$ that is not a root of unity. 
Let $\mathcal {STD}$ 
denote the set of isomorphism classes of 
standardized TD-pairs $A,\,A^*$ of Type~I
together with the trivial TD-pairs:  $A$ (resp. $A^*$) has 
eigenvalues $\{\theta_i\}_{i=0}^d$ 
(resp. $\{\theta_i^*\}_{i=0}^d$) 
as in (\ref{theta}) (resp. (\ref{theta*})) 
for some integer $d \geq 0$ and nonzero $b$ (resp. $b^*$) $\in \C$ 
with a suitable ordering of the eigenspaces $\{V_i\}_{i=0}^d$ 
(resp. $\{V_i^*\}_{i=0}^d$). 
Set $\A=\A_q^{(\varepsilon, \varepsilon^*)}$. 
Let $\mathcal {I}rr(\A)$ denote the set of isomorphism classes of 
finite-dimensional irreducible representations of $\A$ 
that satisfy the property $({\rm C}_1)$. 
Then we have 
\begin{proposition}
\label{prop: STD} The mapping 
$\rho \mapsto A=\rho(z),\, A^*=\rho(z^*)$ gives a bijection 
from  $\mathcal {I}rr(\A)$ to $\mathcal {STD}$. 
The trivial representations, i.e., 1-dimensional representations, 
correspond to the trivial TD-pairs. 
\end{proposition}

Thus the classification of standardized TD-pairs of Type~I is 
reduced to the following problem.

\begin{problem}
\label{problem: A}
Classify up to isomorphism 
the finite-dimensional irreducible representations of $\A$ 
that satisfy the property $({\rm C}_1)$. 
\end{problem}

Let us start with a finite-dimensional irreducible 
representation $\rho: \T \longrightarrow {\rm{End}}(V)$ 
of the augmented TD-algebra $\T$ with type $s$ and 
diameter $d$. 
We assume that $\rho$ satisfies the following properties 
${({\rm C}_1)}_t , \,{({\rm C}_2)}_t$
for some nonzero $t \in \C$: 

\begin{quote}
${({\rm C}_1)}_t$:  ~~$\rho(z_t), \,\rho(z^*_t) $ are both diagonalizable.\\
${({\rm C}_2)}_t$:  ~
\begin{minipage}[t]{13cm}
The restriction $\rho|_{\langle z_t, z^*_t\rangle}: 
\langle z_t, z^*_t\rangle \longrightarrow {\rm{End}}(V)$ is irreducible, \\
where $\langle z_t,z^*_t\rangle$
is the subalgebra of $\T$ genarated by $z_t, \,z^*_t$.
\end{minipage}
\end{quote}

\noindent
Set $A=\rho(z_t), \,A^*=\rho(z^*_t)$.
Then $A, \,A^*$ satisfy the TD-relations. 
Since  
the TD-relations for $A,\, A^*$ imply the TD-structures for $A,\, A^*$,  
we find $A,\, A^*$ are a TD-pair on  $V$. 
By Proposition~\ref{prop: C1}, the TD-pair $A,\,A^*$ has 
distinct eigenvalues 
$\{ \theta_i\}_{i=0}^d, \,\{\theta_i^*\}_{i=0}^d$ as 
in (\ref{theta}), (\ref{theta*}) 
with $b=st, \,b^*=st^{-1}$. 
So $A, \,A^* \in {\rm{End}}(V)$ are a standardized TD-pair of Type I. 
By Lemma \ref{lemma: weight-space decomposition} and 
Proposition~\ref{prop: C1},
the eigenspace decomposition for $\rho(k)$ 
coincides with the split decomposition for the TD-pair $A, \,A^*$. 
So we have $\rho(x)=R, \,\rho(y)=L, \,\rho(k)=sK$, 
where $R, \,L, \,K$ are the raising, lowering, diagonalizable maps, 
respectively, associated with the split decomposition.

\medskip
Conversely, let us start with 
a standardized TD-pair $A,\,A^* \in {\rm{End}}(V) $ of Type I with 
eigenvalues 
\begin{eqnarray*}
\theta_i &=& b\,q^{2i-d} + \varepsilon\,b^{-1} q^{d-2i} 
\qquad  ~~\:(0 \leq i \leq d),\\
\theta_i^* &=& \varepsilon^* b^* q^{2i-d}+ {b^*}^{-1} q^{d-2i} 
\qquad  (0 \leq i \leq d), 
\end{eqnarray*}
respectively as in (\ref{theta}), (\ref{theta*}), 
where we understand $q$ and $(\varepsilon,\varepsilon^*)$ 
are chosen in advance and fixed. 
We have the raising map $R$, the lowering map $L$ and the 
diagonalizable $K$ associated with the split decomposition for 
the TD-pair $A, \,A^*$. 
Consider the augmented TD-algebra $\T=\T_q^{(\varepsilon,\varepsilon^*)}$. 
Define the nonzero scalars \textit{s,~t} $\in \C$ by 
\begin{eqnarray}
\label{st}
b= st , \qquad b^* = s t^{-1}. 
\end{eqnarray}
The scalars \textit{s, t} are determined by $b,~ b^*$ up to the $\pm$ sign : 
$s^2 = b b^*,~ t^2 = b {b^*}^{-1}$. We choose \textit{s, t} as one of the 
solutions of (\ref{st}) and fix them. 
By the relations $({\rm TD})_0', \,({\rm TD})'$ for 
$R, \,L, \,K,$ we obtain a finite-dimensional representation 
$\rho$ of $\T$ with type $s$ and diameter $d$ 
that sends $x, \,y, \,k$ to $R, \,L, \,sK$, respectively: 
$$ \rho: \T \longrightarrow {\rm{End}}(V) ~~~~~ 
(x, \,y, \,k \mapsto R, \,L, \,sK) $$ 
expressed 
in terms of the 1st presentation of $\T$ with respect to the generators 
$x, \,y, \,k, \,k^{-1}$. 
By (\ref{zt}), (\ref{z*t}), it holds that 
$\rho(z_t)=A, \,\rho(z^*_t)=A^*$. 
So we have 
$$ \rho: \T \longrightarrow {\rm{End}}(V) ~~~~~ 
(z_t, \,z^*_t, \,k \mapsto A, \,A^*, \,sK) $$ 
expressed 
in terms of the 2nd presentation of $\T$ with respect to the generators 
$z_t, \,z^*_t, \,k, \,k^{-1}$. 
By the conditions (iv), (i)  for the TD-pair $A,\,A^*$, $\rho$ 
is irreducible and satisfies 
the properties $ ({\rm C}_1)_t $, $({\rm C}_2)_t$.

\medskip
We restate what we saw in the previous two paragraphs 
as a proposition below. 
We are given in advance 
$(\varepsilon, \varepsilon^*) \in \{(1,1), (1,0), (0,0)\}$ 
and 
a nonzero scalar $q$ that is not a root of unity. 
Suppose that we are further given a positive integer $d$ 
and nonzero $b, \,b^* \in \C$ such that 
the scalars 
$\theta_i = b q^{2i-d} + \varepsilon b^{-1} q^{d-2i}$  ($0 \leq i \leq d$) 
in (\ref{theta}) 
are mutually distinct 
and the scalars 
$\theta^*_i=\varepsilon^* b^* q^{2i-d} + {b^*}^{-1} q^{d-2i}$ 
($0 \leq i \leq d$) 
in (\ref{theta*}) 
are mutually distinct.  
By   
$\mathcal {STD}_d^{(b,b^*)}$ 
we denote 
the set of isomorphism classes of standardized 
TD-pairs $A,\,A^*$ of Type~I with eigenvalues 
$\{\theta_i \}^d_{i=0},\, \{ \theta_i^* \}^d_{i=0}$ respectively. 
Note that if a standardized TD-pair $A,\,A^*$ of Type~I belongs 
to $\mathcal {STD}_d^{(b,b^*)}$, then the ordering of the eigenspaces 
$\{V_i\}_{i=0}^d$ of $A$ 
(resp. $\{V_i^*\}_{i=0}^d$ of $A^*$) is uniquely determined 
by the corresponding eigenvalues 
$\theta_i = b q^{2i-d} + \varepsilon b^{-1} q^{d-2i}$ 
(resp. $\theta^*_i=\varepsilon^* b^* q^{2i-d} + {b^*}^{-1} q^{d-2i}$). 
Recall that if $\varepsilon = 1$ (resp. $\varepsilon^* = 1$), then 
$b$ (resp. $b^*$) is changed to $b^{-1}$ (resp. ${b^*}^{-1}$) 
when we reverse the ordering of the eigenspaces of $A$ (resp. $A^*$).  
Thus if $\varepsilon=1$ (resp. $\varepsilon^*=1$), then 
$\mathcal {STD}_d^{(b,b^*)}=\mathcal{STD}_d^{(b^{-1},b^*)}$ 
(resp. $\mathcal {STD}_d^{(b,b^*)}=\mathcal{STD}_d^{(b,{b^*}^{-1})}$): 
$\mathcal{STD}_d^{(b^{-1},b^*)}$ (resp. $\mathcal{STD}_d^{(b,{b^*}^{-1})}$) 
coincides with $\mathcal {STD}_d^{(b,b^*)}$ as sets of isomorphism classes 
of standardized TD-pairs $A,\,A^*$ of Type I but has the ordering of 
the eigenspaces of $A$ (resp. $A^*$) reversed. 
Set $b=st, ~b^*=st^{-1}$ as in (\ref{st}).  
Such scalars $s, \,t$ are determined by $b, \,b^*$ 
uniquely up to the $\pm$ sign. 
We choose one of them and fix it. 
Note that if $(s,t)$ is a solution of $b=st, ~b^*=st^{-1}$, 
then 
\begin{eqnarray}
\label{(s',t')}
(s',t')=(t^{-1},s^{-1}), \,(t,s), \,(s^{-1},t^{-1}) 
\end{eqnarray}
are a solution of $c=s't', ~c^*=s't'^{-1}$ for 
\begin{eqnarray}
\label{(c,c*)}
(c,c^*)=(b^{-1},b^*),\,(b,{b^*}^{-1}),\,(b^{-1},{b^*}^{-1})
\end{eqnarray}
respectively. 
Set $\T=\T_q^{(\varepsilon, \varepsilon^*)}$. 
By $\mathcal {I}rr_d^{s, t}(\T)$ 
we denote 
the set of isomorphism classes of 
finite-dimensional irreducible representations $\rho$ of $\T$ 
with type $s$ and diameter $d$ 
that satisfy the properties 
$({\rm C}_1)_t, \,({\rm C}_2)_t$ for the scalar $t$. 
Then we have 

\begin{proposition}
\label{prop: STDd(b,b*)} 
The mapping $\rho \mapsto A=\rho(z_t), \,A^*=\rho(z_t^*)$ 
gives a bijection 
from $\mathcal {I}rr_d^{s, t}(\T)$ to $\mathcal {STD}_d^{(b,b^*)}$, 
where $b=st, \,b^*=st^{-1}$. 
\end{proposition}
  
Thus Problem \ref{problem: A} 
is reduced to the following problem. 

\begin{problem}
\label{problem: T}
~~~\\
\vspace{-7mm}
\begin{enumerate}
\item[$(i)$] 
Classify up to isomorphism 
the finite-dimensional irreducible representations of $\T$ 
with type $s$ and diameter $d$. 
\item[$(ii)$]
Determine when a finite-dimensional irreducible representation 
$\rho$ of $\T$ with type $s$ and diameter $d$ satisfies 
the properties $({\rm C}_1)_t,\,({\rm C}_2)_t$ .
\end{enumerate}
\end{problem}

We solve Problem~\ref{problem: T} in this paper. 
Problem~\ref{problem: A} is reduced to 
Problem~\ref{problem: T} via $\mathcal {STD}_d^{(b,b^*)}$ as follows.  
The set $\mathcal {STD}$ is the disjoint union of the trivial TD-pairs and 
$\mathcal {STD}_d^{(b,b^*)}$ over $d \in \N$ and 
$(b, b^*) \in (\C \backslash \{0\}) \times (\C \backslash \{0\})/\sim$, 
where $\sim$ is the equivalence relation defined by 
$(b,b^*)\sim (c,c^*)$ if and only if 
\begin{eqnarray}
\label{sim 1st kind}
(c,c^*)&\in& 
\{(b, b^*),\,(b^{-1}, b^*),\,(b, {b^*}^{-1}),\,(b^{-1}, {b^*}^{-1})\}~~~
{\rm for~the~case} ~(\varepsilon, \varepsilon^*)=(1,1), \\
\label{sim 2nd kind}
(c,c^*)&\in& 
\{(b, b^*),\,(b^{-1}, b^*)\}~~~
{\rm for~the~case} ~(\varepsilon, \varepsilon^*)=(1,0), 
\end{eqnarray} 
and $(b, b^*)=(c,c^*)$ for the case $(\varepsilon, \varepsilon^*)=(0,0)$. 
For nonzero $b,\,b^*\in \C$, 
let $\mathcal {I}rr_d^{(b,b^*)}(\A)$ denote the subset of 
$\mathcal {I}rr(\A)$ that is mapped to the subset 
$\mathcal {STD}_d^{(b,b^*)}$ of $\mathcal {STD}$ 
by the bijection $\rho \mapsto A=\rho(z),\, A^*=\rho(z^*)$ 
from $\mathcal {I}rr(\A)$ to $\mathcal {STD}$ 
(see Proposition \ref{prop: STD}). 
Set $b=st, ~b^*=st^{-1}$ as in (\ref{st}).  
Such scalars $s, \,t$ are determined by $b, \,b^*$ 
uniquely up to the $\pm$ sign. We choose one of them and fix it. 
Then by Proposition \ref{prop: STDd(b,b*)}, 
$\mathcal {I}rr_d^{s, t}(\T)$ is mapped to 
$\mathcal {STD}_d^{(b,b^*)}$ by the bijection 
$\rho \mapsto A=\rho(z_t), \,A^*=\rho(z_t^*)$. 
This means that if a finite-dimensional irreducible representation 
$\rho: \T \longrightarrow {\rm{End}}(V)$ 
belongs to $\mathcal {I}rr_d^{s, t}(\T)$, then 
$$\rho'= \rho \circ \iota_t : \A \longrightarrow {\rm{End}}(V)$$
is a finite-dimensional irreducible representation of 
$\A$ that belongs to $\mathcal {I}rr_d^{(b,b^*)}(\A)$, 
where 
$$ \iota_t : \A \longrightarrow \T 
\qquad ( z, \,z^* \mapsto z_t, \,z^*_t) $$ 
is the injective algebra homomorphism from Proposition \ref{prop: iota}. 
Moreover every finite-dimensional irreducible representation of 
$\A$ that belongs to $\mathcal {I}rr_d^{(b,b^*)}(\A)$ arises 
in this way. In other words, 
a finite-dimensional irreducible representation 
$\rho': \A \longrightarrow {\rm{End}}(V)$ that belongs to 
$\mathcal {I}rr_d^{(b,b^*)}(\A)$ can be `extended' via $\iota_t$ 
to a finite-dimensional irreducible representation 
$\rho: \T \longrightarrow {\rm{End}}(V)$ 
that belongs to $\mathcal {I}rr_d^{s, t}(\T)$. Thus we have 
\begin{corollary}
\label{cor: Irr(A)}
If $b=st, ~b^*=st^{-1}$, then 
the mapping $\rho \mapsto \rho \circ \iota_t$ 
gives a bijection from $\mathcal {I}rr_d^{s, t}(\T)$ 
to $\mathcal {I}rr_d^{(b,b^*)}(\A)$, where 
$ \iota_t : \A \longrightarrow \T ~ ( z, \,z^* \mapsto z_t, \,z^*_t)$ 
is the injective algebra homomorphism from Proposition \ref{prop: iota}.
\end{corollary}
\medskip
Since $\mathcal {I}rr(\A)$ is the disjoint union of 
the trivial representations and 
$\mathcal {I}rr_d^{(b,b^*)}(\A)$ over $d \in \N$ and 
$(b, b^*) \in 
(\C \backslash \{0\})\times (\C \backslash \{0\})/\sim$, 
Problem \ref{problem: A} is reduced to Problem \ref{problem: T} 
by Corollary \ref{cor: Irr(A)}. 
Namely, $\mathcal {I}rr(\A)$ is the disjoint union of 
the trivial representations and 
\begin{eqnarray}
\label{rho circ iota}
\{\rho \circ \iota_t~|~ \rho \in \mathcal {I}rr_d^{s, t}(\T) \} 
\end{eqnarray}
over $d \in \N$ and 
$(s, t) \in (\C \backslash \{0\}) \times 
(\C \backslash \{0\}) /\approx$, 
where the equivalence relation $\approx$ is defined by 
$(s, t) \approx (s', t')$ if and only if 
\begin{eqnarray}
\label{approx 1st kind}
(s', t')&\in& 
\{\pm (s,t),\, \pm (t^{-1},s^{-1}), \, \pm (t,s), 
\, \pm (s^{-1},t^{-1})\}~~~
{\rm for~the~case} ~(\varepsilon, \varepsilon^*)=(1,1), \\
\label{approx 2nd kind}
(s', t')&\in& 
\{\pm (s,t),\, \pm (t^{-1},s^{-1})\}~~~
{\rm for~the~case} ~(\varepsilon, \varepsilon^*)=(1,0), 
\end{eqnarray}
and $(s', t')=\pm (s,t)$ 
for the case $(\varepsilon, \varepsilon^*)=(0,0)$.

\medskip
As we see in the next proposition, 
the property $({\rm C_1})$ for $\mathcal {I}rr(\A)$ 
is automaically satisfied 
when $(\varepsilon,\varepsilon^*)=(1,1)$. 
\begin{proposition}
\label{prop: q-Onsager}
If $(\varepsilon,\varepsilon^*)=(1,1)$, then every 
finite-dimensional irreducible representation 
$\rho: \A \longrightarrow {\rm{End}}(V) $ satisfies the property 
$({\rm C_1})$, i.e., $\rho(z),\,\rho(z^*)$ are diagonalizable. 
\end{proposition}

\noindent
Proof. Regard $V$ as an irreducible $\A$ -module via $\rho$ .
For $\theta \in \C$, set $V(\theta)=\{v \in V | zv=\theta v\}$. 
Note that $\theta$ is an eigenvalue of $\rho(z)$ if and only if 
$V(\theta) \neq 0$,
and in this case $V(\theta)$ is the correspording eigenspace.
Using the relation 
$[z,~z^2z^*-\beta z\,z^*z+z^*z^2] = \delta \,[z,~z^*]$ 
with $\beta=q^2+q^{-2}$, $\delta=-(q^2-q^{-2})$, we find 
$(z- \theta^- )(z- \theta)(z- \theta^+) z^* v =0$ 
for all $v \in V(\theta)$, where $\theta=\zeta + \zeta^{-1},\,
\theta^+ = q^2 \zeta + q^{-2} \zeta^{-1}, \,
\theta^- = q^{-2} \zeta + q^2 \zeta^{-1}$ , i.e., 
$$ z^* V(\theta) \subseteq V(\theta^-) + V(\theta) +V(\theta^+).$$
Set $\theta^{(i)} = q^{2i} \zeta + q^{-2i} \zeta^{-1}$. 
Then $\sum_{i \in \Z} V(\theta^{(i)})$,
which is a finite sum by dim $V < \infty$, is invariant under $z,\,z^*$.
Since $z,\,z^*$ generate $\A$ and $V$ is irreducible as an $\A$ -module,
we have $V=\sum_{i \in \Z} V(\theta^{(i)})$.
This implies that $\rho(z)$ is diagonalizable.
Similarly, $\rho(z^*)$ is shown to be diagonalizable. 
\hfill $\Box $

\medskip
So for the case $(\varepsilon,\varepsilon^*)=(1,1)$, 
Problem \ref{problem: A} is equivalent to 
\begin{problem}
\label{problem: q-Onsager}
Classify up to isomorphism 
the finite-dimensional irreducible representations of 
the q-Onsager algebra $\A_{\rm I}$. 
\end{problem}

\noindent
Thus the classification of standardized TD-pairs 
of Type I that are the 1st kind 
is equivalent to that of finite-dimensional 
irreducible representations of the $q$-Onsager algebra 
$\A_{\rm I}$.

\subsection{ Finite-dimensional irreducible $\T$-modules} 
Let $\rho: \T \longrightarrow {\rm{End}}(V)$ be a 
finite-dimensional irreducible representation of the 
augmented TD-algebra $\T$.
We regard $V$ as an irreducible $\T$-module via $\rho$. 
Let us recall Lemma~\ref{lemma: weight-space decomposition} in Section 1.2. 
The action of $k$ on $V$ is diagonalizable with eigenvalues 
$\{ s q^{2i-d} | 0 \leq i \leq d \}$ for some nonzero  
$s \in \C$ and an integer $d \geq 0$.
The scalar $s$ and the integer $d$ are called the type 
and the diameter, respectively. 
Let $V=\bigoplus^d_{i=0} U_i$ denote the eigenspace 
decomposition of the action of $k$ on $V$, where 
$U_i$ is the eigenspace belonging to $s q^{2i-d}$. 
It holds that $x U_i \subseteq U_{i+1}, \;y 
U_i \subseteq U_{i-1} \;(0 \leq i \leq d)$, 
where $U_{-1}=U_{d+1}=0$. We call the direct sum 
$V=\bigoplus^d_{i=0} U_i$ the weight space decomposition 
and $U_0$ the highest weight space. 

\begin{theorem}
\label{thm: shape}
Let $V$ be a finite-dimensional irreducible $\T$-module and 
$V=\bigoplus^d_{i=0} U_i$ the weight space decomposition.
Then 
\[
{\rm dim}\, U_i \leq 
{d \choose i}
\qquad  (0 \leq i \leq d).
\]
In particular $U_0$ has dimension 1. 
\end{theorem}

Theorem~\ref{thm: shape}  will be proved in Section 3. 
Since 
$x \,U_j \subseteq U_{j+1},\; 
y\, U_j \subseteq U_{j-1} \;(0 \leq j \leq d)$
 by Lemma~\ref{lemma: weight-space decomposition}, 
the highest weight space $U_0$ is invariant 
under $y^i x^i$ for every integer $i \geq 0$. 
Since dim $U_0 =1$ by Theorem~\ref{thm: shape}, 
there exists $\sigma_i=\sigma_i(V) \in \C$
such that 
$$ y^i x^i v = \sigma_i v ~~~~~ (v \in U_0) $$
for every integer $i \geq 0$. Observe $\sigma_0=1$ and 
$\sigma_i=0$ if $i>d$, where $d$ is the diameter of the $\T$-module $V$. 
It is shown later that $\sigma_d \neq 0$. 
Let $\mathcal {M}^s_d (\T)$ denote the set of isomorphism classes of  
finite-dimensional irreducible $\T$-modules with type $s$, 
diameter $d$, and  $\Sigma_d$ the set of sequences 
${\{\sigma_i\}}^d_{i=0}$ of scalars $\sigma_i \in \C$ with 
$\sigma_0=1$, $\sigma_d \neq 0$.
Then we have a mapping $\sigma$ from $\mathcal {M}^s_d (\T)$  to 
$\Sigma_d$ that sends $V$ to ${\{ \sigma_i(V) \}}^d_{i=0}$, 
where $\sigma_i(V)$ is the eigenvalue of $y^i x^i$  
on the highest weight space of $V$. 

\begin{theorem}
\label{thm: sigma}
For each nonzero $s \in \C$, the mapping 
$$ \sigma : \mathcal{M}^s_d (\T) \longrightarrow \Sigma_d ~~~
(V \longmapsto {\{ \sigma_i(V)\}}^d_{i=0} ) $$
is a bijection. 
\end{theorem}

The fact $\sigma_d(V) \neq 0$ and the injectiveness of 
$\sigma$ will be proved in Section 3. 
The surjectiveness of $\sigma$ will be proved in Section 5.

\medskip
For a finite-dimensional irreducible $\T$-module $V$
of type $s$ and diameter $d$, we define a monic polynomial 
$P_V(\lambda)$ of degree $d$ in $\lambda$ as follows:
\begin{eqnarray*}
P_V(\lambda)=Q^{-1} \sum^d_{i=0} \sigma_i(V) \prod^d_{j=i+1} 
{(q^j - q^{-j})}^2 (\varepsilon s^{-2} q^{2(d-j)} 
+ \varepsilon^*s^2 q^{-2(d-j)} - \lambda), 
\end{eqnarray*}
where $\sigma_i(V)$ is the eigenvalue of $y^i x^i$  on the 
highest weight space of $V$ and 
$$ Q = Q_d = {(-1)}^d {(q - q^{-1})}^2 {(q^2 - q^{-2})}^2 \cdots
{(q^d - q^{-d})}^2.$$ 
The polynomial $P_V (\lambda)$ is called the \textit{Drinfel'd polynomial} 
of the  $\T$-module  $V$. 
Note that the parameters $q$ and $(\varepsilon , \varepsilon^*)$ 
in the definition of $P_V (\lambda)$ are independent 
of the $\T$-module $V$, since they are chosen and fixed in advance 
for the augmented TD-algebra  $\T$.

\begin{remark}
\label{rem: Drinfeld polynomial} 
\rm{

The following identities directly follow 
from the definition of $P_V (\lambda)$. 
\begin{enumerate}
\item
For $\lambda = \varepsilon s^{-2} + \varepsilon^* s^2$, 
$$P_V (\lambda) = Q^{-1} \sigma_d(V).
~~~~~~~~~~~~~~~~~~~~~~~~~~~~~~~~~~~~~~~~~~~~~~~~~~~~~~~~~~~~~ $$
\item
For $\lambda = t^2 +  \varepsilon \varepsilon^* t^{-2}$ 
with $t$ an arbitrary nonzero scalar, 
$$P_V (\lambda) = Q^{-1} \sum^d_{i=0} \sigma_i(V) 
(\theta_0 - \theta_{i+1}) \cdots 
(\theta_0 - \theta_d)(\theta_0^* - \theta^*_{i+1}) \cdots 
(\theta_0^* - \theta_d^*),$$ 
where 
$\theta_i = s\,t \,q^{2i-d} + \varepsilon s^{-1} t^{-1} q^{d-2i}, \:
\theta^*_i = \varepsilon^* s\, t^{-1} q^{2i-d} + s^{-1} t\, q^{d-2i}$. 
\end{enumerate}

}
\end{remark}

\medskip
Let $\mathcal P^s_d$ denote the set of monic polynomials $P(\lambda)$ 
of degree $d$ in $\lambda$ such that 
$$ P(\lambda) \neq 0 ~~ \mbox{for} ~~ \lambda = \varepsilon s^{-2} + 
\varepsilon^* s^2. $$
Then the mapping that sends ${\{ \sigma_i \}}^d_{i=0} $ to 
$$ P(\lambda) = Q^{-1} \sum^d_{i=0} \sigma_i 
\prod^d_{j=i+1} {(q^j - q^{-j})}^2 (\varepsilon s^{-2} q^{2(d-j)} 
+ \varepsilon^* s^2 q^{-2(d-j)} - \lambda ) $$
gives a bijection from $\Sigma_d$ to $\mathcal P^s_d$. 
So we can restate Theorem \ref{thm: sigma} as follows.

\bigskip
\noindent
${\bf Theorem~ {\bf \ref{thm: sigma}}}'$ 
{\it 
The mapping $V \longmapsto P_V (\lambda)$ 
gives a bijection from 
$\mathcal {M}^s_d(\T)$ to $\mathcal P^s_d$. 
}

\bigskip 
This gives a parametrization of the set $\mathcal {M}^s_d(\T)$ 
in question in Problem~\ref{problem: T}~(i). 

\begin{theorem}
\label{thm: C2}
Let $V$ be a finite-dimensional irreducible 
$\T$-module of type $s$ and diameter $d$. 
Assume that the property ${({\rm C}_1)}_t$ holds for some 
$t \in \C$, i.e., the action of $z_t, \,z_t^*$ on $V$ 
are both diagonalizable. Then $V$ is irreducible as 
a $\langle z_t, z_t^* \rangle$-module 
if and only if $P_V (\lambda) \neq 0$ 
for $\lambda = t^2 + \varepsilon \varepsilon^* t^{-2}$. 
Here $P_V (\lambda)$ is the Drinfel'd polynomial of the 
$\T$-module $V$. 
\end{theorem}

Theorem~\ref{thm: C2} will be proved in Section 4. 
Theorem~\ref{thm: C2} together with Proposition~\ref{prop: C1}  
gives a parametrization of the representations of $\T$ in question 
in Problem~\ref{problem: T}~(ii). 
For an integer $d \geq 1$ and nonzero $s,\, t \in \C$, 
define the sets $\mathcal {M}_d^{s,t}(\T)$ and $\mathcal {P}^{s,t}_d$ 
as follows. 
$\mathcal {M}_d^{s,t}(\T)$ denotes 
the set of isomorphism classes 
of finite-dimensional irreducible $\T$-modules $V$ 
with type $s$, diameter $d$ 
that satisfy the properties 
$(\rm {C}_1)_t$, $(\rm {C}_2)_t$, i.e., 
the action of $z_t, \,z_t^*$ on $V$ 
are both diagonalizable 
and $V$ is irreducible as a $\langle z_t, z_t^* \rangle$-module. 
$\mathcal {P}^{s,t}_d$ denotes the set of monic 
polynomials $P(\lambda)$ of degree $d$ in $\lambda$ 
such that $P(\lambda) \neq 0 $ for 
$\lambda= \varepsilon s^{-2} +\varepsilon^* s^2$ 
and $\lambda =t^2 + \varepsilon \varepsilon^* t^{-2}$. 
Note that 
$\mathcal {M}_d^{s,t}(\T)$ (resp. $\mathcal {P}^{s,t}_d(\T)$) is   
a subset of $\mathcal {M}^s_d(\T)$ (resp. $\mathcal {P}^s_d(\T)$) 
and $\mathcal {M}^s_d(\T)$ is bijectively mapped to 
$\mathcal {P}^s_d$ by $V \longmapsto P_V (\lambda)$ 
by Theorem $\ref{thm: sigma}'$. 
Let $V$ be a finite-dimensional irreducible $\T$-module 
that belongs to $\mathcal {M}^s_d(\T)$. 
Then by Proposition~\ref{prop: C1}, 
the property $(\rm {C}_1)_t$ holds for the $\T$-module $V$ 
if and only if 
\begin{eqnarray}
\label{C1zt}
s\,t \neq \pm \varepsilon q^i ~~~
{\rm for~any~integer} ~i ~~(1-d \leq i \leq d-1),\\ 
\label{C1z*t}
s\,t^{-1} \neq \pm \varepsilon^* q^i ~~~
{\rm for ~any ~integer} ~i ~~(1-d \leq i \leq d-1). 
\end{eqnarray}
Thus if one of the conditions (\ref{C1zt}), (\ref{C1z*t}) 
fails, then $\mathcal {M}_d^{s,t}(\T)$ is empty. 
Suppose each of (\ref{C1zt}), (\ref{C1z*t}) holds. 
Then by Theorem~\ref{thm: C2}, 
the property $(\rm {C}_2)_t$ holds for the $\T$-module $V$ 
if and only if 
$P_V (\lambda) \neq 0$ 
for $\lambda = t^2 + \varepsilon \varepsilon^* t^{-2}$. 
So $\mathcal {M}_d^{s,t}(\T)$ is precisely mapped 
onto $\mathcal {P}^{s,t}_d$  by the bijection 
$V \longmapsto P_V (\lambda)$ from $\mathcal {M}_d^s(\T)$ 
to $\mathcal {P}^s_d$. 
Thus we have 

\begin{corollary} 
\label{cor: M}
If one of the conditions $(\ref{C1zt})$, $(\ref{C1z*t})$ 
fails, then $\mathcal {M}_d^{s,t}(\T)$ is empty.
Suppose each of $(\ref{C1zt})$, $(\ref{C1z*t})$ holds. 
Then the mapping $V \longmapsto P_V (\lambda)$ 
gives a bijection from 
$\mathcal {M}^{s,t}_d(\T)$ to $\mathcal P^{s,t}_d$. 
\end{corollary}

This gives a parametrization of the set $\mathcal {M}^{s,t}_d(\T)$ 
in question in Problem~\ref{problem: T}~(ii). 
Since $\mathcal {M}_d^{s,t}(\T)$ can be natually identified 
with $\mathcal {I}rr_d^{s,t}(\T)$, 
Corollary \ref{cor: M} gives 
a parametrization of $\mathcal {STD}_d^{(b,b^*)}$ 
through Proposition~\ref{prop: STDd(b,b*)}, 
where $b=s\,t, \, b^*=s\,t^{-1}$.

\subsection{Construction of finite-dimensional irreducible $\T$-modules} 
Given 
$(\varepsilon, \varepsilon^*) \in \{(1,1), (1,0), (0,0)\}$ 
and a nonzero scalar $q$ that is not a root of unity, 
let $\T=\T_q^{(\varepsilon, \varepsilon^*)}$ 
denote the augmented TD-algebra. $\T$ is generated by 
$x,y, k^{\pm 1}$ subject to the relations $({\rm TD})'_0$, $({\rm TD})'$ 
in Section 1.2.  
In the next proposition, we give an injective algebra-homomorphism 
$\varphi_s$ of $\T$ into the $U_q(sl_2)$-loop algebra 
${\LL}=U_q(L(sl_2))$ 
for each nonzero scalar $s \in \C$. 
$\LL$ is the associative 
$\C$-algebra with 1 defined by generators 
$e^+_i, e^-_i, k_i, k^{-1}_i ~(i=0,1)$ subject to the relations 
\begin{eqnarray*}
k_0 k_1 &=& k_1 k_0 = 1, \\
k_i k^{-1}_i &=& k^{-1}_i k_i = 1, \\
k_i e^{\pm}_i k^{-1}_i &=& q^{\pm 2} e^{\pm}_i, \\
k_i  e^{\pm}_j k_i^{-1} &=& q^{\mp 2} e^{\pm}_j ~~~~(i \neq j), \\
\lbrack e^+_i, e^-_i \rbrack &=& \frac{{k_i} - k^{-1}_i}{q-q^{-1}}, \\
\lbrack e^+_i,e^-_j\rbrack&=&0 ~~~~(i \neq j), 
\end{eqnarray*}
\vspace{-9mm}
\begin{eqnarray*}
\lbrack e^{\pm}_i, (e^{\pm}_i)^2 e^{\pm}_j - (q^2 + q^{-2})
e^{\pm}_i e^{\pm}_j e^{\pm}_i + e^{\pm}_j (e^{\pm}_i )^2\rbrack=0 ~~~~
(i \neq j).
\end{eqnarray*}
Note that if we replace $k_0 k_1 = k_1 k_0 = 1$ 
in the defining relations for $\LL$ by $k_0 k_1 = k_1 k_0$, 
then we have the quantum affine algebra $U_q(\widehat{sl}_2)$: $\LL$ 
is isomorphic to the quotient algebra of $U_q(\widehat{sl}_2)$ 
by the two-sided ideal generated by $k_0 k_1 -1$.

\begin{proposition}
\label{prop: phi}
For each nonzero $s \in \C$, there exists an 
algebra homomorphism $\varphi_s$ from $\T$ to $\LL$ that sends 
$x,\, y, \,k$ to $x(s),\, y(s),\, k(s)$, respectively, where 
\begin{eqnarray*}
x(s) & =& \alpha(se_0^+ + \varepsilon s^{-1} e_1^- k_1) 
~~~ {\rm with} ~\alpha = - q^{-1} {(q-q^{-1})}^2,\\
y(s) & =& \varepsilon^* s e_0^- k_0 + s^{-1} e_1^+, \\
k(s) & =& s k_0. 
\end{eqnarray*}
Moreover $\varphi_s$ is injective. 
\end{proposition}

The existence of $\varphi_s$ follows from the fact that 
the relations $({\rm TD})'_0$, $({\rm TD})'$  
hold for $x(s)$, $y(s)$, $k(s)$, $k(s)^{-1}$. 
We leave the tedious calculations of checking it to the reader. 
The injectivity of $\varphi_s$ will be proved in Section 2.

\medskip
Let $\LL '$ denote the subalgebra of $\LL$ generated by 
$e^+_0$, $e^{\pm}_1$, $k^{\pm 1}_i$ $(i=0,1)$: $e^-_0$ is missing 
from the set of generators for $\LL '$. 
Let $\B$ denote the subalgebra of $\LL$ generated by 
$e^+_0$, $e^+_1$, $k^{\pm 1}_i$ $(i=0,1)$, 
the Borel subalgebra of $\LL$. 
Observe $\B \subseteq \LL'$. 
Note that the image of $\varphi_s$ is contained in 
$\LL'$ if $( \varepsilon,  \varepsilon^*)=(1,0)$ and it 
coincides with $\B$ if $( \varepsilon,  \varepsilon^*)=(0,0)$.
If $(\varepsilon,  \varepsilon^*)=(1,1)$ (resp. $(1,0), \,(0,0)$), 
each finite-dimensional irreducible $\LL$-module 
(resp. $\LL '$-module, $\B$-module) 
can be regarded as a $\T$-module via the injective algebra 
homomorphism $\varphi_s: \T \longrightarrow \LL$. Such a $\T$-module 
is called a {\it $\T$-module via $\varphi_s$}. 
We determine when a finite-dimensional irreducible $\LL$-module 
(resp. $\LL '$-module, $\B$-module) remains irreducible 
as a $\T$-module via $\varphi_s$, 
and by finding an explicit formula for the Drinfel'd polynomial  
$P_V(\lambda)$, we show that every finite-dimensional 
irreducible $\T$-module with type $s$ arises in this way 
via $\varphi_s$ (see Theorem~$\ref{thm: sigma}'$).

\medskip
We give an overview of finite-dimensional representations of 
$\LL$ that we need to state our explicit construction of 
irreducible $\T$-modules via $\varphi_s$. 
For $a \in \C ~(a \neq 0) $ and $\ell \in \Z ~(\ell > 0)$, 
$V(\ell,a)$ denotes the {\it evaluation module} 
of $\LL$, i.e., $V(\ell,a)$ is an $(\ell + 1)$-dimensional 
vector space over $\C$ with a basis 
$v_0, v_1, \ldots , v_{\ell}$ on which $\LL$ acts as follows: 
\begin{eqnarray*}
k_0 v_i  &=& q^{2i-\ell} \,v_i ,\\
k_1 v_i  &=& q^{\ell-2i} v_i, \\
e^+_0 v_i  &=& a\, q \,[i+1] \,v_{i+1}, \\
e^-_0 v_i  &=& a^{-1} q^{-1} [\ell-i+1]\, v_{i-1}, \\
e^+_1 v_i  &=&  [\ell-i+1] \,v_{i-1}, \\
e^-_1 v_i  &=& [i+1]\, v_{i+1},
\end{eqnarray*}
where $v_{-1} = v_{\ell+1} = 0 $ and 
$ [j] = [j]_q=(q^j - q^{-j})/(q - q^{-1})$. 
$V(\ell,a)$  is an irreducible  $\LL$-module. 
We call $v_0, v_1, \ldots , v_{\ell}$ a {\it standard basis}.

\medskip
Let $\Delta$ denote the {\it coproduct} of
$\LL$: the algebra homomorphism from $\LL$ to 
$\LL \otimes \LL$ defined by  
\begin{eqnarray*}
\Delta (k^{\pm 1}_i) & =& k^{\pm 1}_i \otimes k^{\pm 1}_i ,\\
\Delta (e^+_i) & =& k_i \otimes e^+_i + e^+_i \otimes 1, \\
\Delta (e^-_i k_i) & =& k_i \otimes e^-_i k_i + e^-_i k_i \otimes 1.
\end{eqnarray*}
Given $\LL$-modules $V_1,V_2$, the tensor product $V_1 \otimes V_2$ 
becomes an $\LL$-module via $\Delta$.
Given a set of evaluation modules $V(\ell_i,a_i)~(1 \leq i \leq n)$ 
of $\LL$, the tensor product 
$$ V(\ell_1,a_1) \otimes \cdots \otimes V(\ell_n,a_n) $$
makes sense as an $\LL$-module without being affected by   
the parentheses for the tensor product 
because of the coassociativity of $\Delta$. 

\medskip
With an evaluation module $V(\ell, a)$ of $\LL$, 
we associate the set $S(\ell, a)$ of scalars 
$a\,q^{-\ell+1}$, $a\,q^{-\ell+3}$, $\cdots$, $a\,q^{\ell-1}$: 
$$ S(\ell,a) = \{ a\,q^{2i-\ell+1} ~|~ 0 \leq i \leq \ell-1 \}. $$
The set $S(\ell, a)$ is called a $q$-$string$ of length $\ell$.
Two $q$-strings $S(\ell, a)$, $S(\ell', a')$ are said to be 
{\it adjacent} 
if $S(\ell, a) \cup S(\ell', a')$ is a longer $q$-string, i.e., 
$S(\ell, a) \cup S(\ell', a') = S(\ell'', a'')$ 
for some $\ell''$, $a''$ with $\ell'' > {\rm max} \{\ell,~\ell '\}$. 
It can be easily checked that $S(\ell, a)$, $S(\ell', a')$ are 
adjacent if and only if $a^{-1}a'=q^{\pm i}$ for some 
$$i \in 
\{|\ell -\ell'|+2,\, |\ell -\ell'|+4, \cdots, \ell +\ell' \}.$$ 
Two $q$-strings $S(\ell, a)$, $S(\ell', a')$ are defined to be 
{\it in general position} if they are not adjacent, i.e., 
if either 
\begin{quote}
$(i)$ ~
$S(\ell, a) \cup S(\ell', a')$ is not a $q$-string, \\
or \\
$(ii)$ $\,$
$S(\ell, a) \subseteq S(\ell', a')$ or 
$S(\ell, a) \supseteq S(\ell', a')$. 
\end{quote}
A multi-set ${\{ S(\ell_i,a_i) \}}^n_{i=1}$ 
of $q$-strings 
is said to be {\it in general position} if 
$S(\ell_i,a_i)$ and $S(\ell_j,a_j)$ are in general position for 
any $i,\,j ~(i \neq j,\, 1 \leq i \leq n,\, 1 \leq j \leq n)$. 
The following fact is well-known and easy to prove. 
Let $\Omega$ be a finite multi-set of nonzero scalars from $\C$. 
Then there exists a multi-set 
$\{ S(\ell_i,a_i)\}_{i=1}^n$ 
of $q$-strings in general position such that 
$$ \Omega = \bigcup_{i=1}^n ~S(\ell_i,a_i)$$
as multi-sets of nonzero scalars. 
Moreover such a multi-set of $q$-strings is uniquely determined by $\Omega$. 

\medskip
With a tensor product $V(\ell_1,a_1) \otimes \cdots \otimes V(\ell_n,a_n)$ 
of evaluation modules $V(\ell_i,a_i)~ (1 \leq i \leq n)$, 
we associate the multi-set ${\{ S(\ell_i,a_i) \}}^n_{i=1} $ 
of $q$-strings. 
The following (i), (ii), (iii) are well-known \cite{CP}: 
\begin{quote}
(i)~\,A tensor product 
$V(\ell_1,a_1) \otimes \cdots \otimes V(\ell_n,a_n)$ 
of evaluation modules is irreducible as an $\LL$ -module 
if and only if 
the multi-set ${\{ S(\ell_i,a_i)\}}^n_{i=1}$ of $q$-strings is 
in geneal position. \\
\\
(ii)~\,Set $V=V(\ell_1,a_1) \otimes \cdots \otimes V(\ell_n,a_n)$, 
$V'= V(\ell_1',a_1') \otimes \cdots V(\ell_{n'}',a_{n'}')$ and assume 
that $V$, $V'$ are both irreducible as an $\LL$-module. 
Then $V$, $V'$ are isomorphic as $\LL$-modules if and only if 
the multi-sets ${\{ S(\ell_i,a_i)\}}^n_{i=1}$, 
${\{S(\ell_i',a_i')\}}^{n'}_{i=1}$ 
coincide, i.e., $n=n'$ and $\ell_i=\ell_i'$, $a_i=a_i'$ 
for all $i$ $(1 \leq i \leq n)$ 
with a suitable reordering of 
$S(\ell_1',a_1'), \cdots, S(\ell_n',a_n')$. \\
\\
(iii)~\,Every nontrivial 
finite-dimensional irreducible $\LL$-module of 
type (1,1)  is isomorphic to some 
$V(\ell_1,a_1) \otimes \cdots V(\ell_n,a_n)$.\\
\end{quote}

Two multi-sets 
${\{S(\ell_i,a_i)\}}^n_{i=1}, {\{S(\ell_i',a_i')\}}^{n'}_{i=1}$ 
of $q$-strings are defined to be 
{\it equivalent} if there exists 
$\varepsilon_i \in \{\pm 1\} ~(1 \leq i \leq n)$ 
such that 
${\{S(\ell_i,a_i^{\varepsilon_i})\}}^n_{i=1}$
and 
${\{S(\ell_i',a_i')\}}^{n'}_{i=1}$ 
coincide, i.e., $n=n'$ and 
$\ell_i=\ell_i', \, a_i^{\varepsilon_i}=a_i'$ 
for all $i~(1 \leq i \leq n)$ 
with a suitable reordering of 
$S(\ell_1', a_1'), \cdots, S(\ell_n',a_n')$. 
A multi-set 
${\{S(\ell_i,a_i)\}}^n_{i=1}$
of $q$-strings is defined to be {\it strongly in general position} 
if any multi-set of $q$-strings equivalent to 
${\{S(\ell_i,a_i)\}}^n_{i=1}$ is in general position, i.e., 
the multi-set 
${\{S(\ell_i,a_i^{\varepsilon_i})\}}^n_{i=1}$
is in general position for any choice of 
$\varepsilon_i \in \{ \pm 1 \} ~(1 \leq i \leq n)$.

\begin{lemma}
\label{lemma: q-string}
Let $\Omega$ be a finite multi-set of nonzero scalars from $\C$ 
such that c and $c^{-1}$ appear in $\Omega$ in pairs, i.e., c and 
$c^{-1}$ have the same multiplicity in $\Omega$ for each $c \in \Omega$, 
where we understand that if 1 or -1 appears in $\Omega$, 
it has even multiplicity. 
Then there exists a multi-set 
${\{S(\ell_i,a_i) \}}_{i=1}^n$ of $q$-strings 
strongly in general position such that 
$$ \Omega = \bigcup_{i=1}^n ~
\bigl(S(\ell_i,a_i) \cup S(\ell_i,a_i^{-1})\bigr) $$
as multi-sets of nonzero scalars. 
Such a multi-set of $q$-strings is uniquely determined by $\Omega$ 
up to equivalence. 
\end{lemma}

Lemma \ref{lemma: q-string} will be proved in Section 7. 

\begin{theorem}
\label{thm: 1st kind}
{\rm (Case ($\varepsilon,\varepsilon^*) = (1,1)$)}~~
Let $\T=\T_q^{(1,1)}$ denote the augmented TD-algebra of the 1st kind. 
The following $(i)$, $(ii)$, $(iii)$ hold. 

\begin{enumerate}
\item[$(i)$]
A tensor product 
$V(\ell_1,a_1)\otimes \cdots \otimes V(\ell_n,a_n)$ 
of evaluation modules is irreducible as a $\T$-module via 
$\varphi_s$ if and only if 
$-s^2 \notin S(\ell_i,a_i) \cup S(\ell_i,a_i^{-1})$ 
for all $i ~(1 \leq i \leq n)$ and the multi-set 
${\{ S(\ell_i,a_i)\}}^n_{i=1} $
of $q$-strings is strongly in general position. 
In this case, the $\T$-module 
$V=V(\ell_1,a_1)\otimes \cdots \otimes V(\ell_n,a_n)$ 
via $\varphi_s$ has type $s$ and diameter 
$d=\ell_1+\cdots+\ell_n$ and the Drinfel'd polynomial $P_V(\lambda)$ 
of the $\T$-module $V$ via $\varphi_s$ is 
$$ P_V(\lambda) = \prod^n_{i=1} P_{V(\ell_i,a_i)}(\lambda),$$
where
$$ P_{V(\ell_i,a_i)}(\lambda) = 
\prod_{c \in S(\ell_i,a_i)}(\lambda +c+c^{-1}).$$
\item[$(ii)$] 
Set $V=V(\ell_1,a_1) \otimes \cdots \otimes V(\ell_n,a_n)$, 
$V'= V(\ell_1',a_1') \otimes \cdots \otimes V(\ell_{n'}',a_{n'}')$ 
and assume that $V$, $V'$ are both irreducible 
as a $\T$-module via $\varphi_s$. 
Then $V$, $V'$ are isomorphic as $\T$-modules via $\varphi_s$ 
if and only if 
the multi-sets ${\{ S(\ell_i,a_i)\}}^n_{i=1}$, 
${\{S(\ell_i',a_i')\}}^{n'}_{i=1}$ of $q$-strings are equivalent. 
\item[$(iii)$] 
Every nontrivial 
finite-dimensional irreducible $\T$-module of 
type $s$ is isomorphic to some $\T$-module 
$V(\ell_1,a_1) \otimes \cdots \otimes V(\ell_n,a_n)$ 
via $\varphi_s$ . 
\end{enumerate}
\end{theorem}

Theorem \ref{thm: 1st kind} will be proved in Section 7. 
Note that the Drinfel'd polynomial of an irreducible $\T$-module 
$V(\ell_1,a_1) \otimes \cdots \otimes V(\ell_n,a_n)$ via $\varphi_s$ 
is determined by the multi-set ${\{ S(\ell_i,a_i)\}}^n_{i=1}$ 
of $q$-strings and independent of $\varphi_s$. 
Problem \ref{problem: T}, which is to determine $\mathcal {M}^s_d(\T)$ 
and $\mathcal {M}_d^{s,t}(\T)$, is solved by Theorem \ref{thm: 1st kind} 
as follows in the case of $(\varepsilon,\varepsilon^*) = (1,1)$. 
Assume $(\varepsilon,\varepsilon^*) = (1,1)$. 
The set $\mathcal {M}^s_d(\T)$ is determined in terms of 
tensor products of evaluation modules 
by Theorem \ref{thm: 1st kind} (i), (ii), (iii). 
Recall the bijection 
$V \mapsto P_V (\lambda)$ 
from $\mathcal {M}^s_d(\T)$ to $\mathcal P^s_d$ 
in Theorem~$\ref{thm: sigma}'$. 
The subset $\mathcal {M}_d^{s,t}(\T)$ of $\mathcal {M}^s_d(\T)$ 
is nonempty if and only if the conditions (\ref{C1zt}), (\ref{C1z*t}) 
hold, i.e., 
\begin{eqnarray}
\label{C1 1st kind}
\pm st,\, \pm st^{-1} \notin \{q^i~|~ 
i= -d+1,\, -d+2, \cdots, d-1 \},
\end{eqnarray}
and in this case $\mathcal {M}_d^{s,t}(\T)$ is mapped onto 
$\mathcal P^{s,t}_d$ by the bijection 
$V \mapsto P_V (\lambda)$ (see Corollary~\ref{cor: M}). 
For an irreducible $\T$-module 
$V(\ell_1,a_1) \otimes \cdots \otimes V(\ell_n,a_n)$ via $\varphi_s$, 
we find by Theorem \ref{thm: 1st kind} (i) that 
$P_V (\lambda)$ does not vanish at $\lambda =t^2+t^{-2}$ 
if and only if $-t^2 \notin S(\ell_i,a_i) \cup S(\ell_i,a_i^{-1})$ 
for all $i ~(1 \leq i \leq n)$. Thus we have

\begin{corollary}
\label{cor: M 1st kind}
Assume $(\varepsilon,\varepsilon^*) = (1,1)$. Then 
$\mathcal {M}^s_d(\T)$ and $\mathcal {M}_d^{s,t}(\T)$ 
are determined as follows. 

\begin{enumerate}
\item[$(i)$]
$\mathcal {M}^s_d(\T)$ consists of the isomorphism classes 
of $\T$-modules  
$V(\ell_1,a_1) \otimes \cdots \otimes V(\ell_n,a_n)$ 
via $\varphi_s$ with the properties that  
\begin{description}
\item[$(i.1)$] 
the multi-set 
${\{ S(\ell_i,a_i)\}}^n_{i=1}$
of $q$-strings is strongly in general position, 
\item[$(i.2)$] 
$-s^2 \notin S(\ell_i,a_i) \cup S(\ell_i,a_i^{-1})$ 
for all $i ~(1 \leq i \leq n)$, 
\item[$(i.3)$] 
$d=\ell_1+\cdots+\ell_n$. 
\end{description}
The isomorphism classes of such $\T$-modules 
$V(\ell_1,a_1) \otimes \cdots \otimes V(\ell_n,a_n)$ 
via $\varphi_s$ are in one-to-one correspondence with 
the equivalence classes 
of the multi-sets ${\{ S(\ell_i,a_i)\}}^n_{i=1}$
of $q$-strings 
that have the properties $(i.1)$, $(i.2)$, $(i.3)$ above. 

\item[$(ii)$]
$\mathcal {M}_d^{s,t}(\T)$ is nonempty if and only if 
the condition $(\ref{C1 1st kind})$ holds. 
Suppose the condition $(\ref{C1 1st kind})$ holds. Then 
$\mathcal {M}_d^{s,t}(\T)$ consists of the isomorphism classes 
of $\T$-modules 
$V(\ell_1,a_1) \otimes \cdots \otimes V(\ell_n,a_n)$ 
via $\varphi_s$ with the properties $(i.1)$, $(i.2)$, $(i.3)$ above 
and 
\begin{description}
\item[$(ii.1)$]
$-t^2 \notin S(\ell_i,a_i) \cup S(\ell_i,a_i^{-1})$ 
for all $i ~(1 \leq i \leq n)$. 
\end{description}

\end{enumerate}
\end{corollary}

The next theorem follows from Corollary~\ref{cor: M 1st kind} 
and \cite[Proposition 7.15]{ITe}.
It solves Problem~\ref{problem: q-Onsager}, 
which is to determine 
the finite-dimensional irreducible representations of 
the $q$-Onsager algebra up to isomorphism. 
For an $\LL$-module $V$, let $\rho_V$ denote the representation of $\LL$ 
afforded by the $\LL$-module $V$. Then $\rho_V \circ \varphi_s$ is 
the representation of $\T$ afforded by the $\T$-module $V$ via $\varphi_s$,  
and $\rho_V \circ \varphi_s \circ \iota_t$ is a representation of $\A$, 
where $$ \iota_t : \A \longrightarrow \T 
\qquad ( z, \,z^* \mapsto z_t, \,z^*_t) $$
is the injective algebra homomorphism from Proposition \ref{prop: iota}.

\begin{theorem}
\label{thm: q-Onsager}
Assume $(\varepsilon,\varepsilon^*) = (1,1)$. Let $\A=\A_q^{(1,1)}$ 
denote the q-Onsager algebra. 
The following $(i)$, $(ii)$, $(iii)$ hold. 
\begin{enumerate}
\item[$(i)$]
For an $\LL$-module $V=V(\ell_1,a_1) \otimes \cdots \otimes V(\ell_n,a_n)$ 
and nonzero $s,\,t \in \C$, the representaiton 
$\rho_V \circ \varphi_s \circ \iota_t$ of $\A$ is irreducible if and only if 
\begin{description}
\item[$(i.1)$] 
the multi-set 
${\{ S(\ell_i,a_i)\}}^n_{i=1}$
of $q$-strings is strongly in general position, 
\item[$(i.2)$] 
none of $-s^2,\, -t^2$ belongs to $S(\ell_i,a_i) \cup S(\ell_i,a_i^{-1})$ 
for any $i$ ~$(1 \leq i \leq n)$, 
\item[$(i.3)$] 
none of the four scalars $\pm st,\, \pm st^{-1}$ equals $q^i$ 
for any $i \in \Z$ ~$(-d+1 \leq i \leq d-1)$, 
where $d=\ell_1+\cdots+\ell_n$. 
\end{description}

\item[$(ii)$]
For $\LL$-modules $V=V(\ell_1,a_1) \otimes \cdots \otimes V(\ell_n,a_n)$, 
$V'=V(\ell_1',a_1') \otimes \cdots \otimes V(\ell_{n'}',a_{n'}')$ 
and $(s,t),\, (s',t') 
\in (\C \backslash \{0\})\times (\C \backslash \{0\})$, 
set $\rho=\rho_V \circ \varphi_s \circ \iota_t$ and 
$\rho'=\rho_{V'} \circ \varphi_{s'} \circ \iota_{t'}$. 
Assume that the representations $\rho$, $\rho'$ of $\A$ are both irreducible. 
Then they are isomorphic as representations of $\A$ 
if and only if the multi-sets  
${\{ S(\ell_i,a_i)\}}^n_{i=1}$,  ${\{ S(\ell_i',a_i')\}}^{n'}_{i=1}$ 
are equivalent and $(s,t) \approx (s',t')$ in the sense of 
$({\rm \ref{approx 1st kind}})$ ,i.e., 
$$
(s', t') \in 
\{\pm (s,t),\, \pm (t^{-1},s^{-1}), \, \pm (t,s), 
\, \pm (s^{-1},t^{-1})\}.
$$ 

\item[$(iii)$]
Every nontrivial finite-dimensional irreducible representation of $\A$ 
is isomorphic to $\rho_V \circ \varphi_s \circ \iota_t$ 
for some $\LL$-module $V=V(\ell_1,a_1) \otimes \cdots \otimes V(\ell_n,a_n)$ 
and $(s,t) \in (\C \backslash \{0\})\times (\C \backslash \{0\})$. 
\end{enumerate}
\end{theorem}

\medskip
\noindent
Proof. The assertions (i), (iii) follow from 
Corollary \ref{cor: Irr(A)} and Corollary \ref{cor: M 1st kind}, 
since $\I rr_d^{s,t}(\T)$ is naturally identified 
with $\mathcal{M}_d^{s,t}(\T)$. 
To prove the assertion (ii), suppose that  
the irreducible representations 
$\rho=\rho_V \circ \varphi_s \circ \iota_t$ and 
$\rho'=\rho_{V'} \circ \varphi_{s'} \circ \iota_{t'}$ of $\A$ 
are isomorphic, 
where $V=V(\ell_1,a_1) \otimes \cdots \otimes V(\ell_n,a_n)$, 
$V'=V(\ell_1',a_1') \otimes \cdots \otimes V(\ell_{n'}',a_{n'}')$. 
Set $A=\rho(z)$, $A^*=\rho(z^*)$ and $B=\rho'(z)$, $B^*=\rho'(z^*)$. 
Then $A,\,A^*$ are a TD-pair belonging to $\mathcal {STD}_d^{(b,b^*)}$, 
where $b=st, \,b^*=st^{-1}$, $d=\ell_1+\cdots+\ell_n$, 
and $B,\,B^*$ are a TD-pair belonging to $\mathcal {STD}_{d'}^{(c,c^*)}$, 
where $c=s't', \,c^*=s't'^{-1}$, $d'=\ell'_1+\cdots+\ell'_{n'}$ 
(see Proposition \ref{prop: STDd(b,b*)}). Since $\rho,\,\rho'$ are 
isomorphic, the TD-pair $A,\,A^*$ is isomorphic to the TD-pair $B,\,B^*$, 
so we have $(b,b^*)\sim (c,c^*)$ in the sense of $(\ref{sim 1st kind})$, 
i.e., $(s, t) \approx (s', t')$ in the sense of $(\ref{approx 1st kind})$. 
Moreover by \cite[Proposition 7.15]{ITe}, the Drinfel'd polynomial 
$P_V(\lambda)$ of the $\T$-module $V$ via $\varphi_s$ coincides with 
the Drinfel'd polynomial $P_{V'}(\lambda)$ of the $\T$-module $V'$ via 
$\varphi_{s'}$. By Theorem \ref{thm: 1st kind} (i) and Lemma 
\ref{lemma: q-string}, the multi-sets ${\{ S(\ell_i,a_i)\}}^n_{i=1}$, 
${\{S(\ell_i',a_i')\}}^{n'}_{i=1}$ of $q$-strings are equivalent.  

\medskip
Conversely for the irrreducible representations 
$\rho=\rho_V \circ \varphi_s \circ \iota_t$ and 
$\rho'=\rho_{V'} \circ \varphi_{s'} \circ \iota_{t'}$ of $\A$ 
with $V=V(\ell_1,a_1) \otimes \cdots \otimes V(\ell_n,a_n)$, 
$V'=V(\ell_1',a_1') \otimes \cdots \otimes V(\ell_{n'}',a_{n'}')$, 
suppose that $(s, t) \approx (s', t')$ and 
the multi-sets ${\{ S(\ell_i,a_i)\}}^n_{i=1}$, 
${\{S(\ell_i',a_i')\}}^{n'}_{i=1}$ of $q$-strings are equivalent. 
Set $b=st, \,b^*=st^{-1}$, $c=s't', \,c^*=s't'^{-1}$
and $d=\ell_1+\cdots+\ell_n$, 
$d'=\ell'_1+\cdots+\ell'_{n'}$. Then $(b,b^*)\sim (c,c^*)$ and $d=d'$, 
so $\mathcal {STD}_d^{(b,b^*)}=\mathcal {STD}_{d'}^{(c,c^*)}$. 
Set $A=\rho(z)$, $A^*=\rho(z^*)$. 
Then $A,\,A^*$ is a TD-pair belonging to $\mathcal {STD}_d^{(b,b^*)}$, 
so it belongs to  $\mathcal {STD}_{d'}^{(c,c^*)}$: 
the difference is the orderings of the eigenspaces of $A,\,A^*$. 
Apply Proposition \ref{prop: STDd(b,b*)} to $\mathcal {STD}_{d'}^{(c,c^*)}$. 
Then there exists a unique representation $\rho''$ of $\T$ up to isomorphism 
belonging to $\mathcal {I}rr_{d'}^{s', t'}(\T)$ such that the TD-pair 
$B=\rho'' \circ \iota_{t'} (z)$, $B^*=\rho'' \circ \iota_{t'} (z^*)$ 
is isomorphic to $A,\,A^*$. 
By Theorem \ref{thm: 1st kind} (iii), we may assume 
$\rho''=\rho_{V''} \circ \varphi_{s'}$ for some 
$V''=V(\ell''_1,a''_1) \otimes \cdots \otimes V(\ell''_{n''},a''_{n''})$. 
Apparently, 
$\rho''\circ \iota_{t'}=\rho_{V''} \circ \varphi_{s'} \circ \iota_{t'}$ 
is isomorphic to 
$\rho=\rho_V \circ \varphi_s \circ \iota_t$ 
as representations of $\A$, since 
the TD-pair 
$B=\rho'' \circ \iota_{t'} (z)$, $B^*=\rho'' \circ \iota_{t'} (z^*)$ 
is isomophic to 
$A=\rho(z)$, $A^*=\rho(z^*)$. 
Then by what we have already proved in the 1st half of the proof, 
the multi-set ${\{S(\ell_i'',a_i'')\}}^{n''}_{i=1}$ of $q$-strings 
is equivalent to ${\{ S(\ell_i,a_i)\}}^n_{i=1}$ and hence  
to ${\{S(\ell_i',a_i')\}}^{n'}_{i=1}$. This means 
$\rho''\circ \iota_{t'}=\rho_{V''} \circ \varphi_{s'} \circ \iota_{t'}$ 
is isomorphic to 
$\rho'=\rho_{V'} \circ \varphi_{s'} \circ \iota_{t'}$
as representations of $\A$. 
So $\rho,\,\rho'$ are isomophic as representations of $\A$. 
This completes the proof of the theorem. 
\hfill $\Box $

\medskip
Next we consider the case $(\varepsilon,\varepsilon^*)=(1,0)$. 
Then $\varphi_s(\T) \subseteq \LL'$. 
Note that the subalgebra $\LL'$ of $\LL$ is, by the 
triangular decomposition of $\LL$, isomorphic to the 
algebra generated by the symbols 
$e_0^+, \,e_1^{\pm},\, k_i^{\pm}\, ~(i=0,1) $
subject to the defining relations
\begin{eqnarray*}
k_0 k_1&=&k_1 k_0=1, \\
k_i k_i^{-1}&=&k_i^{-1} k_i=1, \\
k_0 e_0^+ k_0^{-1}&=&q^2 e_0^+, \\
k_1 e_1^{\pm} k_1^{-1}&=&q^{\pm2} e_1^{\pm},\\
\lbrack e_0^+, e_1^- \rbrack &=&0, \\
\lbrack e_1^+,e_1^- \rbrack 
&=& \frac{k_1 - k_1^{-1}}{q-q^{-1}},
\end{eqnarray*}
\vspace{-9mm} 
\begin{eqnarray*}
\lbrack e_i^+, {(e_i^+)}^2 e_j^+ - (q^2 + q^{-2})
e_i^+ e_j^+e_i^+ +e_j^+{(e_i^+)}^2 \rbrack = 0 ~~~ (i \neq j).
\end{eqnarray*}
So the evaluation module $V(\ell, a)$  makes sense as an 
$\LL'$-modle even for $a=0$: 
for the standard basis 
$v_0,v_1, \cdots, v_d$ of $V(\ell,a)$,
\begin{eqnarray*}
k_0 v_i &=& q^{2i-\ell}\, v_i \\
k_1 v_i &=& q^{\ell-2i}\, v_i, \\
e_0^+ v_i &=& a\,q \,[i+1] \,v_{i+1}, \\
e_1^+ v_i &=& [\ell-i+1]\, v_{i-1}, \\
e_1^- v_i &=& [i+1]\, v_{i+1}.
\end{eqnarray*}
For a positive integer $\ell$  and a scalar $a \in \C$, 
allowing $a=0$, $V(\ell, a)$ is irreducible as an $\LL'$-module 
and called an {\it evaluation module for $\LL'$}.
Since the coproduct  $\Delta$ of $\LL$ is closed for $\LL'$, i.e., 
$\Delta (\LL') \subseteq \LL' \otimes \LL'$, the tensor product 
$V(\ell_1, a_1) \otimes \cdots \otimes V(\ell_n, a_n)$ 
of evaluation modules for $\LL'$ becomes an $\LL'$ -module. 
We denote by $V(\ell)$ the evaluaion module $V(\ell, 0)$. 
We allow $\ell=0$ for $V(\ell)$ and understand that $V(0)$ is the 
trivial $\LL'$-module, i.e., the 1-dimensional space 
on which $k_i^{\pm 1}=1$, the identity map, 
$e_0^+=e_1^{\pm}=0$, the zero map. 
Thus $V(\ell, a)$ means the evaluation module for $\LL'$ with 
$\ell \geq 1, a \neq 0 $ and $V(\ell)$ the evaluation module 
$V(\ell,0)$ for $\LL'$ with $\ell \geq 0$. 

\begin{theorem}
\label{thm: 2nd kind}
{\rm (Case ($\varepsilon,\varepsilon^*) = (1,0)$)}~~
Let $\T=\T_q^{(1,0)}$ denote the augmented TD-algebra of the 2nd kind. 
The following $(i)$, $(ii)$, $(iii)$ hold. 

\begin{enumerate}
\item[$(i)$] 
A tensor product 
$V(\ell)\otimes V(\ell_1,a_1)\otimes \cdots \otimes V(\ell_n,a_n)$ 
of evaluation modules for $\LL'$ is irreducible 
as a $\T$-module via $\varphi_s$ if and only if 
$-s^{-2} \notin S(\ell_i,a_i)$ 
for all $i ~(1 \leq i \leq n)$ and the multi-set 
${\{ S(\ell_i,a_i)\}}^n_{i=1} $
of $q$-strings is in general position. 
In this case, the $\T$-module 
$V=V(\ell)\otimes V(\ell_1,a_1)\otimes \cdots \otimes V(\ell_n,a_n)$ 
via $\varphi_s$ has type $s$ and diameter 
$d=\ell+\ell_1+\cdots+\ell_n$ and the Drinfel'd polynomial $P_V(\lambda)$ 
of the $\T$-module $V$ via $\varphi_s$ is 
$$ P_V(\lambda) = 
P_{V(\ell)}(\lambda)\prod^n_{i=1} P_{V(\ell_i,a_i)}(\lambda),$$
where 
\begin{eqnarray*}
P_{V(\ell)}(\lambda) &=& \lambda^{\ell}, \\
P_{V(\ell_i,a_i)}(\lambda)& =&
\prod_{c \in S(\ell_i,a_i)}(\lambda +c).
\end{eqnarray*}

\item[$(ii)$] 
Set $V=V(\ell)\otimes V(\ell_1,a_1) \otimes \cdots \otimes V(\ell_n,a_n)$, 
$V'= V(\ell')\otimes 
V(\ell_1',a_1') \otimes \cdots \otimes V(\ell_{n'}',a_{n'}')$ 
and assume that $V$, $V'$ are both irreducible 
as a $\T$-module via $\varphi_s$. 
Then $V$, $V'$ are isomorphic as $\T$-modules via $\varphi_s$ 
if and only if $\ell=\ell'$ and 
the multi-sets ${\{ S(\ell_i,a_i)\}}^n_{i=1}$, 
${\{S(\ell_i',a_i')\}}^{n'}_{i=1}$ of $q$-strings coincide, i.e., 
$n=n'$,\, $\ell_i=\ell_i', \, a_i=a_i'$ 
for all $i~(0 \leq i \leq n)$ 
with a suitable reordering of 
$S(\ell_1', a_1'), \cdots, S(\ell_n',a_n')$. 

\item[$(iii)$] 
Every nontrivial 
finite-dimensional irreducible $\T$-module of 
type $s$ is isomorphic to some $\T$-module 
$V(\ell)\otimes V(\ell_1,a_1) \otimes \cdots \otimes V(\ell_n,a_n)$ 
via $\varphi_s$ . 
\end{enumerate}
\end{theorem}

Theorem \ref{thm: 2nd kind} will be proved in Section 7. 
Note that the Drinfel'd polynomial of an irreducible $\T$-module 
$V(\ell)\otimes 
V(\ell_1,a_1) \otimes \cdots \otimes V(\ell_n,a_n)$ via $\varphi_s$ 
is determined by $\ell$ and the multi-set ${\{ S(\ell_i,a_i)\}}^n_{i=1}$ 
of $q$-strings, independent of $\varphi_s$. 
Problem \ref{problem: T}, which is to determine $\mathcal {M}^s_d(\T)$ 
and $\mathcal {M}_d^{s,t}(\T)$, is solved by Theorem \ref{thm: 2nd kind} 
as follows in the case of $(\varepsilon,\varepsilon^*) = (1,0)$. 
Assume $(\varepsilon,\varepsilon^*) = (1,0)$. 
The set $\mathcal {M}^s_d(\T)$ is determined in terms of 
tensor products of evaluation modules 
by Theorem \ref{thm: 2nd kind} (i), (ii), (iii). 
Recall the bijection 
$V \mapsto P_V (\lambda)$ 
from $\mathcal {M}^s_d(\T)$ to $\mathcal P^s_d$ 
in Theorem~$\ref{thm: sigma}'$. 
The subset $\mathcal {M}_d^{s,t}(\T)$ of $\mathcal {M}^s_d(\T)$ 
is nonempty if and only if the conditions (\ref{C1zt}), (\ref{C1z*t}) 
hold, i.e., 
\begin{eqnarray}
\label{C1 2nd kind}
\pm st \notin \{q^i~|~ 
i= -d+1,\, -d+2, \cdots, d-1 \},
\end{eqnarray}
and in this case $\mathcal {M}_d^{s,t}(\T)$ is mapped onto 
$\mathcal P^{s,t}_d$ by the bijection 
$V \mapsto P_V (\lambda)$ (see Corollary~\ref{cor: M}). 
For an irreducible $\T$-module 
$V=V(\ell)\otimes 
V(\ell_1,a_1) \otimes \cdots \otimes V(\ell_n,a_n)$ via $\varphi_s$, 
we find by Theorem~\ref{thm: 2nd kind} (i) that 
$P_V (\lambda)$ does not vanish at $\lambda =t^2$ 
if and only if $-t^2 \notin S(\ell_i,a_i)$ 
for all $i ~(1 \leq i \leq n)$. Thus we have  

\begin{corollary}
\label{cor: M 2nd kind}
Assume $(\varepsilon,\varepsilon^*) = (1,0)$. Then 
$\mathcal {M}^s_d(\T)$ and $\mathcal {M}_d^{s,t}(\T)$ 
are determined as follows. 

\begin{enumerate}
\item[$(i)$] 
$\mathcal {M}^s_d(\T)$ consists of the isomorphism classes 
of $\T$-modules  
$V(\ell)\otimes V(\ell_1,a_1) \otimes \cdots \otimes V(\ell_n,a_n)$ 
via $\varphi_s$ with the properties that  
\begin{description}
\item[$(i.1)$] 
the multi-set 
${\{ S(\ell_i,a_i)\}}^n_{i=1}$
of $q$-strings is in general position, 
\item[$(i.2)$] 
$-s^{-2} \notin S(\ell_i,a_i)$ 
for all $i ~(1 \leq i \leq n)$, 
\item[$(i.3)$] 
$d=\ell+\ell_1+\cdots+\ell_n$. 
\end{description}
The isomorphism classes of such $\T$-modules 
$V(\ell)\otimes V(\ell_1,a_1) \otimes \cdots \otimes V(\ell_n,a_n)$ 
via $\varphi_s$ are in one-to-one correspondence with 
the set of pairs of 
$\ell \in \N \cup \{0\}$ and 
the multi-sets ${\{ S(\ell_i,a_i)\}}^n_{i=1}$ of $q$-strings 
that have the properties $(i.1)$, $(i.2)$, $(i.3)$ above. 

\item[$(ii)$] 
$\mathcal {M}_d^{s,t}(\T)$ is nonempty if and only if 
the condition $(\ref{C1 2nd kind})$ holds. 
Suppose the condition $(\ref{C1 2nd kind})$ holds. Then 
$\mathcal {M}_d^{s,t}(\T)$ consists of the isomorphism classes 
of $\T$-modules 
$V(\ell)\otimes V(\ell_1,a_1) \otimes \cdots \otimes V(\ell_n,a_n)$ 
via $\varphi_s$ with the properties $(i.1)$, $(i.2)$, $(i.3)$ above 
and 
\begin{description}
\item[$(ii.1)$]
$-t^2 \notin S(\ell_i,a_i)$ 
for all $i ~(1 \leq i \leq n)$. 
\end{description}

\end{enumerate}
\end{corollary}

The next theorem follows from Corollary~\ref{cor: M 2nd kind} 
and \cite[Proposition 7.15]{ITe}. 
It solves Problem~\ref{problem: A}, 
which is to determine $\mathcal {I}rr(\A)$, the set of isomorphism classes 
of finite-dimensional irreducible representations of 
the TD-algebra $\A=\A_q^{(1,0)}$ of the 2nd kind 
that have the property ${\rm (C_1)}$. 
For an $\LL'$-module $V$, let $\rho_V$ denote the representation of $\LL'$ 
afforded by the $\LL'$-module $V$. Then $\rho_V \circ \varphi_s$ is 
the representation of $\T$ afforded by the $\T$-module $V$ via $\varphi_s$,  
and $\rho_V \circ \varphi_s \circ \iota_t$ is a representation of $\A$, 
where $ \iota_t : \A \longrightarrow \T ~~( z, \,z^* \mapsto z_t, \,z^*_t) $
is the injective algebra homomorphism from Proposition \ref{prop: iota}.

\begin{theorem}
\label{thm: A 2nd kind}
Assume $(\varepsilon,\varepsilon^*) = (1,0)$. Let $\A=\A_q^{(1,0)}$ 
denote the TD-algebra of the 2nd kind. 
The following $(i)$, $(ii)$, $(iii)$ hold. 
\begin{enumerate}
\item[$(i)$] 
For an $\LL'$-module 
$V=V(\ell)\otimes V(\ell_1,a_1) \otimes \cdots \otimes V(\ell_n,a_n)$ 
and nonzero $s,\,t \in \C$, the representaiton 
$\rho_V \circ \varphi_s \circ \iota_t$ of $\A$ is irreducible if and only if 
\begin{description}
\item[$(i.1)$] 
the multi-set 
${\{ S(\ell_i,a_i)\}}^n_{i=1}$
of $q$-strings is in general position, 
\item[$(i.2)$] 
none of 
$-s^{-2},\, -t^2$ belongs to $S(\ell_i,a_i)$ for any $i$ ~$(1 \leq i \leq n)$, 
\item[$(i.3)$] 
none of $\pm st$ equals $q^i$ 
for any $i \in \Z$ ~$(-d+1 \leq i \leq d-1)$. 
\end{description}
\item[$(ii)$] 
For $\LL'$-modules 
$V=V(\ell)\otimes V(\ell_1,a_1) \otimes \cdots \otimes V(\ell_n,a_n)$, 
$V'=V(\ell')\otimes 
V(\ell_1',a_1') \otimes \cdots \otimes V(\ell_{n'}',a_{n'}')$ 
and $(s,t),\, (s',t') \in (\C \backslash \{0\})\times (\C \backslash \{0\})$, 
set $\rho=\rho_V \circ \varphi_s \circ \iota_t$ and 
$\rho'=\rho_{V'} \circ \varphi_{s'} \circ \iota_{t'}$. 
Assume that the representations $\rho$, $\rho'$ of $\A$ are both irreducible. 
Then they are isomorphic as representations of $\A$ 
if and only if $\ell=\ell'$, the multi-sets  
${\{ S(\ell_i,a_i)\}}^n_{i=1}$,  ${\{ S(\ell_i',a_i')\}}^{n'}_{i=1}$ 
coincide and $(s,t) \approx (s',t')$ in the sense of 
$({\rm \ref{approx 2nd kind}})$ ,i.e., 
$$
(s', t') \in 
\{\pm (s,t),\, \pm (t^{-1},s^{-1})\}.
$$ 
\item[$(iii)$] 
Every nontrivial finite-dimensional irreducible representation of $\A$ is 
isomorphic to $\rho_V \circ \varphi_s \circ \iota_t$ 
for some $\LL'$-modules 
$V=V(\ell)\otimes V(\ell_1,a_1) \otimes \cdots \otimes V(\ell_n,a_n)$ 
and $(s,t) \in (\C \backslash \{0\})\times (\C \backslash \{0\})$. 
\end{enumerate}
\end{theorem}

We do not give a proof of Theorem \ref{thm: A 2nd kind}, 
since it can be proved by the same argument 
for the case of the $q$-Onsager algebra.

\medskip
Finally we consider the case $(\varepsilon,\varepsilon^*)=(0,0)$. 
By Proposition \ref{prop: phi}, 
$\varphi_s$ gives an isomorphism 
between the augmented TD-algebra $\T$ and the Borel subalgebra $\B$ of $\LL$. 
The TD-algebra $\A$ is isomorphic to  
the positive part of the Borel subalgebra $\B$. 

\begin{theorem}
\label{thm: 3rd kind}
{\rm (Case ($\varepsilon,\varepsilon^*) = (0,0)$)}~~
Let $\T=\T_q^{(0,0)}$ denote the augmented TD-algebra of the 3rd kind. 
The following $(i)$, $(ii)$, $(iii)$ hold. 

\begin{enumerate}
\item[$(i)$] 
A tensor product 
$V(\ell_1,a_1)\otimes \cdots \otimes V(\ell_n,a_n)$ 
of evaluation modules for $\LL$ is irreducible 
as a $\T$-module via $\varphi_s$ if and only if 
the multi-set 
${\{ S(\ell_i,a_i)\}}^n_{i=1} $
 of $q$-strings is in general position. 
In this case, the $\T$-module 
$V=V(\ell_1,a_1)\otimes \cdots \otimes V(\ell_n,a_n)$ 
via $\varphi_s$ has type $s$ and diameter 
$d=\ell_1+\cdots+\ell_n$ and the Drinfel'd polynomial $P_V(\lambda)$ 
of the $\T$-module $V$ via $\varphi_s$ is 
$$ P_V(\lambda) = 
\prod^n_{i=1} P_{V(\ell_i,a_i)}(\lambda),$$
where 
$$
P_{V(\ell_i,a_i)}(\lambda) =
\prod_{c \in S(\ell_i,a_i)}(\lambda +c).
$$

\item[$(ii)$] 
Set $V=V(\ell_1,a_1) \otimes \cdots \otimes V(\ell_n,a_n)$, 
$V'= V(\ell_1',a_1') \otimes \cdots \otimes V(\ell_{n'}',a_{n'}')$ 
and assume that $V$, $V'$ are both irreducible 
as a $\T$-module via $\varphi_s$. 
Then $V$, $V'$ are isomorphic as $\T$-modules via $\varphi_s$ 
if and only if 
the multi-sets ${\{ S(\ell_i,a_i)\}}^n_{i=1}$, 
${\{S(\ell_i',a_i')\}}^{n'}_{i=1}$ of $q$-strings coincide. 

\item[$(iii)$] 
Every nontrivial 
finite-dimensional irreducible $\T$-module of 
type $s$ is isomorphic to some $\T$-module 
$V(\ell_1,a_1) \otimes \cdots \otimes V(\ell_n,a_n)$ 
via $\varphi_s$ . 
\end{enumerate}
\end{theorem}

Theorem \ref{thm: 3rd kind} is well-known 
but a brief proof will be given in Section 7.  
The polynomial $\lambda^d P_V (\lambda^{-1})$~ 
$(d=\ell_1 + \cdots + \ell_n)$
for the case $(\varepsilon,\varepsilon^*)=(0,0)$ is known as 
the original Drinfel'd polynomial: 
$$ \lambda^d P_V (\lambda^{-1}) 
= \prod^n_{i=1} \prod_{c \in S (\ell_i,a_i)}(1+c \lambda).$$
Corollary \ref{cor: M 3rd kind} and Theorem \ref{thm: A 3rd kind} below, 
which are the main results of \cite[Theorem 1.6, Theorem 1.7]{ITd},  
follow immediately from 
Theorem \ref{thm: 3rd kind} through ${\rm Theorem} \ref{thm: sigma}'$ 
and Corollary \ref{cor: M}, solving Problem \ref{problem: A} and 
Problem \ref{problem: T} in the case of 
$(\varepsilon,\varepsilon^*) = (0,0)$. 

\begin{corollary}
\label{cor: M 3rd kind}
Assume $(\varepsilon,\varepsilon^*) = (0,0)$. Then 
$\mathcal {M}^s_d(\T)$ and $\mathcal {M}_d^{s,t}(\T)$ 
are determined as follows. 

\begin{enumerate}
\item[$(i)$] 
$\mathcal {M}^s_d(\T)$ consists of the isomorphism classes 
of $\T$-modules  
$V(\ell_1,a_1) \otimes \cdots \otimes V(\ell_n,a_n)$ 
via $\varphi_s$ with the properties that  
\begin{description}
\item[$(i.1)$] 
the multi-set 
${\{ S(\ell_i,a_i)\}}^n_{i=1}$
of $q$-strings is in general position, 
\item[$(i.2)$] 
$d=\ell_1+\cdots+\ell_n$. 
\end{description}
The isomorphism classes of such $\T$-modules 
$V(\ell_1,a_1) \otimes \cdots \otimes V(\ell_n,a_n)$ 
via $\varphi_s$ are in one-to-one correspondence with 
the set of 
the multi-sets ${\{ S(\ell_i,a_i)\}}^n_{i=1}$ of $q$-strings 
that have the properties $(i.1)$, $(i.2)$ above. 

\item[$(ii)$] 
$\mathcal {M}_d^{s,t}(\T)$ is nonempty for any nonzero $s,\, t \in \C$. 
$\mathcal {M}_d^{s,t}(\T)$ consists of the isomorphism classes 
of $\T$-modules 
$V(\ell_1,a_1) \otimes \cdots \otimes V(\ell_n,a_n)$ 
via $\varphi_s$ with the properties $(i.1)$, $(i.2)$ above 
and 
\begin{description}
\item[$(2.1)$]  
$-t^2 \notin S(\ell_i,a_i)$ 
for all $i ~(1 \leq i \leq n)$. 
\end{description}

\end{enumerate}
\end{corollary}

\begin{theorem}
\label{thm: A 3rd kind}
Assume $(\varepsilon,\varepsilon^*) = (0,0)$. Let $\A=\A_q^{(0,0)}$ 
denote the TD-algebra of the 3rd kind. 
The following $(i)$, $(ii)$, $(iii)$ hold. 
\begin{enumerate}
\item[$(i)$] 
For an $\LL$-module 
$V=V(\ell_1,a_1) \otimes \cdots \otimes V(\ell_n,a_n)$ 
and nonzero $s,\,t \in \C$, the representaiton 
$\rho_V \circ \varphi_s \circ \iota_t$ of $\A$ is irreducible if and only if 
\begin{description}
\item[$(i.1)$] 
the multi-set 
${\{ S(\ell_i,a_i)\}}^n_{i=1}$
of $q$-strings is in general position, 
\item[$(i.2)$] 
$-t^2 \notin S(\ell_i,a_i)$ for any $i$ ~$(0 \leq i \leq n)$. 
\end{description}

\item[$(ii)$] 
For $\LL$-modules 
$V=V(\ell_1,a_1) \otimes \cdots \otimes V(\ell_n,a_n)$, 
$V'=V(\ell_1',a_1') \otimes \cdots \otimes V(\ell_{n'}',a_{n'}')$ 
and $(s,t),\, (s',t') \in (\C \backslash \{0\})\times (\C \backslash \{0\})$, 
set $\rho=\rho_V \circ \varphi_s \circ \iota_t$ and 
$\rho'=\rho_{V'} \circ \varphi_{s'} \circ \iota_{t'}$. 
Assume that the representations $\rho$, $\rho'$ of $\A$ are both irreducible. 
Then they are isomorphic as representations of $\A$ 
if and only if the multi-sets  
${\{ S(\ell_i,a_i)\}}^n_{i=1}$,  ${\{ S(\ell_i',a_i')\}}^{n'}_{i=1}$ 
coincide and $(s,t)= \pm (s',t')$. 

\item[$(iii)$] 
Every nontrivial finite-dimensional irreducible representation of $\A$ is 
isomorphic to $\rho_V \circ \varphi_s \circ \iota_t$ 
for some $\LL$-module 
$V=V(\ell_1,a_1) \otimes \cdots \otimes V(\ell_n,a_n)$ 
and $(s,t) \in (\C \backslash \{0\})\times (\C \backslash \{0\})$. 
\end{enumerate}
\end{theorem}

\medskip
Let $A, A^* \in {\rm{End}}(V)$ be a TD-pair of Type I 
with eigenspaces 
${\{V_i\}}^d_{i=0},~{\{V_i^*\}}^d_{i=0}$ respectively. 
Then we have the split decomposition (see Section 1.1): 
\[ V= \bigoplus^d_{i=0} U_i , \]
where
\[ U_i=(V_0^* + \cdots + V_i^*) \cap (V_i + \cdots + V_d). \]
By \cite[Corollary 5.7 ]{ITT}, it holds that 
$$ {\rm dim}\, U_i= {\rm dim}\, V_i = {\rm dim}\, V_i^* 
~~~ (0 \leq i \leq d), $$
and 
$$ {\rm dim}\, U_i= {\rm dim}\, U_{d-i} ~~~~~~~~~~
~~~ (0 \leq i \leq d).$$
Note that ${\rm dim}\,U_i$ is invariant under standardization 
of $A,\,A^*$. We want to find the generating function for 
${\rm dim}\,U_i$:
$$ g(\lambda)= \sum^d_{i=0} \,({\rm dim}\, U_i ) \,\lambda^i. $$
We may assume that $A,A^*$ are standardized. 
Then by Theorem \ref{thm: q-Onsager}, Theorem \ref{thm: A 2nd kind}, 
Theorem \ref{thm: A 3rd kind}, 
the TD-pair $A,\,A^*$ is afforded via $\varphi_s \circ \iota_t$ 
by an $\LL$-module 
$$V=V(\ell_1,a_1) \otimes \cdots \otimes V(\ell_n,a_n)$$ 
for the cases $(\varepsilon,\varepsilon^*)=(1,1), \,(0,0)$
and by an $\LL'$-module 
$$V=V(\ell) \otimes V(\ell_1,a_1) \otimes \cdots \otimes V(\ell_n,a_n)$$ 
for the case $(\varepsilon,\varepsilon^*)=(1,0)$. 
The split decomposition of $V$ for $A,\,A^*$ coincides 
with the eigenspace decomposition of the element $k_0$ of $\LL$ 
acting on $V$. Thus we have 

\begin{proposition} 
\label{prop: ch}
{\rm (\cite[Conjecture 13.7 ]{ITT})}
\begin{eqnarray*}
g(\lambda)&=&\prod^n_{i=1}(1+\lambda+\lambda^2+ \cdots + \lambda^{\ell_i}) 
~~{\rm if} ~(\varepsilon,\varepsilon^*)=(1,1),\,(0,0), \\
g(\lambda)&=&\prod^n_{i=0}(1+\lambda+\lambda^2+ \cdots + \lambda^{\ell_i}) 
~~{\rm with ~\ell_0=\ell} 
~~{\rm if} ~(\varepsilon,\varepsilon^*)=(1,0).
\end{eqnarray*}
\end{proposition}

A TD-pair $A,A^*$ is called a {\it Leonard pair} 
if dim $U_i=1$ for all $i ~~(0 \leq i \leq d)$. 
A standardized TD-pair $A,A^*$ of Type I is a Leonard 
pair if and only if it is afforded by an evaluation 
module. In view of this fact, a standardized TD-pair $A,\,A^*$ 
of Type I is regarded as a `tensor product of Leonard pairs'.

\section{Linear bases for $\A$ and $\T$}  
In this section, we give a linear basis for the TD-algebra $\A$ that 
involves the generators $z, z^*$. We also give two linear bases for the 
augmented TD-algebra $\T$; one involves the generators $x,\,y,\,k^{\pm 1}$ and 
the other involves the 
generators $z_t, \,z_t^*,\, k^{\pm 1}$ (see Section 1.2). 
Using these bases, we prove Proposition \ref{prop: iota}, 
Proposition \ref{prop: C1} and Proposition \ref{prop: phi}. 

\medskip
For an integer $r \geq 0$, we denote by $\Lambda_r$ the set of sequences 
$\lambda=(\lambda_0, \lambda_1, \cdots, \lambda_r)$ of integers such that 
$\lambda_0 \geq 0, \,\lambda_i \geq 1 ~(1 \leq i \leq r)$, 
and define $\Lambda$ to be the union of $\Lambda_r ~(r \geq 0)$: 
\begin{eqnarray*}
\Lambda_r &=& \{ \lambda = (\lambda_0, \lambda_1, \cdots, \lambda_r) \in 
\Z^{r+1}~|~ \lambda_0 \geq 0,\, \lambda_i \geq 1 ~(1 \leq i \leq r) \}, \\
\Lambda &=& \bigcup_{r \in \N \cup \{ 0 \}} \Lambda_r. 
\end{eqnarray*}
Call $\lambda = (\lambda_0, \lambda_1, \cdots, \lambda_r) \in \Lambda$ 
{\it irreducible} if there exists an integer $i ~(0 \leq 1 \leq r)$ such that 
$$ \lambda_0 < \lambda_1 < \cdots < \lambda_i \geq \lambda_{i+1} \geq \cdots 
\geq \lambda_r. $$
Note that each $\lambda$ in $\Lambda_0 \cup \Lambda_1$ is irreducible. 
We denote the set of irreducible $\lambda \in \Lambda$ by $\Lambda^{irr}$: 
$$ \Lambda^{irr} = \{\lambda \in \Lambda ~|~ \lambda ~\mbox{is irreducible} \}.$$
\medskip
Let $X,\,Y$ denote noncommuting indeterminates. 
For $\lambda = (\lambda_0, \lambda_1, \cdots, \lambda_r) \in \Lambda$, 
we define the word $\omega_{\lambda}(X,Y)$ by 
\begin{eqnarray*}
\omega_{\lambda}(X,Y) = \Biggl\{
\begin{array}{lll}
X^{\lambda_0} Y^{\lambda_1} \cdots X^{\lambda_r}
&& \mbox{if $r$ is even}, \\
X^{\lambda_0} Y^{\lambda_1} \cdots Y^{\lambda_r}
&& \mbox{if $r$ is odd}, 
\end{array} 
\end{eqnarray*}
where we interpret $X^{\lambda_0}=1$ if $\lambda_0=0$. By the length of the 
word $\omega_{\lambda}(X,Y)$, we mean 
$\lambda_0 + \lambda_1 + \cdots + \lambda_r$ and denote it by $|\lambda|$: 
$$ |\lambda| = \lambda_0 + \lambda_1 + \cdots + \lambda_r.$$ 

\begin{theorem}
\label{thm: A basis}
The following set is a basis of the TD-algebra $\A$ as a $\C$-vector 
space: 
$$ \{ \omega_{\lambda}(z,z^*)~ |~ \lambda \in \Lambda^{irr} \}. $$
\end{theorem}

\begin{theorem}
\label{thm: T basis}
Each of the following sets is a basis of 
the augmented TD-algebra $\T$ as a $\C$-vector space: 
\begin{enumerate}
\item[$(i)$] 
$\{ k^n\, \omega_{\lambda} (x,y) ~|~ n \in \Z, ~\lambda \in \Lambda^{irr} \}.$ 
\item[$(ii)$] 
$\{ k^n\, \omega_{\lambda} (z_t,z_t^*) ~|~ 
n \in \Z, \,\lambda \in \Lambda^{irr} \},$ 
where $t$ is a fixed nonzero scalar of $\C$ and 
$z_t,z_t^*, k^{\pm 1}$ are the second generators of $\T$ that are introduced 
in $(\ref{zt})$, $(\ref{z*t})$ in Section 1.2.
\end{enumerate}
\end{theorem}

We first prove the spanning property for the sets in 
Theorem \ref{thm: A basis}, Theorem \ref{thm: T basis}. 
Our strategy will be to reduce the essential part to 
\cite[Theorem 2.29]{ITa}.
We start with a description of the $\C$-algebra generated by symbols 
$\xi,\, \eta, \,\kappa,\, \kappa^{-1}$ subject to the relations 
${\rm {(TD)}_0'}$: 
$\kappa \kappa^{-1} = \kappa^{-1} \kappa =1, \,
\kappa \xi \kappa^{-1} = q^2 \xi, \,
\kappa \eta \kappa^{-1} = q^{-2} \eta$. 
Let $\Phi$ denote the free algebra over $\C$ 
generated by symbols $\xi,\, \eta$.
Let $\C [\kappa, \kappa^{-1}]$ denote the algebra over 
$\C$ generated by symbols 
$\kappa,\, \kappa^{-1}$ subject to the relations 
$\kappa \kappa^{-1} = \kappa^{-1}\kappa =1$. 
Consider the $\C$-vector space $\C [\kappa, \kappa^{-1}] \otimes \Phi$, 
where $\otimes = \otimes_{\C}$. 
This space has the basis 
$$ \{ \kappa^n \otimes \omega_{\lambda} (\xi,\eta) ~| ~n \in \Z, \,\lambda \in 
\Lambda \}. $$
Define the product of basis elements by 
$$ (\kappa^m \otimes \omega_{\lambda} (\xi,\eta) ) 
( \kappa^n \otimes \omega_{\mu} (\xi,\eta)) 
= \kappa^{m+n} \otimes \omega_{\lambda} (q^{-2n} \xi, q^{2n} \eta) 
\omega_{\mu} (\xi,\eta)$$
and extend it bilinearly to the product of elements of 
$\C [\kappa, \kappa^{-1}] \otimes \Phi$. 
Then $\C [\kappa, \kappa^{-1}] \otimes \Phi$ becomes an associative 
$\C$-algebra. The mapping $f \otimes u \mapsto fu ~
(f \in \C [\kappa, \kappa^{-1}],\,u \in \Phi)$ 
induces a $\C$-algebra isomorphism from
 $\C [\kappa, \kappa^{-1}] \otimes \Phi$ to the 
$\C$-algebra generated by $\xi, \,\eta, \,\kappa,\, \kappa^{-1}$
subject to the relations ${\rm {(TD)}_0'}$: 
$\kappa \kappa^{-1} = \kappa^{-1} \kappa = 1, \,
\kappa \xi \kappa^{-1} = q^2 \xi,\, \kappa \eta \kappa^{-1} = q^{-2} \eta$. 
We henceforth identify these two algebras via the isomorphism and denote this 
algebra by $\C [\kappa, \kappa^{-1}] \Phi$.

\medskip
Define the elements $v_0, v_1 \in \Phi$ and $u_0, u_1 
\in \C [\kappa, \kappa^{-1}] \Phi$ by 
\begin{eqnarray*}
v_0 &=& [ \xi, \,\xi^2 \eta - \beta \xi \eta \xi + \eta \xi^2], \\
v_1 &=& [ \xi \eta^2 - \beta \eta \xi \eta  + \eta^2 \xi,\, \eta ],\\
u_0 &=& \delta'\,(\varepsilon^* \xi^2 k^2 - \varepsilon k^{-2} \xi^2),\\
u_1 &=& \delta'\,(\varepsilon^* k^2 \eta^2 - \varepsilon \eta^2 k^{-2}),
\end{eqnarray*}
where $\beta = q^2 + q^{-2}$ and 
$\delta' = - (q-q^{-1})(q^2-q^{-2})(q^3-q^{-3}) q^4$. 
Let $\I$ denote the two-sided ideal of 
$\C [\kappa, \kappa^{-1}] \Phi$ generated by $v_0 -u_0, \,v_1-u_1$: 
$$ {\I} = \C [\kappa, \kappa^{-1}]\Phi\, (v_0 -u_0)\Phi + 
\C [\kappa, \kappa^{-1}]\Phi(v_1-u_1)\Phi. $$
Since the relations ${\rm (TD)'}$ for the augmented TD-algebra is 
$v_0=u_0,\,v_1=u_1$, 
the quotient algebra $\C [\kappa, \kappa^{-1}] \Phi/ {\I}$ coincides 
with $\T$ and we have the canonical algebra homomorphism 
$$ \pi_{\T}: \C [\kappa, \kappa^{-1}]\Phi \longrightarrow \T ~~~
(\xi, \eta, \kappa, \kappa^{-1} \mapsto x,y,k,k^{-1} ~\mbox{respectively}). $$
%

\medskip
Let $\J$ denote the two-sided ideal of $\Phi$ generated by $v_0,v_1$: 
$$ \J = \Phi v_0 \Phi + \Phi v_1 \Phi. $$
Write $\A_{\rm III}=\Phi/\J$ for the quotient algebra 
(the TD-algebra of the 3rd kind), 
and let us use the bar notation for the canonical algebra homomorphism: 
$$ \pi_{\A_{\rm III}} : \Phi \longrightarrow \A_{\rm III} ~~~ 
(\xi, \eta \mapsto \bar{\xi}, \bar{\eta} ~~ \mbox{respectively}). $$
By \cite[Theorem 2.29]{ITa}, 
the set 
$\{ \omega_{\lambda}(\bar{\xi}, \bar{\eta})~|~ 
\lambda \in \Lambda^{irr} \}$ is a basis for $\A_{\rm III}$. 
Consequently 
$$ \Phi = \W + \J ~~~ \mbox{(direct sum)}, $$
where $\W$ is the subspace of $\Phi$ spanned by 
$$ \{ \omega_{\lambda} (\xi, \eta) ~|~ \lambda \in \Lambda^{irr} \}. $$
%

\medskip
For an integer $n \geq 0$, 
we mean by a word of length $n$ in $\Phi$ a product 
$a_1 a_2 \cdots a_n$ in $\Phi$ such that $a_i \in \{ \xi, \eta \}$ for 
$1 \leq i \leq n$. 
We interpret the word of length $0$ as the identity in $\Phi$. 
Let $\Phi_n$ denote the subspace of $\Phi$ spanned by the words of length $n$. 
For example, $\Phi_0 = \C \,1$. 
We have the direct sum 
$\Phi = \sum_{n \geq 0} \Phi_n$ and 
$\Phi_r \,\Phi_s = \Phi_{r+s}$ for all $r, s \geq 0$. 
For an integer $n \geq 0$, difine 
$\W_n = \Phi_n \cap \W$ and $\J_n = \Phi_n \cap \J$. 
This yields the direct sum 
decompositions 
$\W = \sum_{n \geq 0} \W_n$ and $\J = \sum_{n \geq 0} \J_n$. 
By $\Phi = \W + \J$, we have 
$$ \Phi_n = \W_n + \J_n $$
for $n \geq 0$. Since $v_0, v_1 \in \Phi_4$,
$$ \J_n = \sum \Phi_i v_0 \Phi_j + \sum \Phi_i v_1 \Phi_j, $$
where both sums are over the ordered pairs of nonnegative integers $(i,j)$ 
such that $i+j=n-4$. In particular, $\J_n=0$ for $n \leq 3$. 
Since $v_0=(v_0-u_0)+u_0, \,v_1=(v_1-u_1)+u_1$ and $u_0,u_1 \in 
\C [\kappa, \kappa^{-1}] \Phi_2$,
the above expression for $\J_n$ together with the definition of $\I$ 
implies 
$$ \J_n \subseteq {\I} + \C [\kappa, \kappa^{-1}] \Phi_{n-2} ~~~ (n \geq 4). $$
%

\medskip
To prove that the set in Theorem \ref{thm: T basis} (i) spans $\T$, 
it suffices to show 
$$ \C [\kappa, \kappa^{-1}] \Phi = 
\C [\kappa, \kappa^{-1}] \W + \I. $$
To this end we show 
$\C [\kappa, \kappa^{-1}] \Phi_n \subseteq \C [\kappa, \kappa^{-1}] \W + \I$ 
for $n \geq 0$ and this will be done by induction on $n$. 
Let $n$ be given. Recall $\Phi_n = \W_n + \J_n$. 
If $n \leq 3$, then $\J_n=0$, and so $\Phi_n=\W_n \subseteq \W$ and certainly 
$\C [\kappa, \kappa^{-1}] \Phi_n \subseteq \C [\kappa, \kappa^{-1}] \W + \I$ 
as desired. 
If $n \geq 4$, we argue by $\J_n \subseteq \I + 
\C [\kappa, \kappa^{-1}] \Phi_{n-2}$
\begin{eqnarray*}
\C [\kappa, \kappa^{-1}] \Phi_n &=& \C [\kappa, \kappa^{-1}] \W_n 
+ \C [\kappa, \kappa^{-1}] \J_n \\
&\subseteq& \C [\kappa, \kappa^{-1}] \W + \I + \C [\kappa, \kappa^{-1}] \Phi_{n-2}
\end{eqnarray*}
and this is contained in $\C [\kappa, \kappa^{-1}] \W + \I$ 
by induction on n. We have now proved that the set 
in Theorem \ref{thm: T basis} (i) spans $\T$.

\medskip
To prove the spanning property for the sets in 
Theorem \ref{thm: A basis} and Theorem \ref{thm: T basis} (ii), 
let $\J'$ denote the two-sided ideal of $\Phi$ generated by 
$v_0 - \varepsilon \delta [\xi,\eta]$ and 
$v_1 - \varepsilon^* \delta [\xi,\eta]$ : 
$$ \J' = \Phi(v_0 - \varepsilon \delta [\xi,\eta]) \Phi + 
\Phi (v_1 - \varepsilon^* \delta [\xi,\eta]) \Phi, $$
where $\delta = -{(q^2 - q^{-2})}^2$. Since $\J_n = \sum \Phi_i v_0 \Phi_j + 
\sum \Phi_i v_i \Phi_j$ 
over $(i,j)$ with $i+j=n-4$, we have 
$$ \J_n \subseteq \J' + \Phi_{n-2} ~~~~~ (n \geq 4), $$
noting that $[\xi,\eta] \in \Phi_2$. We claim that 
$$ \Phi = \W + \J'. $$
The inclusion $\supseteq$ is from the construction. 
To get the inclusion $\subseteq$, 
we show $\Phi_n \subseteq \W + \J'$ for $n \geq 0$ and 
this will be done by induction on $n$. Let $n$ be given. 
If $n \leq 3$, then $\J_n=0,\, 
\Phi_n=\W_n + \J_n = \W_n \subseteq \W$ so 
$\Phi_n \subseteq \W + \J'$ as desired. If $n \geq 4$, we argue 
by $\J_n \subseteq \J' + \Phi_{n-2}$
\begin{eqnarray*}
\Phi_n &=& \W_n + \J_n \\
&\subseteq& \W + \J' + \Phi_{n-2}
\end{eqnarray*}
and this is contained in $\W + \J'$ by induction on  $n$. 
We have now proved the claim. 
Write $\A' = \Phi/ \J'$ for the quotient algebra, and let us use the prime 
notation for the canonical algebra homomorphism : 
$$ \pi_{\A'} : \Phi \longrightarrow \A' ~~~ (\xi,\eta \mapsto \xi',\eta' ~~
\mbox{respectively}). $$
The above claim implies that $\A'$ is spanned by 
$$ \{ \omega_{\lambda} (\xi',\eta')~|~\lambda \in \Lambda^{irr} \}.$$
Since the defining relations (TD) for $\A$ are 
$v_0 = \varepsilon \delta [\xi,\eta],\,v_1 = \varepsilon^* \delta [\xi,\eta]$, 
$\A'$ is isomorphic to the TD-algebra 
$\A$ by the correspandence $\xi' \mapsto z,\, \eta' \mapsto z^*$. 
This proves that the set in Thereom \ref{thm: A basis} spans $\A$. 
For a fixed nonzero $t \in \C$, 
let $\langle z,z_t^*\rangle$ denote the subalgebra of $\T$ 
generated by $z_t,z_t^*$. 
We have a surjective algebra homomorphism 
$$ \A' \longrightarrow \langle z_t,z_t^*\rangle
 ~~~ (\xi',\eta' \mapsto z_t,z_t^* ~~\mbox{respectively}) $$
by the relations (TD) for $z_t,z_t^*$. 
So $\langle z_t,z_t^*\rangle$ is spaned by 
$\{ \omega_{\lambda}(z_t,z_t^*)~|~ \lambda \in \Lambda^{irr}\}$. 
By the relations ${({\rm TD})}_0$ for $\T$, it holds that 
$\T = \sum_{n \in \Z} k^n \langle z_t,z_t^*\rangle$. 
Therefore the spanning property holds for the set 
in Theorem \ref{thm: T basis} (ii). 

\medskip
Next we prove the linear independency of the sets in 
Theorem \ref{thm: A basis} and Theorem \ref{thm: T basis}. 
For a nonzero $s \in \C$, let $\varphi_s$ be the algebra homomorphism 
from $\T$ to the $U_q({sl}_2)$-loop algebra $\LL = U_q (L({sl}_2))$ 
as in Proposition \ref{prop: phi}: 
$\varphi_s$ sends $x,\,y,\,k$ to $x(s),\,y(s),\,k(s)$, respectively, 
where 
\begin{eqnarray*}
x(s) & =& \alpha(se_0^+ + \varepsilon s^{-1} e_1^- k_1) 
~~~ {\rm with} ~\alpha = - q^{-1} {(q-q^{-1})}^2,\\
y(s) & =& \varepsilon^* s e_0^- k_0 + s^{-1} e_1^+, \\
k(s) & =& s k_0. 
\end{eqnarray*}
Note that the existence of the algebra homomorphism $\varphi_s$ has been 
established already, although the injectivety of $\varphi_s$ is left to 
be proved. 
For a nonzero $t \in \C$, let $\iota_t$ be the algebra homomorphism from $\A$ 
to $\T$ as in Proposition \ref{prop: iota}: 
$\iota_t$ sends $z,\,z^*$ to 
$z_t=x + tk + \varepsilon t^{-1} k^{-1}$, 
$z_t^* = y + \varepsilon^* t^{-1} k + t k^{-1}$, respectively. 
Note also that the existence of the algebra 
homomorphism $\iota_t$ has been estabished already, 
although the injectivity of $\iota_t$ is left to be proved. 
We set $z_t(s) = \varphi_s \circ \iota_t(z), \,
z_t^*(s) = \varphi_s \circ \iota_t (z^*)$: 
\begin{eqnarray*}
z_t(s) &=& x(s) + t k(s) + \varepsilon t^{-1} {k(s)}^{-1}, \\
z_t^*(s) &=& y(s) + \varepsilon^* t^{-1} k(s) + t {k(s)}^{-1}.
\end{eqnarray*}
\begin{lemma}
\label{lemma: lin indep}
For nonzero scalars $s,t \in \C$, each of the following sets is linearly 
independent in $\LL$. 

\begin{enumerate}
\item[$(i)$] 
$\{ {k(s)}^n \omega_{\lambda}(x(s),y(s))~|~ 
n \in \Z, \, \lambda \in \Lambda^{irr} \}.$ 
\item[$(ii)$] 
$\{ {k(s)}^n \omega_{\lambda}(z_t(s),z_t^*(s))~|~ 
n \in \Z, \, \lambda \in \Lambda^{irr} \}.$ 
\end{enumerate}
\end{lemma}

The linear independency of the sets in 
Theorem \ref{thm: T basis} (resp. Theorem \ref{thm: A basis}) immediately 
follows from Lemma \ref{lemma: lin indep} 
by applying the algebra homomorphism $\varphi_s$ 
(resp. $\varphi_s \circ \iota_t$). 
We prove Lemma \ref{lemma: lin indep}  
by using the triangular decomposition of $\LL$ 
together with the basis of $\A_{\rm III}$ given in \cite{ITa}. 
Let $\langle e_0^+,e_1^+\rangle$ (resp. $\langle e_0^-,e_1^-\rangle$ ) 
be the subalgebra of $\LL$ genereted by 
$e_0^+, \, e_1^+$ (resp. $e_0^-, \, e_1^-$). 
Then by \cite[Theorem 2.29]{ITa}, 
$\langle e_0^+,e_1^+\rangle$ (resp. $\langle e_0^-,e_1^-\rangle$ ) 
is isomorphic to $\A_{\rm III}$ and has the set $B^+$ (resp.$B^-$) 
as a linear basis, 
where 
\begin{eqnarray*}
B^+ &=& \{ \omega_{\lambda} (e_0^+,e_1^+)~|~ \lambda \in \Lambda^{irr} \}, \\
B^- &=& \{ \omega_{\lambda} (e_0^-,e_1^-)~|~ \lambda \in \Lambda^{irr} \}.
\end{eqnarray*}
By the triangular decomposition of $\LL$, the set 
$$ B = \{ \omega^- k_0^n \omega^+ ~|~ n \in \Z, \, \omega^- \in B^-, \, 
\omega^+ \in B^+ \}$$
is a linear basis of $\LL$, 
and so every element of $\LL$ is uniquely expressed as 
a finite sum of 
$$ c_{\mu, n, \lambda} 
\omega_{\mu} (e_0^-,e_1^-) k^n \omega_{\lambda} (e_0^+,e_1^+)$$
with $c_{\mu, n, \lambda} \in \C, \, n \in \Z, \, 
\mu, \lambda \in \Lambda^{irr}$. 
The expression for the element ${k(s)}^n \omega_{\lambda} (x(s),y(s))$ in 
question of Lemma \ref{lemma: lin indep} (i) is, 
by the defining relations of $\LL$, 
$$ s^n k_0^n \omega_{\lambda} (\alpha s e_0^+,s^{-1}e_1^+) $$
plus some other terms 
$c_{\mu', n'  \lambda'} \omega_{\mu'} (e_0^-,e_1^-) k^{n'} \omega_{\lambda'} 
(e_0^+,e_1^+)$ with $|\lambda'|<|\lambda|$. The highest term 
$s^n k_0^n \omega_{\lambda} (\alpha s e_0^+,s^{-1}e_1^+)$ is 
the product of the nonzero scalar 
$s^n \omega_{\lambda}(\alpha s,s^{-1}) \in \C$, $k_0^n $ 
and the element 
$\omega_{\lambda}(e_0^+,e_1^+) \in B^+$. 
Therefore any linear dependency relation among the elements 
in Lemma \ref{lemma: lin indep} (i) 
is the trivial one by induction on the maximal length $|\lambda|$ of $\lambda$ 
that appears in the relation. Similarly the set in (ii) is shown to be 
linearly independent. 
This completes the proof of  
Lemma \ref{lemma: lin indep}.  

\bigskip
\noindent
{\bf Proof of Proposition \ref{prop: iota} and 
Proposition \ref{prop: phi}: the injectivity of $\iota_t$ and $\varphi_s$}. 
The algebra homomorphism $\varphi_s \circ \iota_t$ is 
injective by Theorem \ref{thm: A basis} and 
Lemma \ref{lemma: lin indep} (ii) and so $\iota_t$ is injective. 
The injectivity of $\varphi_s$ follows 
from Theorem \ref{thm: T basis} (i) and 
Lemma \ref{lemma: lin indep} (i). 
\hfill $\Box $ 

\bigskip
\noindent
{\bf Proof of Proposition \ref{prop: C1}}. 
Let $V$ be a finite-dimensional irreducible $\T$-module of type $s$, 
diameter $d$ and 
$V=\bigoplus_{i=0}^d U_i$ the weight space decomposition from 
Lemma \ref{lemma: weight-space decomposition}. 
Since $z_t=x+tk+\varepsilon t^{-1}k^{-1}$, $k|_{U_i}=sq^{2i-d}$, 
we have $z_t=x+\theta_i$ on $U_i$, where 
$\theta_i=stq^{2i-d} + \varepsilon s^{-1} t^{-1} q^{d-2i}$. 
Since $xU_i \subseteq U_{i+1}$, we have 
$(z_t - \theta_0)(z_t - \theta_1) \cdots (z_t - \theta_d) = 0$ on $V$. 
If $\theta_0,\cdots,\theta_d$ are mutually distinct, 
then $z_t$ is diagonalizable on V and it holds that 
$$V_i+V_{i+1}+\cdots+V_d=U_i+U_{i+1}+\cdots+U_d~~~(0 \leq i \leq d),$$
where $V_i$ is the eigenspace of $z_t$ on $V$ that belongs to 
the eigenvalue $\theta_i$. 

\medskip
Conversely, suppose $z_t$ is diagonalizable on $V$. 
Let $\theta_{i_{0}},\theta_{i_{1}},\cdots,\theta_{i_{r}}$ denote the distinct 
members among $\theta_i~(0 \leq i \leq d)$. 
Then $(z_t - \theta_{i_{r}}) \cdots 
(z_t - \theta_{i_{1}})(z_t - \theta_{i_{0}})$ vanishes on $V$, 
in particular on $U_0$. We claim 
$$ (z_t - \theta_{i_{j}}) \cdots (z_t - \theta_{i_{1}})(z_t - \theta_{i_{0}}) 
= f_j (x) ~~~ \mbox{on} ~ U_0 $$
for some monic polynomial $f_j$ of degree $j+1$ $(0 \leq j \leq r)$. 
The claim holds for $j=0$, 
since $z_t - \theta_{i_{0}} = x + \theta_0 - \theta_{i_{0}}$ on $U_0$. 
If the claim holds for $j$, then there exit scalars 
$c_0, c_1, \cdots, c_{j+1}$ with $c_{j+1}=1$ such that 
$$ (z_t - \theta_{i_{j}}) \cdots (z_t - \theta_{i_{1}})(z_t - \theta_{i_{0}}) 
u = \sum_{i=0}^{j+1} c_i x^i u ~~~~~ (u \in U_0) .$$
Since the right-hand side has the $i$-th term $c_ix^i u \in U_i$ and 
$z_t-\theta_{i_{j+1}}=x+\theta_i-\theta_{i_{j+1}}$ on $U_i$, 
the claim holds for $j+1$. 
Thus the claim is proved by induction on $j$. 
Since 
$(z_t - \theta_{i_{r}}) \cdots (z_t - \theta_{i_{1}})(z_t - \theta_{i_{0}})$ 
vanishes on $U_0$, 
the monic polynomial $f_r$ of degree $r+1$ satisfies $f_r(x) U_0=0$. 
This implies $x^{r+1} U_0=0$, since $x^i U_0 \subseteq U_i$ and 
$V$ is the direct sum of $U_i's$. 
On the other hand, we have $V=\T U_0$ by the irreducibility of the 
$\T$-module $V$, so $V$ is spanned by 
$\omega_{\lambda}(x,y) U_0 ~(\lambda \in \Lambda^{irr})$ due to 
Theorem \ref{thm: T basis}. 
For $\lambda=(\lambda_0,\lambda_1,\cdots,\lambda_n) \in \Lambda^{irr}$, 
there exists some $i ~(0 \leq i \leq n)$ such that
$\lambda_0 < \lambda_1 < \cdots < \lambda_i \geq \lambda_{i+1} \geq \cdots \geq \lambda_n$. 
If $\omega_{\lambda} (x,y) U_0 \neq 0$ for such $\lambda$, 
then $i,\,n$ are even and 
it holds that 
$\lambda_{i+1}=\lambda_{i+2} \geq \cdots \geq \lambda_{n-1}=\lambda_n$, 
since 
$x U_j \subseteq U_{j+1}, \, y U_j \subseteq U_{j-1}$ with $U_{-1}=0$. 
Moreover we have $\lambda_i \leq r$, otherwise $\omega_{\lambda}(x,y)U_0=0$ 
by the vanishing property $x^{r+1} U_0=0$ we just proved. 
Therefore if $\omega_{\lambda}(x,y) U_0 \neq 0$, then 
$\omega_{\lambda}(x,y) U_0 \subseteq U_j$, where 
$j= \lambda_i - \lambda_{i-1} + \cdots + \lambda_2 - \lambda_1 + \lambda_0 
\leq \lambda_i \leq r$. 
Thus $V=\T U_0 \subseteq U_0 + U_1 + \cdots + U_r$. 
This implies $r=d$, i.e., $\theta_0, \cdots, \theta_d$ are mutually distinct. 
We have now prove the first half (i) of Proposition \ref{prop: C1}. 
The sencond half (ii) is similarly proved, using $V=\T U_d$. 
\hfill $\Box $

\section{The subspace of height 0 in $\T$}
Let $\T$ be the augmented TD-algebra.  
$\T$ is the algebra generated by $x,\, y,\, k^{\pm}$ 
subject to the relations $({\rm TD})_0'$, $({\rm TD})'$ in Section 1.2. 
We introduce the notion of height 
for a word in $x,\, y,\, k^{\pm}$ 
and discuss the structure of the subspace of $\T$ 
spanned by the words of height 0. 
The main result of this section is Theorem \ref{thm: comm}. 
As applications of Theorem \ref{thm: comm}, 
we prove 
Theorem \ref{thm: shape} and 
the injectivity of $\sigma$ in Theorem \ref{thm: sigma}. 
We keep the notations in Section 2. 

\medskip
Consider the free algebra over $\C$ generated by 
$\xi, \eta, \kappa, \kappa^{-1}$.  
Let $\tau_0$ denote the automorphism of this free algebra that sends 
$\xi, \eta, \kappa,\kappa^{-1}$ to $\eta, \xi, \kappa, \kappa^{-1}$ 
respectively, 
and let $\tau_1$ denote the anti-automorphism that reverses the word order. 
Then $\tau_0$, $\tau_1$ commute and the product 
$\tau=\tau_0 \tau_1=\tau_1 \tau_0$ is an antiautomorphism that sends 
a word $\zeta_1 \zeta_2 \cdots \zeta_n$ 
to $\zeta_n' \cdots \zeta_2' \zeta_1'$ 
($\zeta_i \in \{\xi, \eta, \kappa, \kappa^{-1}\}$), 
where $\zeta_i'=\eta,\, \xi,\, \kappa,\, \kappa^{-1}$ 
for $\zeta_i=\xi,\, \eta,\, \kappa,\, \kappa^{-1}$ respectively. 
Note that $\tau_0^2=\tau_1^2=\tau^2=id$, the identity map. 
Keeping the notations in Section 2, 
let $\Phi$ denote the free algebra generated by $\xi, \eta$, 
and $\C[\kappa, \kappa^{-1}] \Phi$ 
the algebra generated by $\xi, \eta, \kappa, \kappa^{-1}$ 
subject to the relations 
$({\rm TD})_0': \kappa \kappa^{-1} = \kappa^{-1} \kappa =1,\, 
\kappa \xi \kappa^{-1} = q^2 \xi,\, 
\kappa \eta \kappa^{-1} = q^{-2} \eta$. 
Since $({\rm TD})_0'$ is invariant under $\tau$ as a set of relations, 
the map $\tau$ induces 
an anti-antomorphism of the algebra $\C [\kappa, \kappa^{-1}] \Phi$. 
Recall the elements $v_0, v_1 \in \Phi$ and $u_0, u_1 \in \C [k,k^{-1}] \Phi$ 
introduced in Section 2: 
\begin{eqnarray*}
v_0&=&\lbrack \xi,\, \xi^2\eta-\beta \xi\eta\xi + \eta\xi^2\rbrack,\\
v_1&=&\lbrack \xi\eta^{2}-\beta \eta\xi\eta + \eta^{2}\xi,\, \eta\rbrack ,\\
u_0&=&\delta' 
(\varepsilon^*\xi^2\kappa^2-\varepsilon\kappa^{-2}\xi^2),\\
u_1&=&\delta' 
(\varepsilon^*\kappa^2\eta^2-\varepsilon\eta^2\kappa^{-2}), 
\end{eqnarray*} 
where $\beta = q^2 + q^{-2}$, 
$\delta' = - (q-q)(q^{-2} - q^2)(q^{-3} -q^3) q^4$. 
The augmented TD-algebra $\T$ is defined by 
$({\rm TD})':\, v_0=u_0,\, v_1=u_1$ 
together with ${({\rm TD})}_0'$. 
Since $v_0^{\tau}=v_1, u_0^{\tau}=u_1$, 
$({\rm TD})'$ is invariant under $\tau$ and 
the map $\tau$ induces an anti-antomorphism of $\T$. 
Also $\tau$ induces an anti-antomorphism of $\A_{\rm III} = \Phi/\J$, 
where $\J$ is the  two-sided ideal of $\Phi$ generated by $v_0, v_1$. 
We use the same notation $\tau$ 
for these anti-antomorphisms of 
$\C [\kappa, \kappa^{-1}] \Phi,\, \T,\, \A_{\rm III}$.

\medskip
Let $W$ denote the free semi-group generated by $\xi,~\eta$. 
As a set, $W$ is the collection of all words in $\Phi$. 
Let $h$ 
\[
h:~W\longrightarrow \Z
\]
denote the semi-group homomorphism from $W$ to 
the additive group $\Z$ 
defined by $h(\xi)=1,~h(\eta)=-1$. For a word $w \in W$, 
the value $h(w)$ is called the \textit{height} of $w$. Thus a word 
of height $0$ is a word in which $\xi,~\eta$ appear the same number 
of times. Denote by $\Phi^{(i)}$ the subspace of $\Phi$ 
linearly spanned by the words of height $i$: 
\[
\Phi^{(i)}=Span\{w\in W~|~h(w)=i\}. 
\]
Then $\Phi$ is the direct sum of the vector spaces 
$\Phi^{(i)}~(i\in \Z$): 
\begin{eqnarray}
\label{height grading Phi}
\Phi=\bigoplus_{i\in \Z}\Phi^{(i)}. 
\end{eqnarray}
The above decomposion is an algebra grading, i.e., 
$\Phi^{(i)}\Phi^{(j)} \subseteq \Phi^{(i+j)}$. 
Note that $\Phi^{(0)}$ is a subalgebra of $\Phi$. 
The antiautomorphism $\tau$ changes the sign of the height of a word 
and so sends $\Phi^{(i)}$ to $\Phi^{(-i)}$. In particular, $\tau$ 
induces an antiautomorphism of the subalgebra $\Phi^{(0)}$. 
Let $\Phi^{sym}$ denote the subspace of $\Phi^{(0)}$ consisting of 
the fixed points of $\tau$: 
\[
\Phi^{sym}=\{v\in\Phi^{(0)}~|~v^{\tau}=v\}. 
\]
A word $w\in W$ is called \textit{nil} if $w$ can be written as 
$w=w_1w_2$ with $w_1, w_2 \in W$ and $h(w_2) < 0$. Let $\Phi^{nil}$ 
denote the subspace of $\Phi^{(0)}$ 
linearly spanned by the words of height $0$ that are nil: 
\[
\Phi^{nil}=Span \{w \in W~|~h(w)=0,~w~ \textrm{is nil}\}. 
\]
Then $\Phi^{nil}$ is a two-sided ieal of $\Phi^{(0)}$ and invariant 
under the antiautomorphism $\tau$. 
Recall $\Phi_n$ is the subspase of $\Phi$ spanned by the words of 
length $n$ in $\xi, \eta$. 
Set $\Phi_n^{sym}=\Phi_n\cap \Phi^{sym}$ and 
$\Phi_n^{nil}=\Phi_n\cap \Phi^{nil}$. Then we have the direct sum 
decompositions as vector spaces: 
\begin{eqnarray}
\label{decomp sym}
\Phi^{sym}&=&\bigoplus_{n\geq 0}\Phi_n^{sym},\\
\label{decomp nil}
\Phi^{nil}&=&\bigoplus_{n\geq 0}\Phi_n^{nil}.
\end{eqnarray}

\medskip 
The algebra $\C[\kappa,\kappa^{-1}]\Phi$ becomes a graded algebra 
\begin{eqnarray*}
\C[\kappa,\kappa^{-1}]\Phi=
\bigoplus_{i\in \Z}\C[\kappa,\kappa^{-1}]\Phi^{(i)}. 
\end{eqnarray*}
Recall $\T=\C[\kappa,\kappa^{-1}]\Phi/\I$, 
where ${\mathcal I}$ is the two-sided ideal of 
$\C[\kappa,\kappa^{-1}]\Phi$ 
generated by $v_0-u_0,~ v_1-u_1$. 
Note that $v_0-u_0,~ v_1-u_1$ belong to 
$\C[\kappa,\kappa^{-1}]\Phi^{(2)}$, 
$\C[\kappa,\kappa^{-1}]\Phi^{(-2)}$ respectively. 
Set 
\[
{\mathcal I}^{(i)}={\mathcal I} \cap \C[\kappa,\kappa^{-1}]\Phi^{(i)}. 
\]
Then we have 
\begin{eqnarray}
\label{I}
{\mathcal I}=\bigoplus_{i\in \Z}{\mathcal I}^{(i)}.
\end{eqnarray}
For $\mathcal{T}=\C[\kappa,\kappa]\Phi/{\mathcal I}$, 
consider the canonical homomorphism 
\begin{eqnarray}
\label{can hom}
\pi=\pi_{_{\mathcal{T}}}:~ 
\C[\kappa,\kappa^{-1}]\Phi \longrightarrow 
\mathcal{T}~~
(\xi, \eta, \kappa, \kappa^{-1} 
\mapsto x, y, k, k^{-1}~ \textrm{respectively}).  
\end{eqnarray}
Set $\Psi=\pi(\Phi)$, $\Psi^{(i)}=\pi(\Phi^{(i)})$. 
Then by (\ref{I}), 
the algebra $\mathcal{T}$ inherits the algebra grading of 
$\C[\kappa,\kappa^{-1}]\Phi=
\bigoplus_{i\in \Z}\C[\kappa,\kappa^{-1}]\Phi^{(i)}$ 
via $\pi$: 
\[
\mathcal{T}=\C[k, k^{-1}]\Psi=
\bigoplus_{i\in \Z}\C[k, k^{-1}]\Psi^{(i)}. 
\]
This enables us to define the {\it height} function for $\mathcal{T}$: 
a nonzero element of $\mathcal{T}$ is said to have height $i$ 
if it belongs to $\C[k, k^{-1}]\Psi^{(i)}$. 

\medskip
Note that $\Psi=\pi(\Phi)$ is the subalgebra of 
$\mathcal{T}$ generated by $x, y$.  
$\Psi$ has the grading 
\begin{eqnarray}
\label{height grading Psi}
\Psi=\bigoplus_{i\in \Z}\Psi^{(i)}.
\end{eqnarray}
$\Psi^{(i)}$ is the subspace of $\Psi$ spanned by 
the words in $x,\, y$ of height $i$. 
$\Psi^{(0)}$ is a subalgebra of $\Psi$.  
The antiautomorphism $\tau$ of $\mathcal{T}$ sends 
$\Psi^{(i)}$ to $\Psi^{(-i)}$. In particular, $\tau$ induces 
an antiautomorphism of the subalgebra $\Psi^{(0)}$. 
Set 
\[
\Psi^{sym}=\pi(\Phi^{sym}). 
\]
Then $\Psi^{sym}\subseteq \Psi^{(0)}$ and every element of 
$\Psi^{sym}$ is fixed by $\tau$.  
Let $\Psi^{nil}$ denote the image of $\Phi^{nil}$ under $\pi$: 
\[
\Psi^{nil}=\pi(\Phi^{nil}). 
\]
Then $\Psi^{nil}$ is a two-sided ideal of $\Psi^{(0)}$ and 
invariant under $\tau$. 
Note that 
$k, k^{-1}$ commute with every element of $\Psi^{(0)}$. 
So 
$\C[k, k^{-1}]\Psi^{(0)}$ is a subalgebra of $\mathcal{T}$ 
and 
$\C[k, k^{-1}]\Psi^{sym}$ (resp. $\C[k, k^{-1}]\Psi^{nil}$) 
is a subspace (resp. two-sided ideal) of $\C[k, k^{-1}]\Psi^{(0)}$.



\begin{theorem}
\label{thm: comm}
The following $(i)$, $(ii)$ hold.
\begin{enumerate}
\item[$(i)$] 
$\C[k, k^{-1}]\Psi^{(0)}=
\C[k, k^{-1}]\Psi^{sym}+\C[k, k^{-1}]\Psi^{nil}$. 
\item[$(ii)$] 
The quotient algebra 
$\C[k, k^{-1}]\Psi^{(0)}/\C[k, k^{-1}]\Psi^{nil}$ 
is commutative and generated by $k,\, k^{-1}$ and 
$y^ix^i~(i=0,1,2,\ldots)$~
mod $\C[k, k^{-1}]\Psi^{nil}$. 
\end{enumerate}
\end{theorem}

\noindent
Proof. Our strategy will be to reduce
the essential part to 
\cite[Theorem 2.20]{ITa}. Recall the canonical homomophism 
in Section 2 
\[
\pi_{_{\mathcal{A}_{\rm III}}}:~
\Phi \longrightarrow\ \mathcal{A}_{\rm III}=\Phi/\mathcal{J }
\quad (\xi, \eta \mapsto \bar{\xi}, \bar{\eta}~
\textrm{respectively}), 
\]
where $\mathcal{J }$ is the two-sided ideal of $\Phi$ 
generated by $v_0,v_1$. 
Apply $\pi_{_{\mathcal{A}_{\rm III}}}$ to 
$\Phi^{(0)},~\Phi^{sym},~\Phi^{nil}$ and denote the images 
by $\mathcal{A}^{(0)},~\mathcal{A}^{sym},~ 
\mathcal{A}^{nil}$ respectively. 
Then $\mathcal{A}^{nil}$ is 
a two-sided ideal of $\mathcal{A}^{(0)}$ and 
the quotient $\mathcal{A}^{(0)}/\mathcal{A}^{nil}$ is 
a commutatve algebra generated by 
$\bar{\eta}^i\bar{\xi}^i~(i=0,1,2,\ldots)$ 
mod $\mathcal{A}^{nil}$ 
(see \cite[Lemma 3.1]{ITa}). So each element of 
$\mathcal{A}^{(0)}/\mathcal{A}^{nil}$ is a linear combination of 
$(\bar{\eta}^{i_1}\bar{\xi}^{i_1})
(\bar{\eta}^{i_2}\bar{\xi}^{i_2})
\cdots(\bar{\eta}^{i_n}\bar{\xi}^{i_n})$ 
mod $\mathcal{A}^{nil}$. 
Apply the antiautomorphism $\tau$ of $\A_{\rm III}$ to 
$\bar{w}=(\bar{\eta}^{i_1}\bar{\xi}^{i_1})
(\bar{\eta}^{i_2}\bar{\xi}^{i_2})
\cdots(\bar{\eta}^{i_n}\bar{\xi}^{i_n})$. 
Then by the commutativity of 
$\mathcal{A}^{(0)}/\mathcal{A}^{nil}$, we have 
\begin{eqnarray*}
\bar{w}-\bar{w}^{\tau}
&=&(\bar{\eta}^{i_1}\bar{\xi}^{i_1})
(\bar{\eta}^{i_2}\bar{\xi}^{i_2})
\cdots
(\bar{\eta}^{i_n}\bar{\xi}^{i_n}) -
(\bar{\eta}^{i_n}\bar{\xi}^{i_n})
\cdots
(\bar{\eta}^{i_2}\bar{\xi}^{i_2})
(\bar{\eta}^{i_1}\bar{\xi}^{i_1})\\
&\equiv& 0~~mod~\mathcal{A}^{nil}. 
\end{eqnarray*}
Thus 
\[
\bar{w}=\frac{1}{2}(\bar{w}+\bar{w}^{\tau})+
\frac{1}{2}(\bar{w}-\bar{w}^{\tau})
\in \mathcal{A}^{sym}+\mathcal{A}^{nil}, 
\]
and hence 
\[
\mathcal{A}^{(0)}=\mathcal{A}^{sym}+\mathcal{A}^{nil}.
\]
This means that for a word $w$ in $\xi, \eta$ 
of length $n$, height $0$, there exist elements 
$v\in \Phi^{sym},~v'\in \Phi^{nil}$ such that $w - v - v' \in \J$. 
Write $v$ (resp. $v'$) in the form of 
the direct sum (\ref{decomp sym}) (resp. (\ref{decomp nil})): 
$v=\sum_i v_i,\, v'=\sum_i v'_i$ 
with $v_i \in \Phi_i^{sym},\,  v'_i \in \Phi_i^{nil}$.  
Observe $\J=\bigoplus_i \J_i$, where $\J_i=\J \cap \Phi_i$. 
Since $w \in \Phi_n$, we have $w - v_n - v_n' \in \J_n$. 
Thus from the beginning, we may assume $v=v_n,\, v'=v'_n$, 
i.e., for a word $w$ in $\xi,\, \eta$ of length $n$, height $0$, 
there exist elements 
$v\in \Phi_n^{sym},~v'\in \Phi_n^{nil}$ such that
\begin{eqnarray}
\label{w-v-v'}
w-v-v' \in \mathcal{J }_n=\mathcal{J } \cap \Phi_n. 
\end{eqnarray}

\medskip
First we prove Theorem \ref{thm: comm} (i). 
Take any word $\hat w$ in $x,\,y$ of height 0, length $n$. 
Choose a word $w$ in $\xi,\, \eta$ of height 0, length $n$ 
such that $\hat w=\pi(w)$, 
where $\pi$ is the canonical homomorphism from (\ref{can hom}). 
Then by (\ref{w-v-v'}) there exixt elements 
$v\in \Phi_n^{sym},~v'\in \Phi_n^{nil}$ such that
$w-v-v' \in \mathcal{J }_n$. 
Observe $\J_n=\sum \Phi_iv_0\Phi_j+\sum \Phi_iv_1\Phi_j$, 
where the summation is over $i$, $j$ with $i+j=n-4$, 
since $v_0,\, v_1$ have length 4. 
Write the element $w-v-v'\in \mathcal{J }_n$ 
as a linear combination of $w_iv_0w_j,\, w'_iv_1w'_j$ 
for finitely many words $w_i, w_j, w_i', w_j'$ in $\xi, \eta$ 
such that 
$\ell(w_i)+\ell(w_j)=\ell(w'_i)+\ell(w_j')=n-4$,
where the functions $\ell$ stands for the length of a word.  
Recall that 
$\Phi$ is a graded algebra according 
to the height as in (\ref{height grading Phi}).
Since the element $w-v-v'$ has height 0 
and $v_0,\, v_1$ have height $2,\,-2$ respectively,  
we may assume that 
$h(w_i)+h(w_j)=-2,~h(w'_i)+h(w_j')=2$, 
where the function $h$ stands for the height of a word. 
Apply the canonical homomorphism $\pi$ from (\ref{can hom}) 
to $w-v-v'$. 
Then $\pi(v) \in \Psi^{sym},~
\pi(v')\in \Psi^{nil}$. Since $\pi(v_0)=\pi(u_0),~ 
\pi(v_1)=\pi(u_1)$, the terms $w_iv_0w_j,\, w'_iv_1w'_j$ 
in the linear combination for $w-v-v'$ are mapped to 
\begin{eqnarray*}
\pi(w_i)\pi(u_0)\pi(w_j) \in 
\C[k,k^{-1}](\Psi^{(0)}\cap \Psi_{n-2}),\\
\pi(w'_i)\pi(u_1)\pi(w'_j) \in 
\C[k,k^{-1}](\Psi^{(0)}\cap \Psi_{n-2}),
\end{eqnarray*}
where $\Psi_m=\pi(\Phi_m)$. 
Thus $\pi(w)-\pi(v)-\pi(v')$ belongs to 
$\C[k,k^{-1}](\Psi^{(0)}\cap \Psi_{n-2})$, 
while  $\pi(v)+\pi(v')$ belongs to $\Psi^{sym}+\Psi^{nil}$. 
The proof of part (i) is 
completed by induction on $n$. 

\medskip
Next we prove Theorem \ref{thm: comm} (ii). 
By Theorem \ref{thm: T basis}, 
$\C[k, k^{-1}]\Psi^{(0)}/\C[k, k^{-1}]\Psi^{nil}$ 
is linearly spanned by 
$k^n w_{\lambda}(x,y)$ 
mod $\C[k, k^{-1}]\Psi^{nil}$, 
where $n$ runs through $\Z$ and 
$\lambda$ runs through irreducible sequences 
such that the word $w_{\lambda}(x,y)$ has height $0$. 
Since $\lambda=(\lambda_0, \lambda_1, \cdots, \lambda_r)$ 
is irreducible and 
$w_{\lambda}(x,y)$ has height $0$, 
we may assume that $r$ is even and $\lambda_0=0$, 
$\lambda_1=\lambda_2\geq \lambda_3=\lambda_4 
\geq \cdots \geq \lambda_{r-1}=\lambda_r$, 
otherwise $w_{\lambda}(x,y)\in \Psi^{nil}.$ 
Therefore 
$\C[k, k^{-1}]\Psi^{(0)}/\C[k, k^{-1}]\Psi^{nil}$ 
is generated by 
$k^{\pm 1}$ and $y^ix^i~(i=0,1,2,\ldots)$ 
mod $\C[k, k^{-1}]\Psi^{nil}$. 
Note that $k^{\pm 1}$ commutes 
with $y^ix^i$. We want to show 
$y^ix^i,\, y^jx^j$ commute mod $\C[k, k^{-1}]\Psi^{nil}$. 
Set $w=(y^ix^i)(y^jx^j)$. By part (i) we just proved, 
there exist $f,\, g \in \C[k, k^{-1}],~ 
u \in \Psi^{sym},~v \in \Psi^{nil}$ such that 
$w=fu+gv$. Then $w^{\tau}=fu+gv^{\tau}$, since $u^{\tau}=u$ 
and $k, k^{-1}$ commute with every word of height $0$ in $x, y$. 
Note that $\Psi^{nil}$ is invariant under $\tau$, 
so $w-w^{\tau}=g(v-v^{\tau}) \in \C[k, k^{-1}]\Psi^{nil}$. 
Since $w-w^{\tau}=(y^ix^i)(y^jx^j)-(y^jx^j)(y^ix^i)$, 
this means $y^ix^i,~y^jx^j$ commute 
mod $\C[k, k^{-1}]\Psi^{nil}$ and the proof of part (ii) 
is completed.
\hfill $\Box $

\bigskip
\noindent
{\bf Proof of Theorem \ref{thm: shape}}. 
Let $V$ be a finite-dimensional irreducible module 
of the augmented TD-algebra ${\mathcal T}$. 
Let $V=\bigoplus_{i=0}^{d} U_i$ denote the weight-space decomposition 
of the $\T$-module $V$. 
We want to show 
${\rm dim}\,U_i\,\leq \,
{d \choose i}~ (0\leq i \leq d)$.

\medskip
Recall the algebra grading
$\mathcal{T}=\bigoplus_{i\in \Z} \C[k, k^{-1}]\Psi^{(i)}$, 
where $\Psi^{(i)}$ is the linear span of 
the words of height $i$ in $x,\, y$. 
Also recall $x U_i \subseteq U_{i+1},\, y U_i \subseteq U_{i-1}$.
The subalgebra $\C[k, k^{-1}]\Psi^{(0)}$ acts on $U_0$ 
and $\C[k, k^{-1}]\Psi^{nil}$ belongs to the kernel of the action. 
By Theorem \ref{thm: comm} (ii), there exists a common 
eigenvector $v\in U_0$ of $y^ix^i~(0\leq i\leq d)$. 
Since $y^jx^j$ vanishes on $U_0$ for $j\geq d+1$, each element of 
$\C[k, k^{-1}]\Psi^{(0)}$ fixes the 1-dimensional subspace $\C v$ 
by Theorem~\ref{thm: comm}~(ii). 
Since $V$ is irreducible and 
$\T=\bigoplus_{i\in \Z}\C[k, k^{-1}]\Psi^{(i)}$,  
we have 
$V=\mathcal{T}v=\sum_{i=0}^d\Psi^{(i)}v$. 
Then $U_i=\Psi^{(i)} v$, 
since $\Psi^{(i)}v\subseteq U_i$ and 
the sum $V=\sum_{i=0}^d U_i$ is direct.  
In particular, $U_0=\Psi^{(0)}v=\C v$. 
By Theorem \ref{thm: T basis}, 
\begin{eqnarray}
\label{Ui}
U_i=\Psi^{(i)}v=\sum_{\lambda \in \Lambda^{(i)}}
\C\, \omega_{\lambda}(x, y)\, v, 
\end{eqnarray}
where $\Lambda^{(i)}$ denotes the set of 
$\lambda=(\lambda_0,\lambda_1,\cdots,\lambda_r)\in \Lambda^{irr}$ 
such that $r$ is even and 
\[
\sum_{j=0}^{r}(-1)^j\lambda_j=i, \quad
\lambda_0 < \lambda_1 < \cdots < \lambda_r \leq d. 
\]
Since $\Lambda^{(i)}$ contains exactly 
${d \choose i}$ members, 
the proof of Theorem \ref{thm: shape} is completed. 
\hfill $\Box $

\bigskip
\noindent
{\bf Proof of Theorem \ref{thm: sigma}: the injectivity of $\sigma$}. 
Let $V$ be a finite-dimensional irreducible module 
of the augmented TD-algebra ${\mathcal T}$. 
Let $V=\bigoplus_{i=0}^{d} U_i$ denote the weight-space decomposition 
of the $\T$-module $V$. 
Recall $kv=sq^{2i-d} v$ for $v \in U_i$, 
$x U_i \subseteq U_{i+1},\, y U_i \subseteq U_{i-1}$. 
By Theorem \ref{thm: shape} we just proved, dim $U_0=1$. 
Let $\sigma_i=\sigma_i(V)$ denote the 
eigenvalue of $y^i x^i$ on the heighest weight space $U_0$. 
Apparently $\sigma_0=1, \sigma_i=0$ for $i \geq d+1$.

\medskip
We want to show $\sigma_d \neq 0$. 
By (\ref{Ui}) in the proof of Theorem \ref{thm: shape}, it holds that  
$U_d = \sum_{\lambda \in \Lambda^{(d)}} \omega_{\lambda} (x,y) U_0$. 
Since $\Lambda^{(d)} = \{ \lambda = (\lambda_0)~|~ \lambda_0 = d \}$, 
we have $U_d = x^d U_0$. 
In the proof of Theorem \ref{thm: shape}, 
the formula (\ref{Ui}) follows from 
$V=\bigoplus_{i=0}^d \Psi^{(i)} U_0$. 
Apply the same argument starting with
$V=\bigoplus_{i=0}^d \Psi^{(-i)} U_d$.  
Then we end up with 
$U_0 = \sum_{\lambda \in \Lambda^{(d)}} \omega_{\lambda} (y,x) U_d$.
Thus we have $U_0 = y^d U_d$. 
So $U_0 = y^d U_d =y^d x^d U_0$ and 
the eigenvalue $\sigma_d$ of $y^d x^d$ on $U_0$ is nonzero. 
Thus the diameter $d$ of the $\T$-module $V$ is determined 
by the property $\sigma_d \neq 0, \sigma_i=0 ~(i \geq d+1)$ of the sequence 
${\{ \sigma_i \}}_{i=0}^{\infty}$.

\medskip
Next we want to show that the isomorphism class of the $\T$-module $V$ is 
determined by the type $s$ and the sequence ${\{ \sigma_i \}}_{i=0}^d$.
Let $\mathcal{N}$ 
denote the set of elements of $\mathcal{T}$ 
that vanish on $U_0$: 
\[
\mathcal{N}=\{\nu \in \mathcal{T}~|~\nu U_0=0\}.
\]
Then $\mathcal{N}$ 
is a maximal left ideal of $\mathcal{T}$ and 
$V$ is isomorphic to $\mathcal{T}/\mathcal{N}$ as 
$\mathcal{T}$-modules. Hence it is enough to show that 
$\mathcal{N}$ is determined by $s$ and $\{\sigma_i\}_{i=0}^d$. 
With respect to the algebra grading 
$\T=\bigoplus_{i \in \Z} \C [k,k^{-1}] \Psi^{(i)}$, 
write $\nu \in \mathcal{N}$ as 
$\nu =  \sum \nu_i ~(\nu_i \in \C [k, k^{-1}] \Psi^{(i)})$. 
Then $\nu_i U_0 \subseteq U_i$. 
Since 
$V=\bigoplus_{i=0}^d U_i$ and $\nu U_0=0$, 
we have $\nu_i U_0 =0$, i.e., 
$\nu_i \in \mathcal{N}$. Therefore 
$$\mathcal{N} = \bigoplus_{i \in \Z} \mathcal{N}^{(i)},$$
where $\mathcal{N}^{(i)} = 
\mathcal{N} \cap \C [k, k^{-1}] \Psi^{(i)}$. 
Note that 
$\mathcal{N}^{(i)} = \C [\kappa, \kappa^{-1}] \Psi^{(i)}$ for $i < 0$. 
Thus it is enough to show that $\mathcal{N}^{(i)}$ is determined 
by $s$ and $\{\sigma_j\}_{j=0}^{d}$ for $i=0, 1, 2, \cdots$.

\medskip
For $i=0$, $\mathcal{N}^{(0)}$ is the kennel of the action of 
$\C[\kappa, \kappa^{-1}] \Psi^{(0)}$ on $U_0$. 
By Theorem \ref{thm: comm} (ii), 
$\C[k, k^{-1}]\Psi^{(0)}/\C[k, k^{-1}]\Psi^{nil}$ 
is generated by 
$k^{\pm 1}$ and $y^ix^i~ (i=0, 1, 2, \ldots)$ mod $\C[k, k^{-1}]\Psi^{nil}$. 
Apparently 
$\C[k, k^{-1}]\Psi^{nil}$ belongs to $\mathcal{N}^{(0)}$ 
and the action of $y^ix^i$ on $U_0$ is determined by $\sigma_i$.
Also, using the fact that the ${\mathcal T}$-module $V$ is type $s$,
the action of $k^{\pm 1}$ on $U_0$ is determined by $s$ and $d$. 
Therefore 
the action of $\C[k, k^{-1}]\Psi^{(0)}$ 
on $U_0$ is determined by $s$ and $\{\sigma_j\}_{j=0}^d$.  
Since $\mathcal{N}^{(0)}=\mathcal{N}\cap \C[k, k^{-1}]\Psi^{(0)}$ is 
the kernel of the action, $\mathcal{N}^{(0)}$  
is determined by $s$ and $\{\sigma_j\}_{j=0}^d$. 

\medskip
For $i \geq 1$, we claim 
\[
\mathcal{N}^{(i)}=\{\nu \in \C[k, k^{-1}]\Psi^{(i)}~|~
\Psi^{(-i)}\nu \subseteq \mathcal{N}^{(0)}\}. 
\]
For $\nu \in 
\mathcal{N}^{(i)} = \mathcal{N} \cap \C[k, k^{-1}] \Psi^{(i)}$, 
we have 
$\nu U_0 = 0$ and so $\Psi^{(-i)} \nu U_0 = 0$, i.e., 
$\Psi^{(-i)} \nu \subseteq \mathcal{N} 
\cap \C[k, k^{-1}] \Psi^{(0)} = \mathcal{N}^{(0)}$. 
Conversely, choose 
$\nu \in \C[k, k^{-1}]\Psi^{(i)}$ such that 
$\Psi^{(-i)}\nu \subseteq \mathcal{N}^{(0)}$. 
If $\nu U_0\neq 0$, 
then $\mathcal{T}\nu U_0=V$ by the irreducibility 
of the $\T$-module $V$. 
Since 
$\T \nu = \bigoplus_{j \in \Z} \C [k, k^{-1}] \Psi^{(j)} \nu$ and 
$\Psi^{(j)} \nu U_0 \subseteq U_{j+i}$, we have 
$\T \nu U_0 =\bigoplus_{j=-i}^{d-i} \Psi^{(j)} \nu U_0$, 
in particular 
$\Psi^{(-i)}\nu U_0=U_0$, 
which contradicts the assumption 
$\Psi^{(-i)}\nu \subseteq \mathcal{N}^{(0)}$. 
Thus 
$\nu U_0= 0$, i.e., 
$\nu \in \mathcal{N}\cap \C[k, k^{-1}]\Psi^{(i)}=\mathcal{N}^{(i)}$, 
and the claim is proved. 
This means 
$\mathcal{N}^{(i)}$ is determined by $\mathcal{N}^{(0)}$. 
Since $\mathcal{N}^{(0)}$ is determined by $s$ and 
$\{\sigma_j\}_{j=0}^{d}$, so is   $\mathcal{N}^{(i)}$.   
This completes the proof of the injectivity of $\sigma$ 
in Theorem \ref{thm: sigma}. 
\hfill $\Box $

\section{Finite-dimensional irreducible $\A$-modules via $\iota_t$: \\
Proof of Theorem \ref{thm: C2}}

The TD-algebra $\A=\A_q^{(\varepsilon,\varepsilon^*)}$ is 
by Proposition \ref{prop: iota}  
embedded into the augmented TD-algebra 
$\T=\T_q^{(\varepsilon,\varepsilon^*)}$ 
via the injective algebra-homomorphism 
$$ \iota_t : \A \longrightarrow \T 
~~(z \mapsto z_t,\, z^* \mapsto z_t^*)$$
for each fixed $t \in \C ~(t\neq 0)$, where
\begin{eqnarray*}
z_t &=& x + tk + \varepsilon t^{-1} k^{-1},\\
{z_t}^* &=& y + \varepsilon^* t^{-1} k + t k^{-1}.
\end{eqnarray*}
Let $V$ be a fintite-dimensional irreducible $\T$-module of type $s$ 
and diameter $d$. 
As we discussed in Section 1.2, 
the pair $A=z_t|_V,\, A^*={z_t}^*|_V$ 
of linear transformations of $V$ gives rise to a TD-pair 
if and only if 
\begin{quote}
${({\rm C}_1)}_t$:  ~~the action of $z_t,\, z_t^*$ on $V$ 
are both diagonalizable,\\
${({\rm C}_2)}_t$:  ~
\begin{minipage}[t]{13cm}
$V$ is irreducivble as an $\langle z_t, z^*_t\rangle$-module, \\
where $\langle z_t,z^*_t\rangle$ is the subalgebra of $\T$ 
genarated by $z_t, \,z^*_t$.
\end{minipage}
\end{quote}
By Proposition \ref{prop: C1}, the condition 
$({\rm C}_1)_t$ holds if and only if 
$\theta_i \neq \theta_j$ and $\theta_i^* \neq \theta_j$ 
for $i \neq j$ $(0 \leq i,\, j \leq d)$, where 
\begin{eqnarray*}
\theta_i &=& s\, t\, q^{2i-d} + \varepsilon\, s^{-1}t^{-1} q^{d-2i},\\
{\theta_i}^* &=& \varepsilon^* s\, t^{-1}q^{2i-d} + s^{-1}t\, q^{d-2i}.
\end{eqnarray*}
In this section, we prove Theorem \ref{thm: C2}, 
a criterion for $({\rm C}_2)_t$. 
Namely 
assume $({\rm C}_1)_t$. Then the condition $({\rm C}_2)_t$ 
holds if and only if $P_V(t^2 + \varepsilon\, \varepsilon^* t^{-2}) \neq 0$, 
where $P_V (\lambda)$ is the Drinfel'd polynomial of the $\T$-module $V$.

\medskip
We proceed parallel to \cite{ITd}. 
Let $V=\bigoplus_{i=0}^d U_i$ 
denote the weight-space decomposition of the $\T$-module $V$, 
and $F_i$ 
the projection of $V=\bigoplus_{i=0}^d U_i$ onto $U_i$. 
Note that $k$ acts on $V$ as $\sum_{t=0}^d s\, q^{2i-d} F_i$. 
Identifying $z,\, z^*$ with $z_t,\, z_t^*$ via $\iota_t$, 
we write $z=z_t,\, z^*={z_t}^*$ for short. 
Since $({\rm C}_1)_t$ is assumed, the action of $z$ (resp.$z^*$) on 
$V$ has $d+1$ distinct eigenvalues $\theta_0, \cdots, \theta_d$ 
(resp. $\theta_0^*, \cdots, \theta_d^*$) on $V$ 
by Proposition \ref{prop: C1}.
Let $V_i$ (resp. $V_i^*$) denote the eigenspace of $z$ (resp. $z*$) 
on $V$ belonging to $\theta_i$ (resp. $\theta_i^*$). 
Then we have 
\begin{eqnarray*}
V_i + V_{i+1} + \cdots + V_d &=& U_i + U_{i+1} + \cdots + U_d,\\
V_0^* + V_1^* + \cdots + V_i^* &=& U_0 + U_1 + \cdots + U_i
\end{eqnarray*}
for $0 \leq i \leq d$.
In particular, $U_0 = V_0^*,\, U_d = V_d$. 
Let $E_i$ (resp. $E_i^*$) denote the projection of 
$V=\bigoplus_{i=0}^d V_i$ (resp. $V=\bigoplus_{i=0}^d V_i^*$) onto 
$V_i$ (resp.$V_i^*$). Then the mappings 
\begin{eqnarray*}
&&F_i|_{V_i}:~V_i \longrightarrow U_i,\\
&&E_i|_{U_i}:~U_i \longrightarrow V_i
\end{eqnarray*}
are both bijections and inverses each other. 
Also the mappings 
\begin{eqnarray*}
&&F_i|_{V_i^*}:~V_i^*  \longrightarrow U_i,\\
&&E_i^*|_{U_i}:~U_i  \longrightarrow V_i^*
\end{eqnarray*}
are both bijections and inverses each other. 
In particular, by Theorem \ref{thm: shape} 
$${\rm dim}\, V_0=1, ~~{\rm dim}\, V_d^*=1.$$ 
By the argument in the proof of Proposition \ref{prop: q-Onsager}, 
the TD-relations ({\rm TD}) for $z,\, z^*$ imply 
\begin{eqnarray*}
z^* V_i &\subseteq& V_{i-1} + V_i + V_{i+1},\\
z\, V_i^* &\subseteq& V_{i-1}^* + V_i^* + V_{i+1}^*
\end{eqnarray*}
for $0 \leq i \leq d$, where $V_{-1}=V_{d+1}=V_{-1}^*=V_{d+1}^*=0$.

\medskip
Regard $V$ as an $\A$-module via $\iota_t$. 
Let $W$ be an irreducible $\A$-submodule of $V$. 
Set $W_i=W\cap V_i$, $W_i^*=W\cap V_i^*.$ Then 
\begin{eqnarray*}
z^* W_i &\subseteq& W_{i-1}+W_{i}+W_{i+1}, \\
z\,W_i^* &\subseteq& W_{i-1}^*+W_{i}^*+W_{i+1}^*
\end{eqnarray*}
for $0\leq i\leq d$. 
Since $W$ is irreducible as an $\A$-module 
and since $z$, $z^*$ are diagonalizable on $W$, 
the pair $z|_W,\, z^*|_W$ is a TD-pair on $W$. 
This implies that 
the eigenspace decompositions of 
$z|_W, z|_{W^*}$ are 
\begin{eqnarray*}
W &=& \bigoplus_{i=r}^{r+d'} W_i,\\
W &=& \bigoplus_{i=r'}^{r'+d'} W_i^*, 
\end{eqnarray*}
for some integers  $r,\, r'$,  where $d'$ is the diameter of 
the TD-pair $z|_W,\, z^*|_W \in {\rm End} (W)$. 
As we discussed in Section 1.2, the $\A$-module structure on $W$ 
can be extended to a $\T$-module structure on $W$ by using the split 
decomposition of the TD-pair $z|_W, z^*|_W$. 
(Note that the weight-space decomposition of the $\T$-module $W$ may 
be totally different from that of the $\T$-module $V$.)
By applying Theorem \ref{thm: shape} to the irreducible 
$\T$-module $W$, we have 
\[
{\rm dim}\,W_r={\rm dim}\,W_{r'+d'}^*=1. 
\]
First we want to show $r=0,\, r'+d'=d$, i.e., 
$W \supseteq V_0,\; W \supseteq V_d^*$.

\medskip
Since ${\rm dim}\, W_r=1$, we have $W_r=\C v$ for some nonzero element 
$v \in W_r$. Since 
$W \subseteq V_r + \cdots + V_d = U_r + \cdots + U_d$, 
we can express $v$ as 
$$ v = u_r + \cdots + u_d, $$
where $u_i=F_i v \in U_i$. Then $u_r \neq 0$, 
since $v \in W_r \subseteq V_r$ 
and $F_r|_{V_r} : V_r \longrightarrow U_r$ is a bijection. 

\begin{lemma}
\label{lemma: x y on ui} 
The action of $\T$ on $V$ satisfies the following $(i),\, (ii),\, (iii)$. 
\begin{enumerate}
\item[$(i)$] 
$x^j u_r=(\theta_r - \theta_{r+1}) \cdots (\theta_r - \theta_{r+j})\,u_{r+j}
~~~ (1 \leq j \leq d-r)$. 
\item[$(ii)$] 
$y u_r = 0$. 
\item[$(iii)$] 
$y^j u_{r+j} \in \C u_r ~~~ (1 \leq j \leq d-r)$. 
\end{enumerate}
\end{lemma}
Proof. Recall $z=z_t$ and so 
$z|_V = x|_V + \sum_{i=0}^d \theta_i F_i.$ 
Since $u_i \in U_i$, we have $z\, u_i = x\, u_i + \theta_i u_i$,  
so $z v = \theta_r u_r + (x\, u_r + \theta_{r+1} u_{r+1}) + \cdots + 
(x u_{d-1} + \theta_d u_d)$. 
Note $x\, u_i \in U_{i+1}$. 
On the other hand, 
since $v \in W_r \subseteq V_r$, 
we have $z\,v=\theta_r v=\theta_ru_r+\theta_ru_{r+1}+\cdots+\theta_ru_d$. 
Therefore we have 
$x\,u_{i-1}+\theta_iu_i=\theta_ru_i$, 
i.e., $x\,u_{i-1}=(\theta_r-\theta_i)u_i$
and we obtain (i) recursively. 

\medskip
Recall $z^*=z_t^*$ and so 
$z^*|_V=y|_V + \sum_{i=0}^d \theta_i^* F_i$. 
Since $u_i \in U_i$, we have $z^* u_i = y\, u_i + \theta_i^* u_i$, 
so $z^* v = y\, u_r + (y\, u_{r+1} + \theta_r^* u_r) + \cdots + 
(y\, u_d + \theta_{d-1}^* u_{d-1}) + (\theta_d^* u_d)$. 
Note $y\, u_i \in U_{i-1}$.  
On the other hand, 
since $z^* v \in W$ and $F_{r-1}W=0$, 
we have $y\, u_r = 0$, i.e., (ii) holds. 
Since $z^* v \in W$ and $F_r W = F_r W_r = \C u_r$, 
we have 
$y\, u_{r+1} + \theta_r^* u_r \in \C u_r$, i.e., (iii) holos for $j=1$. 
By $z^*|_V=y|_V + \sum_{i=0}^d \theta_i^* F_i$ and 
$y\, U_i \subseteq U_{i-1}$, 
we can write $z^{*j}u_i$ as a linear combination of 
$u_i,\, y\,u_i, y^2u_i, \ldots, y^ju_i$, 
in which the coefficient of $y^ju_i$ is $1$ if $i-j \geq r$. 
In particular for $v=u_r+\cdots+u_d$, 
the projection of $z^{*j}v$ onto $U_r$ by $F_r$ 
can be written as 
$$ F_r z^{*j} v = y^j u_{r+j} + c_1 y^{j-1} u_{r+j-1} + \cdots 
+ c_{j-1} y\, u_{r+1} + c_j u_r$$
for some $c_1,\, \cdots, c_{j-1},\, c_j \in \C$. 
Since $F_r z^{*j} v \in F_r W = F_r W_r = \C u_r$, (iii) holds by 
induction on $j$. 
\hfill $\Box $

\begin{proposition} 
\label{prop: endpt}
It holds that 
$W\supseteq V_0$ and $W\supseteq V_d^*.$
\end{proposition}
Proof.  We only show $W\supseteq V_0$, i.e., $r=0$;\,
$W\supseteq V_d^*$ is proved similarly, 
using $W_{r'+d'}^*$ in place of $W_r$. 
By Lemma \ref{lemma: x y on ui}, 
the action of $\T$ on $V$ satisfies 
\begin{eqnarray*}
& & y u_r = 0,\\
& & y^j x^j u_r \in \C u_r ~~ (j=0,1,2,\cdots).
\end{eqnarray*}
This implies 
$\T u_r \subseteq U_r + \cdots + U_d$, 
since $\T$ is linealy 
spanned by 
$k^n w_{\lambda}(x,y) ~(n \in \Z, \lambda \in \Lambda^{irr})$ 
by Theorem \ref{thm: T basis}. 
Since $V$ is irreducible as a $\T$-module, we have 
$V=\T u_r$ and hence $r=0$. 
\hfill $\Box $ 

\medskip
Thus for a finite-dimensional irreducible $\T$-module $V$ of type $s$, 
diameter $d$ and an irreducible $\A$-module $W \subseteq V$ via $\iota_t$, 
we have $W_0=V_0,\, W_d^*=V_d^*$. 
In particuler, $W=\A V_0$.

\medskip
Next we calculate how the eigenspace $V_0$ of 
$z|_V$ is mapped to $V_0^*$ by the projection 
$E_0^*: V=\bigoplus_{i=0}^dV_i^* \longrightarrow V_0^*.$ 
It holds that on $V$
$$ E_0^* = \prod_{j=1}^d 
\frac{z^* - \theta_j^*}{\theta_0^* - \theta_j^*}, $$
since the right-hand side vanishes on $V_j^* (1 \leq j \leq d)$ 
and is the identity map on $V_0^*$.
Write $V_0$ as $V_0=\C v$ for some nonzero element $v \in V_0$ 
and express $v$ as $v = u_0 + u_1 + \cdots + u_d$, 
where $u_i = F_i v \in U_i$.
Then we obtain 
$$ E_0^* u_i = {\Theta_i^*}^{-1} y^i u_i, \qquad
\Theta_i^*=\prod_{j=1}^i (\theta_0^* - \theta_j^*). $$
This is because 
$\bigl(\prod_{j=1}^i (z^* - \theta_j^*)\bigr)\, u_i = y^i u_i 
\in U_0 = V_0^*$ 
by $(z^* - \theta_j^*) |_{U_j} = y|_{U_j},\, y\, U_j \subseteq U_{j-1}$ 
and because 
$(z^* - \theta_j^*) |_{U_0} = \theta_0^* -\theta_j^*$ for 
$i+1 \leq j \leq d$. 
By Lemma \ref{lemma: x y on ui} with $r=0$, 
$$ u_i = \Theta_i^{-1}x^i u_0, \qquad 
\Theta_i=\prod_{j=1}^i (\theta_0 - \theta_j). $$
Since $y^i x^i u_0 = \sigma_i u_0$, we have 
$$ E_0^* u_i = \Theta_i^{-1}{\Theta_i^*}^{-1}
\sigma_i u_0 $$
and so 
$$ E_0^* v = \sum_{i=0}^d \Theta_i^{-1}{\Theta_i^*}^{-1}
\sigma_i u_0.$$
Note that $u_0 = F_0 v \neq 0$, since 
$F_0|_{V_0} : V_0 \longrightarrow U_0$ is a bijection. 
Thus by Remark \ref{rem: Drinfeld polynomial} in Section 1.3, we have 

\begin{proposition} 
\label{prop: proj}
For a finite-dimensional irreducible $\T$-module $V$ of type $s$ and 
diameter $d$, assume 
the condition $({\rm C}_1)_t$ for a nonzero $t \in \C$. 
Then for $v \in V_0$, it holds that 
$$ E_0^* v = {\Theta}^{-1}Q\,
P_V (t^2 + \varepsilon \varepsilon^* t^{-2})\, u_0,$$
where $u_0 = F_0 v$, 
\begin{eqnarray*} 
\Theta &=& (\theta_0 - \theta_1) \cdots (
\theta_0 - \theta_d)(\theta_0^* - \theta_1^*) \cdots 
(\theta_0^* - \theta_d^*),\\
Q &=& {(-1)}^d {(q-q^{-1})}^2 {(q^2-q^{-2})}^2 \cdots {(q^d-q^{-d})}^2, 
\end{eqnarray*} 
and $P_V(\lambda)$ is the Drinfel'd polynomial of 
the $\T$-module $V$ defined in Section 1.3: 
$$
P_V(\lambda)=Q^{-1} \sum^d_{i=0} \sigma_i(V) \prod^d_{j=i+1} 
{(q^j - q^{-j})}^2 (\varepsilon s^{-2} q^{2(d-j)} 
+ \varepsilon^*s^2 q^{-2(d-j)} - \lambda). 
$$

\end{proposition}

\medskip
\noindent
{\bf Proof of Theorem \ref{thm: C2}}. 
Suppose $P_V (t^2 + \varepsilon \varepsilon^* t^{-2}) \neq 0$. 
Then by Proposition \ref{prop: proj}, we have $E_0^* V_0 \neq 0$. 
Then $E_0^* V_0 = V_0^*$, since  $E_0^* V_0 \subseteq V_0^*$ and 
$dim\, V_0^*=1$. Let $W$ be an irreducible 
$\A$-submodule of $V$ 
via $\iota_t$. Then $W \supseteq V_0$ 
by Proposition \ref{prop: endpt}. Since $E_0^*$ is 
a polynomial of $z^*|_V$, $W$ is $E_0^*$-invariant and so 
$W \supseteq E_0^* W \supseteq E_0^* V_0 = V_0^*$, i.e., 
$W \supseteq U_0$ by $U_0=V_0^*$. 
We want to prove $W=V$. To do so, it is enough to show 
$\omega U_0 \subseteq W$ 
for every word $\omega$ in $x,\,y$, 
since $V=\T U_0$ and $\T U_0$ is linearly 
spanned by such $\omega U_0's$. 
Now $\omega U_0$ belongs to some $U_i$ and $x,\,y$ coincide with 
$z - \theta_i,\, z^* - \theta_i^*$ on $U_i$ respectively. 
Therefore $\omega U_0 \subseteq W$ implies 
$x \omega U_0 \subseteq W$ and 
$y \omega U_0 \subseteq W$, 
since $W$ is invariant under $z - \theta_i, z^* - \theta_i^*$.
This means induction works on the word length. 
Thus $\omega U_0 \subseteq W$ holds for every word $\omega$ in $x,\,y$ 
and it is shown that $P_V (t^2 + \varepsilon \varepsilon^* t^{-2}) 
\neq 0$ implies $W=V$, i.e., $V$ is irreducible as an $\A$-module.

\medskip
Suppose $P_V (t^2 + \varepsilon \varepsilon^* t^{-2}) = 0$. 
Then by Proposition \ref{prop: proj}, we have $E_0^* V_0=0$. 
This means $V_0 \subseteq V_1^* + \cdots + V_d^*$. Set 
$$ V_{i,i+1} = (V_0 + \cdots + V_i) \cap (V_{i+1}^* + \cdots + V_d^*)$$
for $0 \leq i \leq d-1$. 
Note $V_0 = V_{0,1}$. Then by 
$z V_j^* \subseteq V_{j-1}^* + V_j^* + V_{j+1}^*$ and 
$z^* V_j \subseteq V_{j-1} + V_j + V_{j+1}$, we have 
\begin{eqnarray*}
(z - \theta_i) V_{i,i+1} &\subseteq& V_{i-1,i},\\
(z^* - \theta_{i+1}^*) V_{i,i+1} &\subseteq& V_{i+1,i+2},
\end{eqnarray*}
where $V_{-1,0}=V_{d,d+1}=0$. 
Set $V'=V_{0,1} + V_{1,2} + \cdots + V_{d-1,d}$. 
Then $V'$ is $\langle z,z^*\rangle$-invariant. 
Since $V_0 \subseteq V' \subseteq V_1^* + \cdots + V_d^*$, the 
$\langle z,z^*\rangle$-invariant subspace $V'$ is 
a proper subspace of $V$. 
Thus it is shown that if 
$P_V (t^2 + \varepsilon \varepsilon^* t^{-2})=0$, 
then $V$ is not irreducible as an $\A$ -module. This comletes 
the proof of Theorem \ref{thm: C2}. 
\hfill $\Box $

\section{
The product formula for the Drinfel'd polynomial 
$P_V(\lambda)$ of a $\T$-module $V$ via $\varphi_s$: 
Proof of the surjectivity of $\sigma$ in Theorem \ref{thm: sigma}}

The augmented TD-algebra $\T = \T_q^{\varepsilon,\varepsilon^*}$ is by 
Proposition \ref{prop: phi}
embedded into the $U_q(sl_2)$-loop algebra $\LL = U_q(L(sl_2))$ 
via the injective algebra-homomorphism 
$$ \varphi_s: 
\T \longrightarrow \LL ~~~(x,\,y,\,k \mapsto x(s),\,y(s),\,sk_0 ~
\mbox{respectively}) $$
for each fixed nonzero $s \in \C$, where
\begin{eqnarray*}
x(s) &=& \alpha (s e_0^+ + \varepsilon s^{-1} e_1^- k_1), \qquad 
\alpha= -q^{-1}{(q-q^{-1})}^2,\\
y(s) &=& \varepsilon^* s e_0^- k_0 + s^{-1} e_1^+.
\end{eqnarray*}

\medskip
For $(\varepsilon,\varepsilon^*)=(1,1),\,(0,0)$, let
$$V=V(\ell_1,a_1) \otimes \cdots \otimes V(\ell_n,a_n)$$ 
be the tensor product of evaluation modules $V(\ell_i,a_i)$ for $\LL$ 
$(1 \leq \ell_i,\, a_i \in \C \backslash \{0\},\,1 \leq i \leq n)$ 
(see Section 1.4). 
We regard $V$ as a $\T$-module via the embedding $\varphi_s$. 
We call 
such a $\T$-module $V$ a {\it tensor product of evaluation 
modules via $\varphi_s$.}

\medskip
For $(\varepsilon,\varepsilon^*)=(1,0),$ let $\LL'$ denote the subalgebra of 
$\LL$ generated by $e_0^+,\, e_1^+,\, e_1^-,\, k_i^{\pm 1}$ $(i=0,\,1)$ and 
$V=V(\ell_1,a_1) \otimes \cdots \otimes V(\ell_n,a_n)$ the tensor product of 
evaluation modules 
$V(\ell_i,a_i)$ for $\LL'$ $(1 \leq \ell_i,\, a_i \in \C,\,1 \leq i \leq n)$: 
note that $e_0^-$ is missing from the set of generators for $\LL'$ and 
$a_i=0$ is allowed for the evaluation module $V(\ell_i,a_i)$ of $\LL'$ 
(see Section 1.4). 
We regard $V$ as a $\T$-module via the embedding $\varphi_s$, 
since the image of $\T$ by $\varphi_s$ is contained in $\LL'$ 
in the case of $(\varepsilon,\varepsilon^*)=(1,0)$. 
We call 
such a $\T$-module $V$  a {\it tensor product of evaluation 
modules via $\varphi_s$.}

\medskip
We treat such a $\T$-module 
$V=V(\ell_1,a_1) \otimes \cdots \otimes V(\ell_n,a_n)$  via $\varphi_s$ 
in one argument, regardless of $(\varepsilon,\varepsilon^*)$, 
and use the same notation $\LL$ for $\LL'$ in the case of  
$(\varepsilon,\varepsilon^*)=(1,0)$. So in this section, we understand 
in the case of $(\varepsilon,\varepsilon^*)=(1,0)$ that 
$\LL$ denotes the subalgebra of the $U_q(sl_2)$-loop algebra 
$U_q(L(sl_2))$ generated by 
$e_0^+,\, e_1^+,\, e_1^-,\, k_i^{\pm 1}$ ($i=0,\,1$) with  
$e_0^-$ missing from the set of generators, and that 
$a_i=0$ is allowed for the evaluation module $V(\ell_i,a_i)$. 

\medskip
For a $\T$-module $V=V(\ell_1,a_i) \otimes \cdots \otimes V(\ell_n,a_n)$ via 
$\varphi_s$, let $v_0^{(i)}, \cdots , v_{\ell_i}^{(i)}$ denote a standard 
basis of 
$V(\ell_i,a_i)$: we write 
$v_0=v_0^{(i)},\, v_1=v_1^{(i)},\, \cdots ,\, v_{\ell_i}=v_{\ell_i}^{(i)}$ 
for short. The action of $\T$ on 
$V(\ell_i,a_i) = \langle v_0, v_1, \cdots , v_{\ell_i} \rangle$ 
is 
\begin{eqnarray*}
k_0 v_j    &=& q^{2j-\ell_i}\,  v_j , \\
e_0^+ v_j  &=& a_i q\, [j+1]\, v_{j+1} , \\ 
\varepsilon^* e_0^- v_j 
           &=& \varepsilon^* a_i^{-1} q^{-1} [\ell_i-j+1]\, v_{j-1} , \\
e_1^+ v_j  &=& [\ell_i-j+1]\, v_{j-1} , \\
e_1^- v_j  &=& [j+1]\, v_{j+1},
\end{eqnarray*}
where $v_{-1}^{(i)} = v_{\ell_i +1}^{(i)} = 0$, and we understand 
$\varepsilon^* a_i^{-1} = 0$ if $(\varepsilon,\varepsilon^*)=(1,0)$ 
and $a_i=0$. Let $U_i$ denote the subspace of $V$ spanned by 
$v_{j_1} \otimes \cdots \otimes v_{j_n}$, where ($j_1, \cdots, j_n)$ runs 
through $0 \leq j_1 \leq \ell_1,\, \cdots,\, 0 \leq j_n \leq \ell_n$ 
such that 
$j_1 + \cdots + j_n = i$: 
$$ U_i = \bigoplus_{j_1 + \cdots + j_n = i} 
\C v_{j_1} \otimes \cdots \otimes v_{j_n}. $$
Then $k|_{U_i}= s q^{2i-d}$, so 
$$ V = \bigoplus_{i=0}^d U_i ~~~~~~~ (d = \ell_1 + \cdots + \ell_n) $$
is the eigenspace decomposition of $\varphi_s(k)$. 
We call 
$ V = \bigoplus_{i=0}^d U_i $ 
the {\it weight space decomposition} of the $\T$-module V via $\varphi_s$ 
and $U_0$ the {\it highest weight space}. 
Observe that 
\begin{eqnarray*}
& & {\rm dim}\, U_0 = 1 , \\
& & x\, U_i \subseteq U_{i+1}, ~~~ y\, U_i \subseteq U_{i-1}
\end{eqnarray*}
for $0 \leq i \leq d $, where $U_{-1}=U_{d+1}=0$. 
So the 1-dimensional space $U_0$ is invariant under $y^i x^i$. 
Define the sequence 
${\{ \sigma_i\}}^{\infty}_{i=0}$ of scalars $\sigma_i=\sigma_i(V)$ 
by 
$$ y^i x^i |_{U_0} = \sigma_i . $$ 
Then $\sigma_0=1\,, \sigma_i=0 ~(d+1 \leq i )$. 
Note that the $\T$-module $V$ via $\varphi_s$ is not necessarily 
irreducible and $\sigma_d=0$ is possible. Define the Drinfel'd 
polynomial $P_V(\lambda)$ of the $\T$-module $V$ via $\varphi_s$ by 
\begin{eqnarray}
\label{D-poly}
P_V(\lambda) &=& Q^{-1} \sum^d_{i=0} \sigma_i(V) \prod^d_{j=i+1} 
{(q^j - q^{-j})}^2 (\varepsilon s^{-2} q^{2(d-j)} 
+ \varepsilon^* s^2 q^{-2(d-j)} - \lambda),  \\
\label{Q}
 Q &=&Q_d= {(-1)}^d {(q-q^{-1})}^2 {(q^2 - q^{-2})}^2 
\cdots {(q^d - q^{-d})}^2 .  
\end{eqnarray}
Since $\sigma_0=1$, 
$P_V (\lambda)$ is a monic polynomial of degree $d$. Observe 
$$ \sigma_d (V) = Q\cdot P_V (\varepsilon s^{-2} + \varepsilon^* s^2). $$

\medskip
More generally the Drinfel'd polynomial $P_V(\lambda)$ is defined 
in the same way 
for a finite-dimensional $\T$-module $V$ that has 
the following properties:  
\begin{quote}
$({\rm D})_0$:~  $k$ is diagonalizable on $V$ with 
$V=\bigoplus_{i=0}^d U_i$,~
$k|_{U_i}=s q^{2i-d} ~~ (0 \leq i \leq d)$ \\
~~~~~~~~~for some nonzero constant $s$. \\
\,$({\rm D})_1$:~  ${\rm dim}\, U_0=1$. 
\end{quote}
By the relations $({\rm TD})_0:\,
k k^{-1}=k^{-1}k=1,\, kxk^{-1}=q^2x,\, kyk^{-1}=q^{-2}y$,
it holds that $x\, U_i \subseteq U_{i+1},\, 
y\, U_i \subseteq U_{i-1} ~(0 \leq i \leq d)$, 
where $U_{-1}=U_{d+1}=0$. Thus $\sigma_i(V)$'s are defined 
as before 
and hence $P_V(\lambda)$ by (\ref{D-poly}), (\ref{Q}). 
The eigenspace decomposition 
and the subspace $U_0$ in $({\rm D})_0$ are called the 
{\it weight-space decomposition} and the 
{\it highest weight space} of the $\T$-module $V$ 
respectively. 
The nonzero scalar $s$ and 
the nonnegative integer $d$ in $({\rm D})_0$ 
are called 
the {\it type} 
and the {\it diameter} of the $\T$-module $V$ respectively. 
We further consider the following property for a $\T$-module $V$ 
that satisfies $({\rm D})_0$, $({\rm D})_1$ with diameter $d$: 
\begin{quote}
$({\rm D})_2$:~  $\sigma_d(V)\neq 0$. 
\end{quote}
\begin{lemma}
\label{lemma: PbarW}
Let $V$ be a finite-dimensional $\T$-module that satisfies 
the properties $({\rm D})_0$, $({\rm D})_1$. 
Consider the $\T$-submodule $W=\T U_0$,  
where $U_0$ is the highest weight space of 
the $\T$-module $V$. 
Let $M$ be a maximal $\T$-submodule of $W$. 
Set $\overline{W}=W/M$. 
Then the $\T$-submodule $W$ and 
the quotient $\T$-module $\overline{W}$ 
satisfy $({\rm D})_0$, $({\rm D})_1$ as well. 
Furthermore if $V$ satisfies $({\rm D})_2$ with diameter $d$, 
then so do the $\T$-modules $W$ and $\overline{W}$ and 
it holds that 
\begin{eqnarray*}
(i)~~
\sigma_i(V)=\sigma_i(W)=\sigma_i(\overline{W})~~
(0 \leq i \leq d),
~~~~~~~~~~~~~~~~~~~~~~~~~~~~~~~~~~~~~~~\\
(ii)~~
P_V(\lambda)=P_W(\lambda)=P_{\overline{W}}(\lambda).
~~~~~~~~~~~~~~~~~~~~~~~~~~~~~~~~~~~~~~~~~~~~~~~~~~~~~
\end{eqnarray*}
\end{lemma}

Lemma \ref{lemma: PbarW} follows 
from Lemma \ref{lemma: weight-space decomposition}, 
since $\overline{W}$ is irreducible as a $\T$-module.

\medskip
In what follows, we fix a nonzero scalar $s \in \C$ arbitrarily and 
we only treat 
finite-dimensional $\T$-modules via $\varphi_s$ that satisfy 
the above properties $({\rm D})_0$, $({\rm D})_1$. 
In this case, the weight space decomposition of a $\T$-module $V$ 
coincides with that of the $\LL$-module $V$, 
since $\varphi_s(k)=s\,k_0$. 
Note that the tensor product of evaluation modules $V(\ell_i, a_i)$ 
$(1 \leq i \leq n)$ satisfies $({\rm D})_0$, $({\rm D})_1$ 
and has type $s$, diameter $d=\ell_1+ \cdots + \ell_n$. 
If $V,\,V'$ are $\T$-modules via $\varphi_s$, then 
the tensor product $V \otimes V'$  
becomes a $\T$-module via $\Delta \circ \varphi_s$, 
where $\Delta: \LL \longrightarrow \LL \otimes \LL$ is the coproduct. 
Furthermore,  
if the $\T$-modules $V,\,V'$ via $\varphi_s$ satisfy 
the properties $({\rm D})_0$, $({\rm D})_1$, 
so does the tensor product $V \otimes V'$ 
as a $\T$-module via $\varphi_s$ 
and so the Drinfel'd polynomial $P_{V \otimes V'} (\lambda) $ 
is defined. We have the following product formula.   
\begin{theorem}
\label{thm: product formula}
Let $V,\,V'$ be finite-dimensional $\T$-modules via $\varphi_s$ 
that satisfy the properties $({\rm D})_0$, $({\rm D})_1$. 
Assume that $V'$ is afforded by a tensor 
product of evaluation modules via $\varphi_s$. 
Then the following $(i)$, $(ii)$ holds. 

\begin{enumerate}
\item[$(i)$]
The Drinfel'd polynomial 
$P_{V \otimes V'} (\lambda) $ of the $\T$-module $V \otimes V'$ via 
$\varphi_s$ is 
$$ P_{V \otimes V'} (\lambda) = P_V (\lambda) P_{V'} (\lambda).$$
\item[$(ii)$]
The Drinfel'd polynomial 
$P_{V(\ell,a)} (\lambda)$ of the $\T$-module $V(\ell,a)$ via $\varphi_s$ is 
$$ P_{V(\ell,a)} (\lambda) = \prod_{c \in S(\ell,a)} 
(\lambda + c + \varepsilon \varepsilon^* c^{-1}), $$ 
where 
$$ S(\ell,a) = \{ a\, q^{2i-\ell+1} ~|~ 0 \leq i \leq \ell-1 \}.  $$
We understand that if $(\varepsilon , \varepsilon^*)=(1,0)$ 
and $a=0$, $S(\ell,a)$ is the multiset with 0 appearing $\ell$ times 
and $P_{V(\ell,0)} = \lambda^\ell$.
\end{enumerate}

\end{theorem}

\medskip
To prove Theorem \ref{thm: product formula}, 
we prepare two lemmas and a proposition. Let 
$V,\,V'$ be $\T$-modules via $\varphi_s$ 
as in Theorem \ref{thm: product formula} 
and have weight-space decompositions 
\begin{eqnarray*}
V  &=& \bigoplus^d_{i=0} U_i, \\
V' &=& \bigoplus^{d'}_{i=0} U_i',
\end{eqnarray*}
respectively. Then the $\T$-module $V \otimes V'$ via $\varphi_s$ has 
weight-space decomposition
$$ V \otimes V' = \bigoplus^{d+d'}_{i=0} \widetilde{U_i}, $$
where
$$ \widetilde{U_i} = \bigoplus_{i_1+i_2=i} U_{i_1} \otimes U_{i_2} 
~~~~~ (0 \leq i \leq d+d').$$
\begin{lemma}
\label{lemma: x(s) y(s)}
Set $x(s)=\varphi_s(x),\, y(s)=\varphi_s(y)$. 
Then the action of $x(s),\,y(s)$ on $U_i \otimes V'$ are 
\begin{eqnarray*}
x(s)|_{U_i \otimes V'}
&=&x(s)|_{U_i} \otimes 1_{V'} + 1_{U_i} \otimes x(q^{2i-d}s)|_{V'}, \\
y(s)|_{U_i \otimes V'}
&=&y(s)|_{U_i} \otimes 1_{V'} + 1_{U_i} \otimes y(q^{2i-d}s)|_{V'}.
\end{eqnarray*}
\end{lemma}
Proof. These identities follow directly from 
$x(s)=\alpha (s e^+_0 + \varepsilon s^{-1} e^-_1 k_1),\,  
y(s)=\varepsilon^* s e^-_0 k_0 + s^{-1} e^+_1$ 
and the coproduct 
$\Delta$ that sends  
$e^+_i,\,e^-_i k_i,\,k_i$ 
to 
$e^+_i \otimes 1 + k_i \otimes e^+_i,\:
e^-_i k_i \otimes 1 + k_i \otimes e^-_i k_i,\:
k_i \otimes k_i$ respectively. 
\hfill $\Box $ 

\begin{lemma}
\label{lemma: V otimes V(1,a)}
Assume $V'=V(1,a)$, an evaluation module of diameter $1$, and let 
$V(1,a)=\langle v_0,\, v_1 \rangle $ be a standard basis. 
For $u \in U_m ~(0 \leq m \leq d)$ and $1 \leq i$, we have 
\begin{eqnarray*}
(i) ~~
x^i (u \otimes v_0) 
&=& (x^i u) \otimes v_0 + \alpha\, q\,[\,i\,]\, c_i(m)(x^{i-1} u) \otimes v_1, 
~~~~~~~~~~~~~~~~~~~~~~~~~~~~~~~~~\\
x^i (u \otimes v_1)  
&=& (x^i u) \otimes v_1,
\end{eqnarray*}
~~~~~~~~~where ~
$c_i(m) = a\,s\,q^{i+2m-d-1} + \varepsilon s^{-1} q^{-i-2m+d+1}, $
\begin{eqnarray*}
(ii) ~~
y^i (u \otimes v_0) 
&=& (y^i u) \otimes v_0, \\
y^i (u \otimes v_1)  
&=& (y^i u) \otimes v_1+[\,i\,]\, c_i^*(m)(y^{i-1} u) \otimes v_0, 
~~~~~~~~~~~~~~~~~~~~~~~~~~~~~~~~~~~
\end{eqnarray*}
~~~~~~~~~where ~$ c_i^*(m)=
\varepsilon^* a^{-1} s\,q^{-i+2m-d+1}+s^{-1} q^{i-2m+d-1}.$ 

\end{lemma}
Proof. 
Recall $e^+_0 v_0=q\,a\,v_1,\, e^-_1 v_0=v_1,\, e^+_0 v_1 = e^-_1 v_1 = 0,\, 
\varepsilon^* e^-_0 v_1=\varepsilon^* a^{-1} q^{-1} v_0,\, e^+_1 v_1=v_0,\,
\varepsilon^* e^-_0 v_0=e^+_1 v_0=0,\, k_0 v_0= q^{-1}v_0,\,\, k_0v_1=q v_1.$ 
We proceed by induction on $i$. For $i=1$, we have 
by Lemma \ref{lemma: x(s) y(s)} 
\begin{eqnarray*}
x(u \otimes v_0) 
&=& \bigl(x(s)u\bigr) \otimes v_0 +
 u \otimes \bigl(x(q^{2m-d} s)\, v_0\bigr) \\
&=& (xu) \otimes v_0 + \alpha\, (q^{2m-d} s\,q\,a + \varepsilon 
    q^{-2m+d} s^{-1} q)\, u \otimes v_1 \\
&=& (xu) \otimes v_0 + \alpha\, q\, c_1(m)\, u \otimes v_1, \\
x(u \otimes v_1) 
&=& \bigl(x(s)u\bigr) \otimes v_1 + 
u \otimes \bigl(x(q^{2m-d} s)\, v_1\bigr) \\
&=& (xu) \otimes v_1,
\end{eqnarray*}
and 
\begin{eqnarray*}
y(u \otimes v_0) 
&=& \bigl(y(s)u\bigr) \otimes v_0 + 
u \otimes \bigl(y(q^{2m-d} s)\, v_0\bigr) \\
&=& (yu) \otimes v_0, \\
y(u \otimes v_1) 
&=& \bigl(y(s)u\bigr) \otimes v_1 + 
u \otimes \bigl(y(q^{2m-d} s)\, v_1\bigr) \\
&=& (yu) \otimes v_1 + (\varepsilon^* q^{2m-d} s\, a^{-1} + 
    q^{-2m+d} s^{-1})\, u \otimes v_0 \\
&=& (yu) \otimes v_1 + c^*_1(m)\, u \otimes v_0.
\end{eqnarray*}
For $i \geq 2$, we have by Lemma \ref{lemma: x(s) y(s)} 
and induction on $i$ 
\begin{eqnarray*}
x^i(u \otimes v_0) 
&=& x^{i-1} 
\Bigl((xu) \otimes v_0 + \alpha\, q\, c_1(m)\, u \otimes v_1\Bigr) \\
&=& (x^i u) \otimes v_0 + 
\alpha\, q[i-1]\, c_{i-1}(m+1)\,(x^{i-1} u) \otimes v_1 \\
& & + \alpha\, q\, c_1(m)\,(x^{i-1} u) \otimes v_1 \\
&=& (x^i u) \otimes v_0 + 
\alpha\, q\,[\,i\,]\, c_i(m)\,(x^{i-1} u) \otimes v_1,
\end{eqnarray*}
since $[i-1] c_{i-1}(m+1) + c_1(m) = [\,i\,]\, c_i (m)$ , and 
\begin{eqnarray*}
y^i(u \otimes v_1) 
&=& y^{i-1} 
\Bigl((yu) \otimes v_1 + c^*_1(m)\, u \otimes v_0\Bigr) \\
&=& (y^i u) \otimes v_1 + 
[i-1]\, c^*_{i-1}(m-1)\,(y^{i-1} u) \otimes v_0 \\
& & + c^*_1(m)\, (y^{i-1} u) \otimes v_0 \\
&=& (y^i u) \otimes v_1 + [\,i\,]\, c^*_i(m)\,(y^{i-1} u) \otimes v_0,
\end{eqnarray*}
since $[i-1]\, c^*_{i-1}(m-1) + c^*_1(m)=[\,i\,]\,c^*_i(m)$. Also we have 
$x^i(u \otimes v_1)=x^{i-1}((xu) \otimes v_1)=(x^i u) \otimes v_1,\; 
y^i(u \otimes v_0)=y^{i-1} ((yu) \otimes v_0) = (y^i u) \otimes v_0$ 
by induction on $i$. 
\hfill $\Box$

\begin{proposition}
\label{prop: tilde sigma}
Assume $V'=V(1,a)$, an evaluation module of diameter $1$. 
For the $\T$-modules 
$V$ and $V \otimes V'$ via $\varphi_s$, set $\sigma_i=\sigma_i(V)$ and  
$\widetilde{\sigma_i}=\sigma_i(V \otimes V')$. 
Then for $i \geq 1$, we have 
$$ \widetilde{\sigma_i} = \sigma_i -{(q^i - q^{-i})}^2 
\bigl(a + \varepsilon\,\varepsilon^*a^{-1} + 
\varepsilon\, s^{-2} q^{2(d+1-i)} + \varepsilon^* s^2 
q^{-2(d+1-i)}\bigr) 
\sigma_{i-1},$$
where $d$ is the diameter of the $\T$-module $V$. We understand that 
$\varepsilon^* a^{-1}=0$ if $(\varepsilon,\varepsilon^*)=(1,0)$ and $a=0$. 
\end{proposition}
Proof. 
Let $V=\bigoplus^d_{i=0} U_i$ denote the weight-space 
decomposition of the $\T$-module $V$ and $V(1,a)=\langle v_0,v_1\rangle$ 
a standard basis of $V'$. Choose a nonzero vector $u_0 \in U_0$. 
Then $u_0 \otimes v_0$ 
spans the highest weight space of $V \otimes V'$. 
We have by Lemma \ref{lemma: V otimes V(1,a)} 
\begin{eqnarray*}
y^i x^i(u_0 \otimes v_0)
&=& y^i\, \Bigl((x^i u_0) \otimes v_0 + 
\alpha\, q\,[\,i\,]\, c_i(0)\,(x^{i-1} u_0) \otimes v_1\Bigr)\\
&=& (y^i x^i u_0) \otimes v_0 \\
& & + \alpha\, q\,[\,i\,]\, c_i(0)\,\Bigl((y^i x^{i-1} u_0) \otimes v_1 + 
[\,i\,]\, c^*_i (i-1)\,(y^{i-1} x^{i-1} u_0) \otimes v_0\Bigr) \\
&=& \sigma_i u_0 \otimes v_0 + 
\alpha\, q\,{[\,i\,]}^2 c_i(0)\, c^*_i (i-1)\, \sigma_{i-1} u_0 \otimes v_0,
\end{eqnarray*}
since $y^i x^{i-1} u_0=0$. So it holds that 
\begin{eqnarray*}
\widetilde{\sigma_i}
&=& \sigma_i + 
\alpha\, q\,{[\,i\,]}^2 c_i(0)\, c^*_i (i-1)\, \sigma_{i-1} \\
&=& \sigma_i - 
{(q^i - q^{-i})}^2\bigl(a + \varepsilon \varepsilon^* a^{-1} + 
\varepsilon s^{-2} q^{2(d+1-i)} + 
\varepsilon^* s^2 q^{-2(d+1-i)}\bigr) \sigma_{i-1}.
\end{eqnarray*}
\hfill $\Box$


\bigskip
\noindent
{\bf Proof of Theorem \ref{thm: product formula}}. 
We first treat the case of $V'=V(1,a)$, an evaluation module of diameter 1. 
By Proposition \ref{prop: tilde sigma}, 
$$ \widetilde{\sigma_i} = \sigma_i - {(q^i -q^{-i})}^2 \sigma_{i-1} 
\Bigl( (\varepsilon s^{-2} q^{2(d+1-i)} + 
\varepsilon^* s^2 q^{-2(d+1-i)}- \lambda) + 
(\lambda + a + \varepsilon \varepsilon^* a^{-1}) \Bigr), $$
and we have, with $\widetilde{Q} = 
{(-1)}^{d+1} {(q-q^{-1})}^2{(q^2-q^{-2})}^2 \cdots 
{(q^{d+1} -q^{-d-1})}^2 $,
\begin{eqnarray*}
P_{V \otimes V(1,a)} (\lambda) 
&=& {\widetilde{Q}}^{-1} 
\sum^{d+1}_{i=0} \,
\widetilde{\sigma}_i 
\prod^{d+1}_{j=i+1} 
(q^j - q^{-j})^2 
\bigl(\varepsilon s^{-2} q^{2(d+1-j)} + 
\varepsilon^* s^2 q^{-2(d+1-j)} - \lambda \bigr) \\  
&=& {\widetilde{Q}}^{-1}
\sum^{d+1}_{i=0}
 \sigma_i 
\prod^{d+1}_{j=i+1} 
(q^j - q^{-j})^2 
\bigl(\varepsilon s^{-2} q^{2(d+1-j)} + 
\varepsilon^* s^2 q^{-2(d+1-j)} - \lambda \bigr) \\
&-&  {\widetilde{Q}}^{-1} 
\sum^{d}_{i-1=0} 
\sigma_{i-1} 
\prod^{d+1}_{j=i} 
(q^j - q^{-j})^2 
\bigl(\varepsilon s^{-2} q^{2(d+1-j)} + 
\varepsilon^* s^2 q^{-2(d+1-j)} - \lambda \bigr) \\
&-&  {\widetilde{Q}}^{-1} 
\sum^{d}_{i-1=0} 
\sigma_{i-1} 
(\lambda + a + \varepsilon\, \varepsilon^* a^{-1})
(q^{d+1}-q^{-d-1})^2  \\ 
&&\qquad \qquad \times \prod^{d}_{j-1=i} 
(q^{j-1}-q^{-j+1})^2 
\bigl(\varepsilon s^{-2} q^{2(d-j+1)} + 
\varepsilon^* s^2 q^{-2(d-j+1)} - \lambda \bigr). 
\end{eqnarray*}
This equals 
$(\lambda + a +\varepsilon\, \varepsilon^* a^{-1})\,
P_V (\lambda)$, 
since $\sigma_{d+1} = \sigma_{d+1} (V) = 0 $ 
and so the first and second terms cancel out. 
This argument is valid even if $V$ is 
the trivial module, i.e., ${\rm dim}\,V=1$, 
${e_i^{\pm} |}_V = 0,\,{k_i^{\pm 1} |}_V = 1$. 
In this case, $V \otimes V(1,a) \simeq V(1,a)$ 
and it is easily checked that $P_V(\lambda)=1$ and 
\begin{eqnarray}
\label{PV(1,a)}
P_{V(1,a)}(\lambda)=\lambda+a+\varepsilon\, \varepsilon^* a^{-1}. 
\end{eqnarray}
Thus in the case of $V'=V(1.a)$, we have 
\begin{eqnarray}
\label{P V otimes V(1,a)}
P_{V \otimes V(1,a)} (\lambda) = P_V (\lambda) P_{V(1,a)} (\lambda).
\end{eqnarray}

\medskip
Next we treat the case $V'=V(\ell,a)$, 
an evaluation module of diameter $\ell$. 
We want to show 
\begin{eqnarray}
\label{P V otimes V(ell,a)}
P_{V \otimes V(\ell,a)} (\lambda) = P_V (\lambda) P_{V(\ell,a)} (\lambda)
\end{eqnarray}
for every integer $\ell \geq 1$ by induction on $\ell$.  
To do so, we prepare a lemma below that gives 
an embedding of $V'=V(\ell,a)$ 
into $V(\ell-1,a\,q^{-1}) \otimes V(1,a\,q^{\ell-1})$ 
as an $\LL$-submodule. 
Start with the evaluation modules 
$V(\ell-1,a\,q^{-1}),\,V(1,a\,q^{\ell-1})$ for $\LL$. 
Let $V(\ell-1,a\,q^{-1})=\langle u_0,u_1,\cdots,u_{\ell-1}\rangle,\,
V(1,a\,q^{\ell-1})=\langle v_0,v_1 \rangle$ 
be standard bases of the evaluation modules. 
By direct calculations, we have the following lemma. 
\begin{lemma}
\label{lemma: ev module}
Consider the tensor product 
$V(\ell-1,a\,q^{-1}) \otimes V(1,a\,q^{\ell-1})$ 
of evaluation modules as an $\LL$-module 
via the coproduct $\Delta$. 
Set 
$$ w_i=q^{-i} u_i \otimes v_0 + u_{i-1} \otimes v_1 
\in V(\ell-1,a\,q^{-1}) \otimes V(1,a\,q^{\ell-1}) $$ 
for $0 \leq i \leq \ell$, where $u_{-1}=u_{\ell}=0$. 
Then 
$\C w_0$ is the highest weight space of the $\LL$-module 
$V(\ell-1,a\,q^{-1}) \otimes V(1,a\,q^{\ell-1})$. 
Set $W=\LL w_0$. 
Then $$ W \simeq V(\ell,a) $$ as $\LL$-modules 
with 
$$ W = \langle w_0, w_1, \cdots, w_{\ell} \rangle$$ 
a standard basis for $W$. 
\end{lemma}

\medskip
Consider the $\LL \otimes \LL$-modules 
$V \otimes V(\ell,a)$ 
and 
$V \otimes 
\bigl(V(\ell-1,a\,q^{-1}) \otimes V(1,a\,q^{\ell-1})\bigr)$.  
Regard them as $\LL$-modules 
via the coproduct 
$\Delta: \LL \longrightarrow \LL \otimes \LL$  
and then as $\T$-modules via $\varphi_s$. 
Choose nonzero vectors $u,\,w$ from the highest weight spaces of 
$V$, $V(\ell-1,a\,q^{-1}) \otimes V(1,aq^{\ell-1})$ respectively. 
Then $\C \,u \otimes w$ is the highest weight space of 
$V \otimes \bigl(V(\ell-1,a\,q^{-1}) \otimes V(1,aq^{\ell-1})\bigr)$ 
as an $\LL$-module and hence as a $\T$-module via $\varphi_s$. 
Set $W=\LL w$. 
The properties $({\rm D})_0$, $({\rm D})_1$ hold for the $\T$-module  
$V \otimes \bigl(V(\ell-1,a\,q^{-1}) \otimes V(1,aq^{\ell-1})\bigr)$ 
and its $\T$-submodule 
$V \otimes W$. 
The Drinfel'd polynomials of these $\T$-modules coincide, since they 
share the common highest weight space and have the same diameter.   
On the other hand, $V \otimes V(\ell,a)$ is isomophic to 
$V \otimes W$ as $\LL$-modules by Lemma \ref{lemma: ev module} and so 
as $\T$-modules via $\varphi_s$. In particular, the Drinfel'd polynomial 
of $V \otimes V(\ell,a)$ 
coincides with that of $V \otimes W$ as $\T$-modules via $\varphi_s$. 
Therefore $V \otimes V(\ell,a)$ and 
$V \otimes \bigl(V(\ell-1,a\,q^{-1}) \otimes V(1,aq^{\ell-1})\bigr)$ 
have the same Drinfel'd polynomial 
as $\T$-modules via $\varphi_s$. 
The Drinfel'd polynomial of the $\T$-module 
$V \otimes \bigl(V(\ell-1,a\,q^{-1}) \otimes V(1,aq^{\ell-1})\bigr)$ 
is the product of those of the $\T$-modules 
$V \otimes V(\ell-1,a\,q^{-1})$ and $V(1,aq^{\ell-1})$ 
by (\ref{P V otimes V(1,a)}), 
since 
$V \otimes \bigl(V(\ell-1,a\,q^{-1}) \otimes V(1,aq^{\ell-1})\bigr)$ 
is isomorphic to 
$\bigl(V \otimes V(\ell-1,a\,q^{-1})\bigr) \otimes V(1,aq^{\ell-1})$ 
as $\T$-modules via $\varphi_s$. 
By induction on $\ell$, the formula (\ref{P V otimes V(ell,a)}) 
holds for $\ell -1$, so we have   
$P_{V \otimes V(\ell-1,a\,q^{-1})} = 
P_V P_{V(\ell-1,a\,q^{-1})}$. 
Therefore $P_{V \otimes V(\ell,a)} = 
P_V P_{V(\ell-1,a\,q^{-1})}\,
P_{V(1,aq^{\ell-1})}$. 
On the other hand, 
$P_{V (\ell-1,\,aq^{-1})} P_{V(1,a\,q^{-1})}=
P_{V (\ell-1,a\,q^{-1}) \otimes V(1,aq^{-1})}
=P_W = P_{V(\ell,a)}$ by the same argument. 
This proves the formula (\ref{P V otimes V(ell,a)}).

\medskip
Finally we treat the general case $V'=V'' \otimes V(\ell,a)$, 
where $V''$ is afforded by a tensor product of evaluation modules  
via $\varphi_s$. 
By (\ref{P V otimes V(ell,a)}),
$P_{V \otimes V'}=
P_{V \otimes V'' \otimes V(\ell,a)}=
P_{V \otimes V''} P_{V(\ell,a)}$. 
By induction on ${\rm dim}\, V''$, 
$P_{V \otimes V''}=P_V P_{V''}$. 
So $P_{V \otimes V'}=P_V P_{V''} P_{V(\ell,a)}$. 
By (\ref{P V otimes V(ell,a)}),
$P_{V''} P_{V(\ell,a)}=P_{V'' \otimes V(\ell,a)} 
= P_{V'}$. 
So $P_{V \otimes V'}=P_V P_{V'}$. 
This completes the proof of Theorem \ref{thm: product formula} (i). 

\medskip
By Lemma \ref{lemma: ev module} and 
Theorem \ref{thm: product formula} (i), 
we have 
$$
P_{V(\ell,a)}(\lambda)=\prod_{c \in S(\ell,a)} P_{V(1,c)}(\lambda). 
$$
By (\ref{PV(1,a)}), 
$P_{V(1,c)}(\lambda)=\lambda+c+\varepsilon\, \varepsilon^* c^{-1}.$ 
This completes the proof of 
Theorem \ref{thm: product formula} (ii).  
\hfill $\Box$

\bigskip
\noindent
{\bf Proof of the surjectivity of $\sigma$ in Theorem \ref{thm: sigma}}.
Given an arbitrary monic polynomial $P(\lambda)$ of degree $d$ 
and an arbitrary nonzero $s \in \C$ 
such that $P(\varepsilon s^{-2}+\varepsilon^* s^2) \neq 0$,
we show that there exists an irreducible 
$\T$-module $V$ of type s and diameter $d$ that has 
Drinfel'd polynomial $P(\lambda)$, i.e., $P_V (\lambda)=P (\lambda)$. 
Let $\lambda_1, \lambda_2, \cdots, \lambda_d$ denote the roots of 
$P(\lambda)$, allowing repetition. 
For each $i$ $(1 \leq i \leq d)$, choose $a_i \in \C$ 
such that 
$$ \lambda_i + a_i + \varepsilon \varepsilon^* a_i^{-1} = 0 . $$
If $(\varepsilon,\varepsilon^*)=(1,1)$, the equation 
$\lambda_i + a_i + a_i^{-1}=0$ has nonzero solutions for $a_i$ 
: we choose one of then and fix it. 
If $(\varepsilon,\varepsilon^*)=(1,0)$ or (0,0), we understand that the 
equation is $\lambda_i + a_i = 0$ and $a_i = - \lambda_i$. 
Observe that if $(\varepsilon,\varepsilon^*)=(0,0)$, 
then $\lambda_i \neq 0 ~(1 \leq i \leq d)$ by the condition 
$P(\varepsilon s^{-2}+\varepsilon^* s^2) \neq 0$, so 
$a_i \neq 0 ~(1 \leq i \leq d)$. 
Consider the $\T$-module $V$ via $\varphi_s$, 
where $\T=\T_q^{(\varepsilon,\varepsilon^*)}$ and
$$V=V(1,a_1) \otimes V(1,a_2) \otimes \cdots \otimes V(1,a_d).$$
By Theorem \ref{thm: product formula}, 
it holds that $P_V(\lambda)=P(\lambda)$. 
Choose a nonzero vector $w$ from the highest weight space of $V$ 
and set $W=\T w$. Let $M$ be a maximal $\T$-submodule of $W$. 
Observe $w \notin M$, since $M \neq W$. The quotient 
$\T$-module $\overline{W}=W/M$ is irreducible. 
Since $P_V(\varepsilon s^{-2}+\varepsilon^* s^2) \neq 0$, 
we have $\sigma_d(V) \neq 0$ 
by Remark \ref{rem: Drinfeld polynomial}. 
By Lemma \ref{lemma: PbarW}, 
$P_V(\lambda)=P_{\overline{W}}(\lambda)$. 
Thus $\overline{W}$ is the desired $\T$-module. 
\hfill $\Box$

\section
{Irreducibility of a $\T$-module 
$\tilde{V}=V \otimes V(\ell,a)$ via $\varphi_s$}

For the augmented TD-algebra $\T=\T_q^{\varepsilon,\varepsilon^*}$, 
we have so far established the bijectivity of the mapping 
$$ {Irr}_d^s (\T) \longmapsto \mathcal P_d^s ~~~ 
( V \longmapsto P_V(\lambda)) ,$$
namely, the set of finite-dimensional irreducible $\T$-modules 
of type $s$ and diameter $d$ 
are parametrized up to isomorphism by monic polynomials of 
degree $d$ that do not vanish at $\varepsilon s^{-2} + \varepsilon^* s^2$. 
Given a polynomial $P(\lambda) \in \mathcal P_d^s$, 
we want to construct explicitly 
a $\T$-module via $\varphi_s$ that belongs to ${Irr}_d^s (\T)$ and has 
$P(\lambda)$ as its Drinfel'd polynomial. 
In this section, we prepare a key proposition to the construction of such 
$\T$-modules via $\varphi_s$. 
The construction itself will be discussed in the next section. 
We consider $\T$-modules 
\begin{eqnarray*}
V&=&V(\ell_1,a_1) \otimes \cdots \otimes V(\ell_n,a_n),\\
\widetilde{V}&=&V \otimes V(1,a)
\end{eqnarray*}
via $\varphi_s$, where 
$1 \leq n,\, 1 \leq \ell_i ~(1 \leq i \leq n)$ and 
$a,\,a_i$ are allowed to be zero if $(\varepsilon,\varepsilon^*)=(1,0)$. 
The diameters of $V,\,\widetilde{V}$ are 
$d=\ell_1 + \cdots + \ell_n$, $d+1$ respectively. 
Observe that 
$\sigma_{d+1} (\widetilde{V})=\widetilde{Q}\, P_{\widetilde{V}} 
(\varepsilon s^{-2}+\varepsilon^* s^2)$ 
for some nonzero scalar $\widetilde{Q}$ by 
Remark \ref{rem: Drinfeld polynomial}. 
We have $P_{\widetilde{V}}(\lambda)=
P_V (\lambda) P_{V(1,a)} (\lambda)$ 
by Theorem \ref{thm: product formula}. 
So again by Remark \ref{rem: Drinfeld polynomial}, 
$\sigma_{d+1}(\widetilde{V}) \neq 0$ if and only if 
$\sigma_d(V) \neq 0$ and $\sigma_1(V(1,a)) \neq 0$. 
By Theorem \ref{thm: sigma}, observe also 
that $\sigma_d(V) \neq 0$ holds if the $\T$-module V is irreducible. 

\begin{proposition}
\label{prop: irr}
Assume that a $\T$-module V via $\varphi_s$ is irreducible and 
has diameter $d$. Assume also that the $\T$-module 
$\widetilde{V}=V \otimes V(1,a)$ via $\varphi_s$ 
satisfies $\sigma_{d+1}(\widetilde{V}) \neq 0$. 
If the $\T$-module $\widetilde{V}$ 
via $\varphi_s$ has a nonzero $\T$-submodule $W$ 
that does not contain the highest weight space $\widetilde{U}_0$ of 
$\widetilde{V}$, then the Drinfel'd polynomial 
$P_V(\lambda)$ of the $\T$-module $V$ via $\varphi_s$ vanishes at 
$\lambda=-a\,q^2 - \varepsilon\, \varepsilon^* a^{-1}q^{-2}$: 
$$ P_V (-a\, q^2 - \varepsilon\, \varepsilon^* a^{-1} q^{-2})=0, $$
where we understand $\varepsilon^* a^{-1}=0$ 
if $(\varepsilon,\varepsilon^*)=(1,0)$ and $a=0$. 
\end{proposition}

The remainder of this section is devoted to the proof of 
Proposition \ref{prop: irr}. 
Without loss of generality, we can assume that $W$ is irreducible as a 
$\T$-submodule, since we may replace $W$ by a minimal $\T$-submodule 
contained in $W$. Let 
\begin{eqnarray*}
\widetilde{V} &=& \bigoplus_{i=0}^{d+1} \widetilde{U}_i, \\
W &=& \bigoplus_{i=r}^{r+d'} U_i (W) ~~~~~ 
( U_i (W) \subseteq \widetilde{U_i})
\end{eqnarray*}
denote the weight-space decompositions of $\widetilde{V}$, $W$ respectively, 
where $d'$ is the diameter of $W$. Since $W$ is irreducible, 
the highest weight space $U_r(W)$ has dimension 1 by 
Theorem \ref{thm: shape} and 
so is spanned by a nonzero vector $w_0$: 
$$ U_r (W) = \langle w_0 \rangle. $$ 
Since $W \nsupseteq \widetilde{U}_0$, we have $r \geq 1$. 
Since $x\, U_i(W) \subseteq U_{i+1}(W)$ and 
$y\, U_i(W) \subseteq U_{i-1}(W)$, 
where $u_{r-1}(W)=0,\, U_{r+d'+1}(W)=0$, we have 
\begin{eqnarray}
\label{yw0}
& & y\, w_0=0, \\
\label{yixiw0}
& & y^i x^i w_0=\sigma_i(W)w_0
\end{eqnarray}
for $i=0,\,1,\,2,\,\cdots$. Let 
$$ V=\bigoplus^{d}_{i=0} U_i $$
denote the weight space decomposition of $V$ and let 
$$ V(1,a)=\langle v_0,v_1\rangle $$
be a standard basis of $V(1,a)$. Then
$$ \widetilde{U}_i=U_i \otimes 
\langle v_0\rangle + U_{i-1} \otimes \langle v_1\rangle$$
for $0 \leq i \leq d+1$, where $U_{-1}=U_{d+1}=0$. In particular 
\begin{eqnarray} 
\label{w0}
 w_0 = u_r \otimes v_0 + u_{r-1} \otimes v_1 
\end{eqnarray}
for some $u_r \in U_r,\, u_{r-1} \in U_{r-1}$. 
\begin{lemma}
\label{lemma: ci(m)}
For i, m $\in \Z$, set 
\begin{eqnarray*}
c_i(m) &=& a\, s\, q^{i+2m-d-1} + \varepsilon\, s^{-1} q^{-i-2m+d+1} ,\\
c_i^*(m) &=& \varepsilon^* a^{-1} s\, q^{-i+2m-d+1} + s^{-1} q^{i-2m+d-1}.
\end{eqnarray*}
Then for $1 \leq i$, the following $(i)\sim(v)$ hold. 
\begin{enumerate}
\item[$(i)$] 
$y\, u_{r-1} = 0$, \\
$y\, u_r = -c_1^* (r-1)\, u_{r-1}$. 
\item[$(ii)$] 
$\sigma_i(W)\, u_{r-1} = 
y^i x^i u_{r-1} + \alpha\, q\, [\,i\,]\, c_i (r)\, y^i x^{i-1} u_r. $
\item[$(iii)$] 
$\sigma_i(W)\, u_r =
 y^i x^i u_r + 
\alpha\, q\,[\,i\,]^2 c_i(r)\, c_{-i}^* (r-1)\, y^{i-1} x^{i-1} u_r 
 + [\,i\,]\, c_{-i}^* (r-1)\, y^{i-1} x^i u_{r-1}.$ 
\item[$(iv)$] 
$y^{i+1} x^i u_r = -[i+1]\, c_{-(i-1)}^* (r-1)\, \sigma_i(W) u_{r-1}.$ 
\item[$(v)$] 
$y^i x^i u_{r-1} = \sigma_i(W) u_{r-1} + 
\alpha\, q\,[\,i\,]^2 c_i(r)\, c_{-(i-2)}^* (r-1)\, \sigma_{i-1}(W)\, u_{r-1}.$
\end{enumerate}
\end{lemma}
Proof. 
By Lemma \ref{lemma: V otimes V(1,a)}, 
we have for $w_0=u_r \otimes v_0 + u_{r-1} \otimes v_1$, 
$$y\,w_0=(y\,u_r) \otimes v_0 + (y\,u_{r-1}) 
\otimes v_1 + c_1^*(r-1)\, u_{r-1} \otimes v_0. $$
Since $y\,w_0=0$ and $y\,u_r \in U_{r-1}$, we obtain 
$y\,u_r+c_1^*(r-1)\,u_{r-1}=0,\;y\,u_{r-1}=0$ and (i) holds.
\medskip
Again by Lemma \ref{lemma: V otimes V(1,a)}, we have 
\begin{eqnarray*}
x^i w_0 &=& (x^iu_r) \otimes v_0 + \alpha\, q\,[\,i\,]\, c_i(r)\,(x^{i-1}u_r) 
\otimes v_1 + (x^i u_{r-1}) \otimes v_1, \\
y^i x^i w_0 &=& (y^i x^i u_r) \otimes v_0 \\
& & + \alpha\, q\,[\,i\,]\, c_i(r)\,\bigl((y^i x^{i-1}u_r) \otimes v_1 + 
[\,i\,]\, c_i^* (r+i-1)\,(y^{i-1}x^{i-1} u_r) \otimes v_0\bigr) \\
& & + (y^i x^i u_{r-1}) \otimes v_1 + 
[\,i\,]\, c_i^* (r+i-1)\,(y^{i-1}x^i u_{r-1}) \otimes v_0. 
\end{eqnarray*}
Since $y^ix^iw_0=\sigma_i(W)w_0$ and $y^ix^iu_r,\,y^{i-1}x^{i-1}u_r,\,
y^{i-1}x^i u_{r-1} \in U_r,\; y^i x^{i-1}u_r,\, y^i x^i u_{r-1} \in U_{r-1}$, 
we obtain 
\begin{eqnarray*}
\sigma_i(W)\,u_r &=& y^i x^i u_r + 
\alpha\, q\,[\,i\,]^2 c_i(r)\, c_i^* (r+i-1)\, y^{i-1}x^{i-1} u_r \\
& & + [\,i\,]\, c_i^* (r+i-1)\, y^{i-1}x^i u_{r-1}, \\
\sigma_i(W)\, u_{r-1} &=& 
\alpha\, q\,[\,i\,]\, c_i(r)\, y^i x^{i-1} u_r + y^i x^i u_{r-1}.
\end{eqnarray*}
Since $c_i^*(r+i-1)=c_{-i}^*(r-1)$, (ii) and (iii) hold. 

\medskip
By (ii) and (iii) we obtain 
$$\sigma_i(W)\, y\, u_r - 
[\,i\,]\, c_{-i}^*(r-1)\, \sigma_i(W)\, u_{r-1}= y^{i+1}x^i u_r.$$ 
Since $y\, u_r=-c_1^*(r-1)\,u_{r-1}$ by (i) and 
$c_1^*(r-1)+[\,i\,]\,c_{-i}^*(r-1)=
[i+1]\,c^*_{-(i-1)}(r-1)$, we have (iv). 
Observe $y^i x^{i-1} u_r = -[\,i\,]\,c_{-(i-2)}^* (r-1)\, 
\sigma_{i-1}(W)\, u_{r-1}$ 
is valid for $i \geq 1$ by (iv), (i) 
and put this identity into (ii) to obtain (v). 
\hfill $\Box$

\medskip 
\begin{lemma}
\label{lemma: r}
It holds that $(i)$ $u_{r-1} \neq 0$, $(ii)$ $u_r \neq 0$ and $(iii)$ $r=1$. 
\end{lemma}
Proof. 
Suppose $u_{r-1}=0$. Then by Lemma \ref{lemma: ci(m)} (iii), 
$$ y^i x^i u_r = 
- \alpha\, q\,[\,i\,]^2 c_i(r)\, c_{-i}^*(r-1)\, y^{i-1} x^{i-1} u_r + 
\sigma_i(W)\, u_r $$
for $1 \leq i$, so we have $y^i x^i u_r \in \langle u_r \rangle$ 
for $0 \leq i$ by induction on $i$. Moreover $y\,u_r=0$
by Lemma \ref{lemma: ci(m)} $(i)$. 
Since $\T$ is spanned by 
$k^n \omega_{\lambda}(x,y) ~(n \in \Z,\, \lambda \in \Lambda^{irr})$ 
by Theorem \ref{thm: T basis}, 
it follows from $y^i x^i u_r \in \langle u_r\rangle ~(0 \leq i)$ 
and $y\, u_r=0$ that 
$$ \T u_r \subseteq \bigoplus_{r \leq i} U_i .$$
Since $1 \leq r, \T u_r$ is a proper $\T$-submodule of $V$.
This contradicts the irreducibility of $V$. Thus (i) holds. 

\medskip
Suppose $u_r=0$. Then by Lemma \ref{lemma: ci(m)} (ii), 
$$y^i x^i u_{r-1} = \sigma_i (W) u_{r-1}.$$
Since $\T$ is spanned by 
$k^n \omega_{\lambda}(x,y) ~(n \in \Z,\, \lambda \in \Lambda^{irr})$, 
it follows from $y^i x^i u_{r-1} 
\in \langle u_{r-1}\rangle ~(0 \leq i)$ and 
$y\, u_{r-1}=0$ that 
$$ \T u_{r-1} \subseteq \bigoplus_{r-1 \leq i} U_i . $$
So $V$ has a nonzero $\T$-submodule contained in 
$\bigoplus_{r-1 \leq i} U_i$. 
Since $V$ is irreducible, we obtain $r-1=0$, i.e., 
$w_0=u_0 \otimes v_1 ~(u_0 \neq 0, u_1=0)$. 
By Lemma \ref{lemma: ci(m)} (i), 
$y\, u_1=-c_1^*(0)\, u_0$ and so $c_1^*(0)=0$. 
By Lemma \ref{lemma: ci(m)} (iii)  with $i=1$, 
$c_{-1}^*(0)\,x\, u_0=0$. Note $x\, u_0 \neq 0$, otherwise 
$\sigma_d(V)=0$, which contradicts the assumption that $V$ is 
irreducible as a $\T$-module. Thus $c_{-1}^*(0)=0$. 
From $c_1^*(0)=0$ and $c_{-1}^*(0)=0$, we have 
\begin{eqnarray*}
& & \varepsilon^* a^{-1} s\, q^{-d} + s^{-1} q^d=0 ,\\
& & \varepsilon^* a^{-1} s \,q^{-d+2} + s^{-1} q^{d-2}=0.
\end{eqnarray*}
This implies $\varepsilon^*=1,\, a=-s^2 q^{-2d}=-s^2 q^{4-2d}$ and 
we have $q^4=1$. This contradicts the assumption that $q$ is not a root 
of unity. Hence (ii) holds. 

\medskip
By Lemma \ref{lemma: ci(m)} (i), (v), 
we have $y\, u_{r-1}=0$ and 
$y^i x^i u_{r-1} \in \langle u_{r-1}\rangle ~(0 \leq i)$. 
The same argument of the previous paragraph is valid and 
$V$ has a nonzero $\T$-submodule $\T u_{r-1}$ 
contained in $\bigoplus_{r-1 \leq i} U_i$. Hence we obtain 
$r-1=0$, i.e., (iii) holds. 
\hfill $\Box$ 

\begin{lemma}
\label{lemma: fi}
For $0 \leq i$, 
$$ \sigma_i(W) = f_i\, \sigma_{i-1} (W) + \sigma_i (V), $$
where
$$ f_i=(q^i-q^{-i})^2 
\bigl(\varepsilon s^{-2} q^{2(d-i)}+\varepsilon^* s^2 q^{-2(d-i)} 
+ a\, q^2 + \varepsilon\, \varepsilon^* a^{-1} q^{-2}\bigr).$$
\end{lemma}

\noindent
Proof. 
Since $r=1$, we have $y^i x^i u_{r-1}= 
\sigma_i(V)\, u_{r-1}$. 
By Lemma \ref{lemma: r}, $u_{r-1} \neq 0$. 
By Lemma \ref{lemma: ci(m)} $(v)$, we obtain 
\begin{eqnarray*}
\sigma_i (V) 
&=& \sigma_i(W) + \alpha\, q\,[\,i\,]^2 c_i(1)\, c_{-(i-2)}^* (0)\, 
\sigma_{i-1}(W) \\
&=& \sigma_i(W) - f_i\, \sigma_{i-1} (W). 
\end{eqnarray*}
\hfill $\Box$

\medskip
\begin{lemma}
\label{lemma: d'}
It holds that $d'=d-1$, 
where $d$, $d'$ are the diameters of $V$, $W$ respectively. 
\end{lemma}
Proof. 
Obviously 
$\sigma_i(W)=\sigma_{i-1}(W)=0$ for $d'+2 \leq i$. 
So we have 
$\sigma_i(V)=0$ for $d'+2 \leq i$ by Lemma \ref{lemma: fi}. 
This implies $d \leq d'+1$, 
since $\sigma_d(V) \neq 0$. 
On the other hand, 
the weight-space decompositions of 
$\widetilde{V}$, $W$ are 
$\widetilde{V}=\widetilde{U}_0 + \cdots + \widetilde{U}_{d+1}$, 
$W=U_r(W) + \cdots + U_{r+d'}(W)$ $ (r=1)$ with 
$U_i(W) \subseteq \widetilde{U}_i$. 
So $r+d' \leq d+1$, i.e., $d' \leq d$. 
Therefore either $d=d'+1$ or $d=d'$. 

\medskip
Suppose $d=d'$. Then $0 \neq U_{d+1}(W) \subseteq \widetilde{U}_{d+1}$. 
Since $\widetilde{V}$ is a tensor product of evaluation modules, 
we generally have dim $\widetilde{U}_0=dim\, \widetilde{U}_{d+1}=1$. 
So $U_{d+1}(W)=\widetilde{U}_{d+1}$, in particular $\T \widetilde{U}_{d+1}$ 
is contained in $W$. On the other hand, we assumed 
$\sigma_{d+1}(\widetilde{V}) \neq 0$ for Proposition \ref{prop: irr}. 
This implies 
$\T \widetilde{U}_{d+1} \supseteq \widetilde{U}_0$. 
Therefore $W$ contains $\widetilde{U}_0$, which is a contradiction. 
\hfill $\Box$

\bigskip
\noindent
{\bf Proof of Proposition \ref{prop: irr}}.
Set $\sigma_i=\sigma_i(V)$. Using Lemma \ref{lemma: fi} repeatedly, 
we have for $0 \leq i$
\begin{eqnarray*}
\sigma_i(W) &=& (f_i\, f_{i-1} \cdots f_1)\, \sigma_0 + 
(f_i\, f_{i-1} \cdots f_2)\, \sigma_1 + \cdots + f_i\, \sigma_{i-1} + \sigma_i.
\end{eqnarray*}
By Lemma \ref{lemma: d'}, $d'=d-1$. So $\sigma_d(W)=0$. Thus 
\begin{eqnarray}
\label{fi}
\sum_{j=0}^d f_d \,f_{d-1} \cdots f_{j+1}\, \sigma_j=0. 
\end{eqnarray}
Define the polynomial $f_i (\lambda)$ of degree 1 in $\lambda$ by 
$$
f_i(\lambda)=(q^i-q^{-i})^2 (\varepsilon s^{-2} q^{2(d-i)} 
+ \varepsilon^* s^2 q^{-2(d-i)} - \lambda)
$$
for $1 \leq i$. Then by definition 
$$ P_V(\lambda)=Q^{-1} \sum_{i=0}^d \sigma_i\, f_{i+1}
 (\lambda) \cdots f_d (\lambda), $$
where $Q={(-1)}^d {(q-q^{-1})}^2 {(q^2 - q^{-2})}^2 \cdots {(q^d-q^{-d})}^2$.
Since $f_i=f_i (-a\,q^2-\varepsilon\, \varepsilon^* a^{-1}q^{-2})$, 
we have by (\ref{fi}) 
$$ P_V (-a\,q^2-\varepsilon\, \varepsilon^* a^{-1}q^{-2})=0. $$
This completes the proof of Proposition \ref{prop: irr}. 
\hfill $\Box$

\section{Construction of finite-dimensional irreducible \\
$\T$-modules via $\varphi_s$}
In this section, we prove Theorem \ref{thm: 1st kind}, 
Theorem \ref{thm: 2nd kind}, Theorem \ref{thm: 3rd kind}. 
Each theorem consists of three parts (i), (ii), (iii). 
The second part (ii) immediately follows from 
Theorem \ref{thm: sigma}' and Theorem \ref{thm: product formula}. 
We only prove (i) and (iii) for each of the theorems. 
Throughout this section, $s$ stands for a nonzero scalar of 
$\C$ chosen arbitrarily. 

%
\medskip
For the augmented TD-algebra $\T=\T_q^{(\varepsilon,\varepsilon^*)}$, 
we consider the following $\T$-module $V$ via $\varphi_s$ (see Section 1.4): 
if $(\varepsilon,\varepsilon^*)=(1,1)$ or $(0,0)$, 
$$ V=V(\ell_1,a_1) \otimes \cdots \otimes V(\ell_n,a_n), $$
where $1 \leq n,\, 1 \leq \ell_i,\, a_i \neq 0 ~(1 \leq i \leq n)$, 
and if $(\varepsilon,\varepsilon^*)=(1,0)$, 
$$ V=V(\ell) \otimes V(\ell_1,a_1) \otimes \cdots \otimes V(\ell_n,a_n), $$
where $0 \leq n,\, 0 \leq \ell,\, 
1 \leq \ell_i,\, a_i \neq 0 ~(1 \leq i \leq n)$.
With such a $\T$-module $V$ via $\varphi_s$, we associate the multi-set 
$ {\{ S(\ell_i,a_i) \}}_{i=1}^n$ of $q$-strings, 
where 
$$ S(\ell_i,a_i)=
\{a_i q^{-\ell_i+1} , a_i q^{-\ell_i+3}, \cdots, a_i q^{\ell_i-1}\}.$$
Consider a $\T$-module $V'$ via $\varphi_s$ of the same kind: 
\begin{eqnarray*}
V' &=& V(\ell_1',a_1') \otimes \cdots \otimes V(\ell_m',a_m') ~~~~~~~~~~
~~~~~ \mbox{if} ~(\varepsilon,\varepsilon^*)=(1,1),\, (0,0), \\
V' &=& V(\ell') \otimes V(\ell_1',a_1') \otimes \cdots \otimes V(\ell_m',a_m') 
~~~~~ \mbox{if} ~(\varepsilon,\varepsilon^*)=(1,0).
\end{eqnarray*}
For $(\varepsilon,\varepsilon^*)=(1,1)$, such $\T$-modules $V,\,V'$ via 
$\varphi_s$ are said to be {\it equivalant} if the associated multi-sets of 
$q$-strings are equvalent, i.e., $m=n$ and there exist 
$\varepsilon_i \in \{ 1,-1 \} ~(1 \leq i \leq n)$ such that 
$S(\ell_i', a_i')=S(\ell_i,a_i^{\varepsilon_i}) ~(1 \leq i \leq n)$ 
with a suitable rearrangement of the ordering of 
$S(\ell_1', a_1'), \cdots, S(\ell_n', a_n')$. 
For $(\varepsilon,\varepsilon^*)=(0,0)$, such $\T$-modules $V,\,V'$ via 
$\varphi_s$ are said to be {\it equivalent} if $m=n$ and 
$S(\ell_i', a_i')=S(\ell_i, a_i) ~(1 \leq i \leq n)$ 
with a suitable rearrangement of the ordering of
$S(\ell_1', a_1'), \cdots, S(\ell_n', a_n')$. 
For $(\varepsilon,\varepsilon^*)=(1,0)$, such $\T$-modules $V,\,V'$ via 
$\varphi_s$ are said to be {\it equivalent} if $\ell=\ell',\, m=n$ and 
$S(\ell_i', a_i')=S(\ell_i, a_i) ~(1 \leq i \leq n)$ with  
a suitable rearrangement of the ordering of 
$S(\ell_1', a_1'), \cdots, S(\ell_n', a_n')$. 
\begin{lemma}
\label{lemma: equivalent module}
If a $\T$-module $V$ via $\varphi_s$ is irreducible, then every $\T$-module 
$V'$ via $\varphi_s$ that is equivalent to $V$ is isomorphic to $V$ as 
$\T$-modules via $\varphi_s$, in particular irreducible. 
\end{lemma}
Proof. 
Since $V$ and $V'$ are equivalent, 
$V$ and $V'$ have the same 
Drinfel'd polynomial by Theorem \ref{thm: product formula}. 
In particular $\sigma_d(V)=\sigma_d(V')$, 
where $d$ is the diameter of the $\T$-modules $V$, $V'$.
Let $U_0'$ denote the highest weight space of the $\T$-module $V'$ 
via $\varphi_s$. Set $W=\T U_0'$ and let $M$ be a maximal $\T$-submodule 
of $W$. Then $V'$ and $W/M$ have the same Drinfel'd polynomial 
by Lemma \ref{lemma: PbarW}. 
Hence $V$ and $W/M$ have the same Drinfel'd polynomial. 
By Theorem $\ref{thm: sigma}'$, the irreducible $\T$-modules $V$, 
$W/M$ are isomorphic, in particular ${\rm dim}\, V= {\rm dim}\, W/M$. 
But $V$ and $V'$ are equivalent, 
in particular ${\rm dim}\, V= {\rm dim}\, V'$.
Thus ${\rm dim}\, V'= {\rm dim}\, W/M$ and we have $V'=W,\, M=0$. 
This means that $V$ and $V'$ are isomorphic as $\T$-modules via $\varphi_s$.  
\hfill $\Box$

\medskip
\noindent
{\bf Proof of the `only if' part of (i)}.  
The `only if' part of 
Theorem \ref{thm: 1st kind} (i) 
follows from Lemma \ref{lemma: equivalent module}. 
Suppose $-s^2 \in S(\ell_i,a_i) \cup S(\ell_i,a_i^{-1})$ 
for some $i ~(1 \leq i \leq n)$. Then  
$P_V (\varepsilon s^{-2}+ \varepsilon^* s^2)=0$ 
by Theorem \ref{thm: product formula}, 
which contradicts the irreducibility of $V$ 
by Theorem $\ref{thm: sigma}'$. 
Suppose the multi-set 
${\{ S(\ell_i,a_i)\}}^n_{i=1} $
of $q$-strings is not strongly in general position. 
Then there exists a multi-set 
${\{ S(\ell_i',a_i')\}}^m_{i=1}$
of $q$-strings that is eqivalent to ${\{ S(\ell_i,a_i)\}}^n_{i=1}$ 
and not in general position. Set 
$V' = V(\ell_1',a_1') \otimes \cdots \otimes V(\ell_m',a_m')$. 
Since ${\{ S(\ell_i',a_i')\}}^m_{i=1}$ is not in general position, 
$V'$ is not irreducible 
as an $\LL$-module, consequently as a $\T$-module via $\varphi_s$. 
Sine ${\{ S(\ell_i',a_i')\}}^m_{i=1}$ is eqivalent to 
${\{ S(\ell_i,a_i)\}}^n_{i=1}$ and $V$ is irreducible as 
a $\T$-module via $\varphi_s$, $V'$ is  
also irreducible as a $\T$-module via $\varphi_s$ 
by Lemma \ref{lemma: equivalent module}. Thus we get a contradiction.  
The `only if' part of 
Theorem \ref{thm: 2nd kind} (i), Theorem \ref{thm: 3rd kind} (i) 
can be proved similarly. 
\hfill $\Box$

\subsection{Proof of the `if' part of Theorem \ref{thm: 1st kind} (i)}
We start with an observation of the $U_q({sl}_2)$-loop algebra $\LL$. 
\begin{lemma}
\label{lemma: tau}
There exists an algebra anti-homomorphism of $\LL$ that sends 
$e_i^+,\, e_i^- k_i,\, k_i,\, k_i^{-1}$ 
to $e_i^- k_i,\, e_i^+,\, k_i,\, k_i^{-1}$, respectively $(i=0,1)$. 
Such an anti-homomorphism is unique and we denote it by $\tau$: 
$$ \tau: \LL \longrightarrow \LL ~~~  (e_i^+,\, e_i^-k_i,\,  k_i^{\pm 1} 
\mapsto e_i^-k_i,\, e_i^+,\, k_i^{\pm 1}~ \mbox{respectively} ). $$ 
It holds that 
$\tau^2 =1$ and $(\tau \otimes \tau)\, \Delta = \Delta\, \tau$, 
where $\Delta: \LL \longrightarrow \LL \otimes \LL$ 
is the coproduct of $\LL$. 
\end{lemma}

The assertions of Lemma \ref{lemma: tau} can be checked 
by straightforward calculations. 

\medskip
For an $\LL$-module $V$, the dual vector space of $V$ 
$$ \mbox{Hom}\,(V,\C) = 
\{ f: V \longrightarrow \C ~|~\mbox{f is a linear mapping} \} $$
becomes an $\LL$-module by 
$$ (Xf)\,(v) = f\, (\tau(X)\, v) \qquad \qquad (v \in V) $$
for $f \in \mbox{Hom}\,(V,\C)$, $X \in \LL$. For $\LL$-modules $V$, $V'$, 
we identify 
$\mbox{Hom}\,(V \otimes V', \C)$ with 
$\mbox{Hom}\,(V,\C) \otimes \mbox{Hom}\,(V',\C)$ as vector 
spaces by 
$$ (f \otimes g)\,(v \otimes v')= f(v)\, g(v'). $$
It can be easily checked by the relation 
$(\tau \otimes \tau)\, \Delta = \Delta\, \tau$ that this identification 
gives an $\LL$-module isomophism 
$$ \mbox{Hom}\,(V \otimes V', \C) \simeq 
\mbox{Hom}\,(V,\C) \otimes \mbox{Hom}\,(V',\C), $$
where $\LL$ acts on $V \otimes V'$ and  
$\mbox{Hom}\,(V,\C) \otimes \mbox{Hom}\, (V',\C)$ via 
the coproduct $\Delta$. 
\begin{lemma}
\label{lemma: Hom}
For evaluation modules, we have the following isomorphisms as $\LL$-modules. 
\begin{eqnarray*}
&& (i) ~~
{\rm Hom}\, (V(\ell,a), \C) \simeq V(\ell,a^{-1}). \\
&& (ii) ~~
{\rm Hom}\, \bigl(V(\ell_1,a_1) \otimes \cdots \otimes V(\ell_n,a_n), \C\bigr) 
\simeq V(\ell_1,a_i^{-1}) \otimes \cdots \otimes V(\ell_n,a_n^{-1}). 
~~~~~~~~~~~~~~~~~~~~~~~~~~~~~~~~~~~~~~ 
\end{eqnarray*}
\end{lemma}
Proof. 
Let $V(\ell,a)=\langle v_0,\,v_1, \cdots, v_{\ell} \rangle$ be a standard basis and 
${\rm Hom}\,(V(\ell,a),\C)=\langle f_0,\,f_1, \cdots, f_l\rangle$ 
the dual basis: $f_i\,(v_j)= \delta_{ij}$. Set 
$$ g_i = q^{-i(\ell-i+1)}\binom{\ell}{i}\, f_i,$$
where 
$\binom{\ell}{i}=[\ell] ! / [\ell-1] !\, [\,i\,] !$, 
the $q$-binomial coefficient. Then we have 
$e_0^+ g_i = a^{-1} q\, [i+1]\, g_{i+1}$, 
$e_0^- g_i=a\, q^{-1} [\ell-i+1]\, g_{i-1},\; 
e_1^+ g_i = [\ell-i+1]\, g_{i-1},\; 
e_1^- g_i = [i+1] \,g_{i+1},\; k_0 g_i= q^{2i-\ell} g_i$, where 
$g_{-1}=g_{\ell+1}=0$. So if 
$V(\ell,a^{-1})=\langle w_0,w_1,\cdots,w_{\ell}\rangle$ is a standard basis, 
then ${\rm Hom}\, (V(\ell,a),\C)$ is isomorphic to 
$V(\ell,a^{-1})$ as $\LL$ -modules by 
the correspondence of $g_i$ to $w_i ~(0 \leq i \leq \ell)$.
Part (ii) follows from part (i) and the fact 
${\rm Hom}\,(V \otimes V', \C) \simeq 
{\rm Hom}\,(V,\C) \otimes {\rm Hom}\,(V',\C)$ as $\LL$-modules. 
\hfill $\Box$ 

\medskip
We now prove the `if' part of 
Theorem \ref{thm: 1st kind} (i). 
Namely, we are in the case of 
$(\varepsilon,\varepsilon^*)=(1,1)$ and given a $\T$-module 
$$ V =V(\ell_1,a_1) \otimes \cdots \otimes V(\ell_n,a_n) $$
via $\varphi_s$ such that 
$-s^2 \notin S(\ell_i,a_i) \cup S(\ell_i,a_i^{-1}) ~~(1 \leq i \leq n)$ 
and the multi-set 
$\{S(\ell_i,a_i) \}_{i=1}^n$ of $q$-strings is 
strongly in general position. We want to show that $V$ is irreducible 
as a $\T$-module via $\varphi_s$. 
Observe that the ordering of the tensor product does not change 
the isomorphism class of $V$ as an $\LL$-module and consequently 
as a $\T$-module via $\varphi_s$, since the multiset 
$\{S(\ell_i,a_i) \}_{i=1}^n$
of $q$-strings is in general position. 
First we show 
\begin{lemma}
\label{lemma: isom}
For any choice of 
$\varepsilon_i \in \{ 1,-1 \} ~(1 \leq i \leq n)$, 
$V$ is isomorphic to 
$$V(\ell_1,a_1^{\varepsilon_1}) \otimes \cdots \otimes 
V(\ell_n,a_n^{\varepsilon_n})$$ 
as $\T$-modules via $\varphi_s$.
\end{lemma}
Proof.  
We proceed by induction on $n$.
First let  $n=1$. Then 
$V(\ell_1,a_1)$, $V(\ell_1,a_1^{-1})$ have the same Drinfel'd 
polynomial by Theorem \ref{thm: product formula} and
the Drinfel'd polynomial does not vanish 
at $s^2+s^{-2}$, since $-s^2 \notin S(\ell_1,a_1) \cup S(\ell_1,a_1^{-1})$. 
So $V(\ell_1,a_1)$, $V(\ell_1,a_1^{-1})$ have 
nonzero $\sigma_{\ell_1}$ by Remark \ref{rem: Drinfeld polynomial}. 
In the case of evaluation modules, 
the property $\sigma_{\ell_1} \neq 0$ 
implies the irreducibility of the $\T$-modules via $\varphi_s$. 
Thus $V(\ell_1,a_1)$, $V(\ell_1,a_1^{-1})$ 
are irreducible as $\T$-modules via $\varphi_s$ and have the same Drinfel'd 
polynomial that does not vanish at $s^2+s^{-2}$. 
By Theorem $\ref{thm: sigma}'$,  
$V(\ell_1,a_1)$, $V(\ell_1,a_1^{-1})$  are isomorphic 
as $\T$-modules via $\varphi_s$. 

\medskip
For $n \geq 2$, set 
$V'=V(\ell_1,a_1) \otimes \cdots \otimes V(\ell_{n-1},a_{n-1})$. 
Then by induction on $n$, $V'$ is isomorphic to 
$V''=V(\ell_1,a_1^{\varepsilon_1}) \otimes \cdots \otimes 
V(\ell_{n-1},a_{n-1}^{\varepsilon_{n-1}})$ as $\T$-modules via $\varphi_s$. 
Let $\psi : V' \longrightarrow V''$ denote an isomorphism between 
the $\T$-modules $V',\,V''$ via $\varphi_s$. 
The generators $x,\,y,\,k,\,k^{-1}$ of $\T$ are mapped by $\varphi_s$ 
to 
$x(s)=\alpha(s\, e_0^+ + s^{-1} e_1^- k_1),\, 
y(s)= s\, e_0^- k_0 + s^{-1} e_1^+,\, s k_0,\, s^{-1} k_1$, respectively, 
and those elements of $\LL$ are mapped by $\Delta$ to 
\begin{eqnarray*}
\Delta \, \bigl(x(s)\bigr) &=& 
x(s) \otimes 1 + \alpha\, s\, k_0 \otimes e_0^+ 
+ \alpha\, s^{-1} k_1 \otimes e_1^- k_1,\\
\Delta\, \bigl(y(s)\bigr) &=&  
y(s) \otimes 1 + s\, k_0 \otimes e_0^- k_0 + s^{-1} k_1 \otimes e_1^+,\\
\Delta\, (s\, k_0)&=&s\, k_0 \otimes k_0,\, \\
\Delta (s^{-1}k_1)&=&s^{-1}k_1 \otimes k_1, 
\end{eqnarray*}
respectively.  It can be easily checked that 
the vector-space isomorphism 
$$\psi \otimes id: 
V' \otimes V(\ell_n, a_n) \longrightarrow V'' \otimes V(\ell_n, a_n)$$ 
commutes with the action of each of the elements 
$\Delta \, \bigl(x(s)\bigr),\,\Delta\, \bigl(y(s)\bigr),\, 
\Delta\, (s\, k_0),\, \Delta (s^{-1}k_1)$.  
So we get 
$$ V' \otimes V(\ell_n, a_n) \simeq V'' \otimes V(\ell_n, a_n)$$ 
as $\T$-modules via $\varphi_s$. 
Since ${\{ S(\ell_i,a_i^{\varepsilon_i}) \}}_{i=1}^{n-1} 
\cup \{S(\ell_n,a_n)\}$  
is in general position, 
$$ V'' \otimes V(\ell_n,a_n) \simeq V(\ell_n,a_n) \otimes V'' $$
as $\LL$-modules and consequently as $\T$-modules via $\varphi_s$. 
By the same argument, we have 
\begin{eqnarray*}
V(\ell_n,a_n) \otimes V'' 
&\simeq& V(\ell_n,a_n^{\varepsilon_n}) \otimes V'' \\
&\simeq& V'' \otimes V(\ell_n,a_n^{\varepsilon_n}) 
\end{eqnarray*}
as $\T$-modules via $\varphi_s$. Thus 
$V' \otimes V(\ell_n,a_n) \simeq V'' \otimes V(\ell_n,a_n^{\varepsilon_n})$ as 
$\T$-modules $\varphi_s$ and the proof is completed.
\hfill $\Box$ 


\medskip
Next we intoroduce a partial ordering on $\C \backslash \{ 0 \}$ by 
\begin{eqnarray}
\label{partial ordering}
a \leq b ~\Leftrightarrow~ 
b = a\, q^{2i} ~\mbox{for some integer}~ i \geq 0.
\end{eqnarray}
Consider $i_0 ~(1 \leq i_0 \leq n)$ such that $a_{i_0} q^{\ell_{i_0}-1}$ or 
$a_{i_0}^{-1} q^{\ell_{i_0}-1}$ is maximal 
with respect to the partial ordering on the set of nonzero scalars 
$a_i q^{\ell_{i}-1},\, a_{i}^{-1} q^{\ell_{i}-1} ~(1 \leq i \leq n)$.
Among such $i_0's$, choose one for which $\ell_{i_0}$ is smallest. 
Since the ordering of the tensor product does not matter 
about the isomorphism class of $V$ as a $\T$-module via $\varphi_s$,   
we may assume that $i_0=n$, and by Lemma \ref{lemma: isom} that 
$a_n q^{\ell_n -1}$ is 
maximal among $a_i^{\pm} q^{\ell_i-1} ~(1 \leq i \leq n)$. So 
\begin{eqnarray}
\label{ellq max}
&& a_n q^{\ell_n +1} \notin 
S(\ell_i,a_i) \cup S(\ell_i,a_i^{-1})~~ (1 \leq i \leq n), \\
\label{ell min}
&& \ell_n \leq \ell_i 
~~{\rm if}~ a_n q^{\ell_n -1} \in 
S(\ell_i,a_i) \cup S(\ell_i,a_i^{-1}).
\end{eqnarray}

\medskip
We proceed by induction on dim $V$ to prove 
the `if' part of Theorem \ref{thm: 1st kind} (i). Set 
$$V' = V(\ell_1,a_1) \otimes \cdots \otimes V(\ell_{n-1}, a_{n-1}).$$ 
Then $V=V' \otimes V(\ell_n,a_n)$. Since the $\LL$-module $V(\ell_n,a_n)$ 
is by Lemma \ref{lemma: ev module} embedded in the $\LL$-module 
$V(\ell_n -1,a_n q^{-1}) \otimes V(1,a_n q^{\ell_n -1})$ 
as the $\LL$-submodule spanned by the highest weight space, 
the $\LL$-module $V$ can be regarded as embedded in the $\LL$-module 
$$ \widetilde{V} = 
\bigl(V' \otimes V(\ell_n -1,a_n q^{-1}) \bigr)
\otimes V(1,a_n q^{\ell_n -1}).$$
We understand 
$V' \otimes V(\ell_n -1, a_n q^{-1}) = V'$ if $\ell_n =1$.
Our strategy is 
to apply Proposition \ref{prop: irr} to 
the $\T$-module $\widetilde{V}$ via $\varphi_s$. 
To do so,  
we need to check the prerequisites for it, namely 
that $V' \otimes V(\ell_n -1, a_n q^{-1})$ is irreducible 
as a $\T$-module via $\varphi_s$, and that 
$\sigma_d (\widetilde{V}) \neq 0$ holds,
where $d=\ell_1 + \cdots + \ell_n$, the diameter of the 
$\T$-module $\widetilde{V}$ via $\varphi_s$. 
To show that $V' \otimes V(\ell_n -1, a_n q^{-1})$ is irreducible 
as a $\T$-module via $\varphi_s$, it is enough to check 
the induction hypotheses for it, i.e., 
the first induction hypothesis that 
$-s^2$ is contained neither in 
$S(\ell_i,a_i^{\pm 1})~(1 \leq i \leq n-1)$ nor in 
$S(\ell_n -1, (a_n q^{-1})^{\pm 1})$, 
and 
the second induction hypothesis 
that the multi-set 
${\{ S(\ell_i,a_i) \}}_{i=1}^{n-1} \cup  \{ S(\ell_n -1,a_n q^{-1}) \}$ 
of $q$-strings associated with $V' \otimes V(\ell_n -1,a_n q^{-1})$ 
is strongly in general position. 
The first induction hypothesis is satisfied, since 
$-s^2 \notin S(\ell_i,a_i^{\pm 1})~(1 \leq i \leq n)$ was assumed 
at the beginning and it generally holds that 
$S(\ell_n -1,a_n q^{-1})=S(\ell_n,a_n) ~\backslash~ \{ a_n q^{\ell_n -1} \}$.  
The second induction hypothesis is satisfied, since the multi-set 
${\{ S(\ell_i,a_i) \}}_{i=1}^{n}$ was assumed at the beginning 
to be strongly in general position 
and $n$ was chosen to satisfy 
(\ref{ellq max}), (\ref{ell min}). 
Thus by induction on dimension, 
$V' \otimes V(\ell_n -1,a_n q^{-1})$ is irreducible as a $\T$-module 
via $\varphi_s$. 
To show $\sigma_d (\widetilde{V}) \neq 0$, it is enough to check that 
$P_{\widetilde{V}}(\lambda)$ does not vanish at $\lambda = s^{-2}+s^2$ 
(see Remark \ref{rem: Drinfeld polynomial}). 
Since $-s^2 \notin S(\ell_i,a_i^{\pm 1})~(1 \leq i \leq n)$, 
$P_{\widetilde{V}}(s^{-2}+s^2) \neq 0$ by Theorem \ref{thm: product formula}. 
We are now ready to apply Proposition \ref{prop: irr} to $\widetilde{V}$. 

\medskip
Suppose that the $\T$-module $V$ via $\varphi_s$ has 
a nonzero $\T$-submodule $W$ that does not contain the highest weight 
space of $V$.  Embed $V$  
into the $\T$-module $\widetilde{V}$ via $\varphi_s$ 
as a $\T$-submodule in such a way that 
$V$ and $\widetilde{V}$ share the highest weight space in common 
(see the last pragragh).  
Then by Proposition \ref{prop: irr}, 
the Drinfel'd polynomial of 
$V' \otimes V(\ell_n -1,a_n q^{-1})$ vanishes at 
$-a_n q^{\ell_n +1} - a_n^{-1} q^{-\ell_n-1}$. 
This contradicts (\ref{ellq max}) by Theorem \ref{thm: product formula}. 
Therefore we conclude that 
every nonzero $\T$-submodule $W$ of the $\T$-module $V$ via $\varphi_s$ 
contains the highest weight space of $V$. 

\medskip
Finally consider the $\LL$-module Hom\,$(V,\C)$. 
By Lemma \ref{lemma: Hom}, Hom\,$(V,\C)$ is isomorphic to 
$V(\ell_1,a_1^{-1}) \otimes \cdots \otimes V(\ell_n,a_n^{-1})$ 
as $\LL$-modules. So Hom $(V,\C)$ and $V$ are  
isomorphic as $\T$-modules via $\varphi_s$ by Lemma \ref{lemma: isom}.
For a subspace $W$ of $V$, define the subspace $W^{\perp}$ of Hom~$(V,\C)$ by 
$$ W^{\perp} = \{ f \in {\rm Hom}\,(V,\C) ~|~ f(w) = 0 ~( w \in W) \}.$$
If $W$ is invariant by the action of $\T$ via $\varphi_s$, 
then so is $W^{\perp}$, 
because the action of $X \in \LL$ on Hom $(V,\C)$ is 
defined by $(Xf)\,(v)=f(\tau (X)\, v) ~(f \in {\rm Hom}\,(V,\C),\, v \in V)$ 
and $\tau (\T)=\T$ holds by $\tau\,\, (x(s))=\alpha\, y(s),\, 
\tau(y(s))=\alpha^{-1}\,x(s),\, \tau(k_0)=k_0$. 
Moreover by the proof of Lemma~\ref{lemma: Hom}, 
the highest weight 
space of Hom\,$(V,\C)$ does not vanish on the highest weight space of $V$, 
i.e., $f(v) \neq 0$ for highest weight vectors $f$,\,$v$ of Hom\,$(V,\C)$, 
$V$, respectively. Now let $W$ be a nonzero $\T$-submodule of 
the $\T$-module $V$ via $\varphi_s$.  
Then $W$ contains the highest weight space of $V$ as shown 
in the last paragraph. This implies that $W^{\perp}$ 
is a $\T$-submodule of Hom\,$(V,\C)$ via $\varphi_s$ and does not 
contain the highest weight space of Hom\,$(V,\C)$. 
Recall Hom $(V,\C)$ and $V$ are  
isomorphic as $\T$-modules via $\varphi_s$. 
Thus $W^{\perp}=0$. 
Therefore W=V and the proof of the `if' part of 
Theorem \ref{thm: 1st kind} (i) is completed.

\subsection{Proof of the `if' part of 
Theorem \ref{thm: 2nd kind} (i), Theorem \ref{thm: 3rd kind} (i)} 
We start with observations about the quantum algebra $U_q(sl_2)$. 
The quantum algebra ${\mathcal U}=U_q (sl_s)$ is the associative $\C$-algebra 
with 1 genarated by $X^{\pm}, K^{\pm}$ subject to the relations
\begin{eqnarray*}
& & K K^{-1} = K^{-1} K = 1, \\
& & K X^{\pm} K^{-1} = q^{\pm 2} X^{\pm},\\
& & [X^+, X^-] = \frac{K-K^{-1}}{q-q^{-1}}.
\end{eqnarray*}
$V(\ell)$ denotes 
the $(\ell+1)$-dimensional irreducible ${\mathcal U}$-module: 
$V(\ell)= \langle v_0,v_1, \cdots, v_{\ell}\rangle$ and 
\begin{eqnarray*}
&&K\, v_i = q^{2i-\ell}\,v_i,\\
&&X^+ v_i=[i+1]\, v_{i+1}, \\
&&X^- v_i=[\ell-i+1]\,v_{i-1},  
\end{eqnarray*} 
where  $v_{-1}=v_{\ell+1}=0$. 
We consider a finite-dimensional ${\mathcal U}$-module $V$ 
that has the following weight-space decomposition: 
\begin{eqnarray}
\label{weight space decopm}
V = \bigoplus_{i=0}^d U_i,~~~K|_{U_i} = q^{2i-d} ~(0 \leq i \leq d). 
\end{eqnarray}
Since $V$ is completely reducible, we have 
$$ V = \bigoplus_{j=0}^{[d/2]} V^{(d-2j)}, $$
where $V^{(\ell)}$ denotes the homogeneous component that is a direct sum of 
irreducible ${\mathcal U}$-modules isomorphic to $V(\ell)$; 
$V^{(\ell)}$ is allowed to be zero. 
Set 
$$ U_i^{(d-2j)} = U_i \cap V^{(d-2j)} ~~~ 
(0 \leq i \leq d,\, 0 \leq j \leq [d/2]).$$
Then 
\begin{eqnarray*}
V^{(d-2j)} &=& \bigoplus_{i=j}^{d-j} U_i^{(d-2j)} ~~~ 
(0 \leq j \leq [d/2]),\\
U_i &=& \bigoplus_{j=0}^{i'} U_i^{(d-2j)} ~~~
(0 \leq i \leq d),
\end{eqnarray*}
where $i'= \mbox{min}\, \{ i, d-i \}$. For $j \leq i < d-j$, 
the mappings 
%
\begin{eqnarray*}
\label{X+}
&&X^+:~ U_i^{(d-2j)} \longrightarrow U_{i+1}^{(d-2j)}, \\
\label{X-}
&&X^-:~ U_{i+1}^{(d-2j)} \longrightarrow U_i^{(d-2j)}
\end{eqnarray*}
are inverses each other up to a nonzero scalar multiple and 
$X^+,\, X^-$ vanish on $U_{d-j}^{(d-2j)}, U_j^{(d-2j)}$, respectively. 
In particular, 
\begin{eqnarray}
\label{ker}
 U_j^{(d-2j)} = {\rm ker}\, {(X^+)}^{d-2j+1}|_{U_j} ~~~~~~ 
(0 \leq j \leq [d/2]). 
\end{eqnarray}
\begin{lemma}
\label{lemma: X}
Let V be a finite-dimensional 
${\mathcal U}$-module that satisfy $(\ref{weight space decopm})$. 
Let $W$ be a subspace of $V$ as a vector space. 
Assume $W$ is invariant by the actions of 
$X^+$ and $K$: 
$$ K\, W \subseteq W ,~~~ X^+\, W \subseteq W. $$
If it holds that 
$$ {\rm dim}\, (W \cap U_i) = {\rm dim}\, (W \cap U_{d-i}) \quad
(0 \leq i \leq d), $$
then $X^-\,W \subseteq W$, i.e., $W$ is a ${\mathcal U}$-submodule. 
\end{lemma}
Proof. 
Set $W_i = W \cap U_i ~~ (0 \leq i \leq d)$. 
Then since $W$ is $K$-invariant, we have 
$$ W = \bigoplus_{i=0}^d W_i, $$
allowing $W_i$ to be zero. Set 
$W_i^{(d-2j)}=W_i \cap U_i^{(d-2j)} ~~ 
(0 \leq i \leq d,\, 0 \leq j \leq [d/2])$. 
We claim 
\begin{eqnarray}
\label{Wi}
W_i = \bigoplus_{j=0}^i W_i^{(d-2j)} ~~~ 
(0 \leq i \leq [d/2]).
\end{eqnarray}

The claim holds for $i=0$, since $W_0 \subseteq U_0 = U_0^{(d)}$. 
Suppose the claim holds up to $i$. Observe the mapping 
$$ {(X^+)}^{d-2i}: U_i \longrightarrow U_{d-i} ~~~~~ 
(0 \leq i \leq [d/2]) $$
is a bijection. By $X^+ W \subseteq W$, the image of $W_i$ by 
${(X^+)}^{d-2i}$ is included in $W_{d-i}$. 
Since dim\,$W_i=$ dim $W_{d-i}$, the mapping 
$$ {(X^+)}^{d-2i}: W_i \longrightarrow W_{d-i} ~~~~~ 
(0 \leq i \leq [d/2])$$
is a bijection. Consider the mapping 
$$ {(X^+)}^{d-2i-1}: W_{i+1} \longrightarrow W_{d-i}.$$
The subspace $X^+ W_i$ of $W_{i+1}$ is bijectively mapped onto 
$W_{d-i}$ by ${(X^+)}^{d-2i-1}$. So we have 
\begin{eqnarray}
\label{oplus}
W_{i+1} = X^+ W_i \oplus {\rm ker}\, {(X^+)}^{d-2i-1} |_{W_{i+1}}.
\end{eqnarray}
Since $ W_i=\bigoplus_{j=0}^i W_i^{(d-2j)}$ by the induction 
hypothesis for the claim (\ref{Wi}) and 
$X^+ W_i^{(d-2j)} \subseteq W_{i+1}^{(d-2j)}$, we have 
$$X^+ W_i \subseteq \bigoplus_{j=0}^i W_{i+1}^{(d-2j)}.$$
On the other hand, 
since ker\,${(X^+)}^{d-2i-1} |_{U_{i+1}} = U_{i+1}^{(d-2i-2)}$ 
$(i+1 \leq [d/2])$ by (\ref{ker}), we have 
$$ {\rm ker}\, {(X^+)}^{d-2i-1} |_{W_{i+1}} = W_{i+1}^{(d-2i-2)}.$$
Thus by (\ref{oplus}), we obtain  
$ W_{i+1} \subseteq \bigoplus_{j=0}^{i+1} W_{i+1}^{(d-2j)}$. 
Since the opposite inclusion 
is obvious, the claim holds for $i+1$ and we finish the proof of the 
claim (\ref{Wi}). 

\medskip
Since $W_i$ is bijectively mapped onto $W_{d-i}$ by 
${(X^+)}^{d-2i} ~ (0 \leq i \leq [d/2])$, it follows from 
(\ref{Wi}) that 
\begin{eqnarray}
\label{W d-i}
 W_{d-i} &=& \bigoplus_{j=0}^i W_{d-i}^{(d-2j)} ~~~
(0 \leq i \leq [d/2]), \\
\label{X d-2i}
 W_{d-i}^{(d-2j)}&=& {(X^+)}^{d-2i} W_i^{(d-2j)} ~~~
(0 \leq j \leq i \leq [d/2]). 
\end{eqnarray}
Define the subspace $W^{(d-2j)}$ by 
$$ W^{(d-2j)} = \bigoplus_{i=j}^{d-j} W_i^{(d-2j)} ~~~ 
(0 \leq j \leq [d/2]). $$
Then by (\ref{Wi}), (\ref{W d-i}), we obtain 
\begin{eqnarray}
\label{hom decomp}
W = \bigoplus_{j=0}^{[d/2]} W^{(d-2j)}. 
\end{eqnarray}
For $j \leq i < d-j$, the mappings 
$X^+:\, U_i^{(d-2j)} \longrightarrow U_{i+1}^{(d-2j)}$ and 
$X^-:\, U_{i+1}^{(d-2j)} \longrightarrow U_i^{(d-2j)}$ 
are inverses each other up to a nonzero scalar multiple. 
The image of $W_i^{(d-2j)}$ by $X^+$ is contained in 
$W_{i+1}^{(d-2j)}$,  in particular 
${\rm dim}\,W_i^{(d-2j)} \leq {\rm dim}\,W_{i+1}^{(d-2j)}$ 
($j \leq i < d-j$). On the other hand by (\ref{X d-2i}),  
${\rm dim}\,W_i^{(d-2j)} = {\rm dim}\,W_{d-i}^{(d-2j)}$ 
$(0 \leq j \leq i \leq [d/2])$. 
So ${\rm dim}\,W_i^{(d-2j)} = {\rm dim}\,W_{i+1}^{(d-2j)}$ 
($j \leq i < d-j$). 
Therefore the mapping 
$$ X^+:\, W_i^{(d-2j)} \longrightarrow W_{i+1}^{(d-2j)}$$ 
is a bijection for $j \leq i < d-j$ and 
the inverse of this mapping  
coincides with $X^-|_{ W_{i+1}^{(d-2j)} }$ 
up to a nonzero scalar multiple. Thus we obtain 
$X^- W_{i+1}^{(d-2j)} = W_i^{(d-2j)} ~(j < i+1 \leq d-j)$. 
Since $X^- W_j^{(d-2j)} \subseteq X^- U_j^{(d-2j)}=0$, it holds that 
$X^- W_i^{(d-2j)} \subseteq W ~
(0 \leq j \leq [d/2],\,j \leq i \leq d-j)$. 
Hence $X^- W \subseteq W$  by (\ref{hom decomp}) and 
the proof of Lemma \ref{lemma: X} is completed. 
\hfill $\Box$

\bigskip
\noindent
{\bf Proof of Theorem \ref{thm: 3rd kind} (i)}.  
Theorem \ref{thm: 3rd kind} (i)
is well-known but we give a brief proof as a warm-up. 
We are in the case of $(\varepsilon,\varepsilon^*)=(0,0)$ and given 
a $\T$-module 
$$ V=V(\ell_1,a_1) \otimes \cdots \otimes V(\ell_n,a_n) $$
via $\varphi_s$. 
To prove Theorem \ref{thm: 3rd kind} (i), 
it is enough to show that every $\T$-submodule 
$W$ of $V$ via $\varphi_s$ is $\LL$-invariant. 
The generators $x,\,y,\,k,\,k^{-1}$ of $\T$ act on $V$ via $\varphi_s$ 
as $\alpha s e_0^+,\, s^{-1} e_1^+,\, 
s\, k_0,\, s^{-1} k_1,$ i.e., $\T$ is embedded in $\LL$ 
via $\varphi_s$ as the Borel subalgebra genarated by 
$e_0^+,\, e_1^+,\, k_0^{\pm 1}$. 
$V$ has 
weight-space decomposition as in $(\ref{weight space decopm})$: 
$$ V=\bigoplus_{i=0}^d U_i ~~~~~ (k_0 |_{U_i} = q^{2i-d}), $$
where $d=\ell_1 + \cdots + \ell_n$. 
For a $\T$-submodule $W$ of the $\T$-module $V$ via $\varphi_s$, set 
$W_i=W \cap U_i$. Then 
$$ W = \bigoplus_{i=0}^d W_i. $$
Since the mapping 
${(e_0^+)}^{d-2i}:\, U_i \longrightarrow U_{d-i}$ 
is a bijection and 
$W_i$ is mapped into $W_{d-i}$ by ${(e_0^+)}^{d-2i}$, we have 
${\rm dim}\,W_i \leq {\rm dim}\, W_{d-i} ~ (0 \leq i \leq [d/2]).$ 
Since the mapping ${(e_1^+)}^{d-2i}:\, U_{d-i} \longrightarrow U_i$ is a 
bijection and $W_{d-i}$ is mapped into $W_i$ by ${(e_1^+)}^{d-2i}$, 
we have ${\rm dim}\, W_{d-i} \leq {\rm dim}\,W_i ~(0 \leq i \leq [d/2])$. 
Thus it holds that 
$$ {\rm dim}\, W_i = {\rm dim}\, W_{d-i} ~~~ (0 \leq i \leq d). $$

Consider the algebra homomorphism from the quantum algebra 
${\mathcal U}=U_q (sl_2)$ to the $U_q(sl_2)$-loop algebra $\LL$ that sends 
$X^+,\, X^-,\, K^{\pm 1}$ to $e_0^+,\, e_0^-,\, k_0^{\pm 1}$, respectively. 
Regard $V$ as a ${\mathcal U}$-module via this algebra homomorphism. 
Then $X^+ W \subseteq W,\, K W \subseteq W$. 
Since dim $W_i=$ dim $W_{d-i} ~(0 \leq i \leq d)$, 
we have by Lemma \ref{lemma: X} 
$X^- W \subseteq W$, i.e., $e_0^- W \subseteq W$. 
Similarly, consider the algebra homomorphism 
from ${\mathcal U}$ to $\LL$ that sends 
$X^+,\, X^-,\, K^{\pm 1}$ to $e_1^+,\, e_1^-,\, k_1^{\pm 1}$, 
respectively. 
Regard $V$ as a ${\mathcal U}$-module via this algebra homomorphism. 
Then the weight-space decomposition of this ${\mathcal U}$-module $V$ is 
$V= \bigoplus_{i=0}^d U_{d-i}$ 
$( K|_{U_{d-i}}=q^{2i-d})$, where $V=\bigoplus_{i=0}^d U_i$ 
is the weight-space decomposition of the $\LL$-module $V$. 
Since ${\rm dim}\,W_{d-i}={\rm dim}\,W_i ~(0 \leq i \leq d)$ and 
$X^+ W \subseteq W,\, K\,W \subseteq W$, 
we have by Lemma \ref{lemma: X}, $X^- W \subseteq W$, i.e., 
$e_1^- W \subseteq W$. 
Thus $W$ is $\LL$-invariant and the proof of 
Theorem \ref{thm: 3rd kind} (i) is completed. 
\hfill $\Box$

\medskip
We are now ready to prove the `if' part of 
Theorem \ref{thm: 2nd kind} (i). 

\medskip
\noindent
{\bf Proof of the `if' part of 
Theorem \ref{thm: 2nd kind} (i)}. 
We are in the case of $(\varepsilon,\varepsilon^*)=(1,0)$ and 
given a $\T$-module
$$ V=V(\ell) \otimes V(\ell_1,a_1) \otimes \cdots \otimes V(\ell_n,a_n) $$
via $\varphi_s$ such that $-s^{-2} \notin S(\ell_i,a_i)$ 
$(1 \leq i \leq n)$ 
and the multi-set ${\{ S(\ell_i,a_i) \}}_{i=1}^n$ is in general position. 
We want to show that the $\T$-module $V$ is irreducible. 
Note that $P_{V(\ell)} (s^{-2}) \neq 0,\, 
P_{V(\ell_i,a_i)} (s^{-2}) \neq 0$ by Theorem \ref{thm: product formula}, 
so 
\begin{eqnarray}
\label{sigma d V}
\sigma_d(V) \neq 0 
\end{eqnarray} 
by Remark \ref{rem: Drinfeld polynomial}, 
where $d=\ell+\ell_1+\cdots+\ell_n$.  
We may assume $n \geq 1$, otherwise $V=V(\ell)$ is obviously 
irreducible as a $\T$-module, 
since 
$\sigma_{\ell}(V) \neq 0$. 
Consider $i_0$ such that $a_{i_0} q^{\ell_{i_0}-1}$ is maximal in the set of 
scalars $a_i q^{\ell_i -1}$ $(1 \leq i \leq n)$ with respect to the partial 
ordering introduced in $(\ref{partial ordering})$ in Section 7.1. 
Among such $i_0' s$, choose one that 
has the smallest $\ell_{i_0}$. 
Since ${\{ S(\ell_i,a_i) \}}_{i=1}^n$ is in general position, the ordering 
of $V(\ell_i,a_i)$ $(1 \leq i \leq n)$ in the tensor product of $V$ does not 
matter about the isomorphism class of $V$ as an $\LL'$-module 
and consequently as a $\T$-module via $\varphi_s$. 
So we may assume $i_0=n$. 
Then 
\begin{eqnarray}
\label{ellq max 2nd kind}
&&a_n q^{\ell_n +1} \notin S(\ell_i,a_i) ~~~ (1 \leq i \leq n), \\
\label{ell min 2nd kind}
&&\ell_n \leq \ell_i ~{\rm if}~ a_n q^{\ell_n -1} \in S(\ell_i,a_i). 
\end{eqnarray}

\medskip
We proceed by induction on dim\,$V$ to prove that $V$ is irreducible 
as a $\T$-module. Set 
$$ V' = V(\ell) \otimes V(\ell_1,a_1) \otimes \cdots \otimes 
V(\ell_{n-1}, a_{n-1}). $$
Then $V=V' \otimes V(\ell_n, a_n)$. 
Since the $\LL$-module $V(\ell_n, a_n)$ is 
by Lemma \ref{lemma: ev module} 
embedded in the 
$\LL$-module $V(\ell_n -1,a_n q^{-1}) \otimes V(1,a_n q^{\ell_n -1})$ as 
the $\LL$-submodule spanned by the highest weight space, the $\LL'$-module 
$V$ can be regarded as embedded in the $\LL'$-module 
$$ \widetilde{V} = \bigl(V' \otimes V(\ell_n -1,a_n q^{-1})\bigr) \otimes 
V(1, a_n q^{\ell_n -1}), $$
sharing the highest weight space in common. 
We understand $V' \otimes V(\ell_n -1,a_n q^{-1})=V'$ if $\ell_n=1$. 
To apply Proposition \ref{prop: irr}  
to $\widetilde{V}$, 
we check the prerequisites for it, namely that 
$V' \otimes V(\ell_n -1, a_n q^{-1})$ is irreducible as a $\T$-module 
via $\varphi_s$, 
and that 
$\sigma_d (\widetilde{V}) \neq 0$ holds, 
where $d=\ell+\ell_1 + \cdots + \ell_n$ is the diameter of the $\T$-module 
$\widetilde{V}$ via $\varphi_s$. 
Observe that $S(\ell_n -1,a_n q^{-1})=S(\ell_n,a_n) 
\,\backslash \,\{a_n q^{\ell_n -1}\}$. 
So $-s^{-2} \notin S(\ell_i,a_i)$ $(1 \leq i \leq n-1)$ and 
$-s^{-2} \notin S(\ell_n -1,a_n q^{-1})$. 
Moreover the multi-set 
${\{ S(\ell_n,a_n) \}}_{i=1}^{n-1} \cup \{ S(\ell_n -1,a_n q^{-1})\}$ 
of $q$-strings associated with 
$V' \otimes V(\ell_n -1,a_n q^{-1})$ 
is in general position 
by (\ref{ellq max 2nd kind}), (\ref{ell min 2nd kind}). 
Therefore by induction on dimension, 
$V' \otimes V(\ell_n -1,a_n q^{-1})$ 
is irreducible as a $\T$-module via $\varphi_s$. 
Since $P_V(s^{-2}) \neq 0$ as we observed before, and since 
$P_{\widetilde{V}}(\lambda) = P_V(\lambda)$ 
by Theorem \ref{thm: product formula}, 
we have 
$P_{\widetilde{V}}(s^{-2}) \neq 0$, i.e., 
$\sigma_d (\widetilde{V}) \neq 0$. 
Thus the prerequisites are satisfied 
for Proposition \ref{prop: irr} 
to be applied  to $\widetilde{V}$. 
On the other hand, the conclusion of Proposition \ref{prop: irr} 
$$ P_{V' \otimes V(\ell_n -1,a_n q^{-1})}(- a_n q^{\ell_n +1}) = 0 $$
fails by Theorem \ref{thm: product formula}, 
since $a_n q^{\ell_n +1} \notin S(\ell_i,a_i)$ $(1 \leq i \leq n)$ 
by (\ref{ellq max 2nd kind}). 
This implies that any nonzero $\T$-submodules 
$W$ of $\widetilde{V}$ via $\varphi_s$ contains the highest weight space 
of $\widetilde{V}$. 
Since the $\T$-module $V$ via $\varphi_s$ is embedded in the $\T$-module 
$\widetilde{V}$ via $\varphi_s$, sharing the highest weight space in common, 
we conclude that any nonzero $\T$-submodule 
$W$ of $V$ contains the highest weight space of $V$. 

\medskip
Let $W$ be a minimal $\T$-submodule of the $\T$-module $V$ via $\varphi_s$. 
Note that $W$ is irreducible as a $\T$-module. 
Let $V = \bigoplus_{i=0}^d U_i$ denote the weight-space 
decomposition of the $\T$-module $V$ via $\varphi_s$. 
Then 
$$ W = \bigoplus_{i=0}^d W_i, ~~ W_i = W \cap U_i ~~ (1 \leq i \leq n).$$
By what we just proved in the last paragraph, we have $W_0 \neq 0$. 
Moreover $W_d \neq 0$ by $(\ref{sigma d V})$. 
Since ${\rm dim}\,U_0={\rm dim}\,U_d=1$, we obtain 
$W_0=U_0$, $W_d=U_d$. 
We claim 
\begin{eqnarray}
\label{Wi W d-i}
W_i = W_{d-i} ~~~~~ (0 \leq i \leq d). 
\end{eqnarray}
Let $\A$ denote the TD-algebra for $(\varepsilon,\varepsilon^*)=(1,0)$. 
Consider $\varphi_s \circ \iota_t : \A \longrightarrow \LL'$ 
and regard $V$ as an 
$\A$-module via $\varphi_s \circ \iota_t$. 
By Theorem \ref{thm: C2} and (\ref{C1zt}), 
the generators $z,\,z^*,$ of $\A$ act on $W$ as a TD-pair, 
if we choose $t$ suitably. 
The split decomposition of $W$ for the TD-pair coincides 
with the weight-space decomposition of $W$.  
Thus we obtain (\ref{Wi W d-i}) by \cite[Corollary 5.7]{ITT}.

\medskip
The generators $x,\, y,\, k,\, k^{-1}$ of $\T$ act on $V$ via $\varphi_s$ as 
$\alpha (s e_0^+ + s^{-1} e_1^- k_1),\, 
s^{-1} e_1^+,\, s k_0,\, s^{-1} k_1$ respectively. 
Consider the algebra homomorphism from ${\mathcal U}$ to $\LL'$ 
that sends $X^+, X^-, K^{\pm 1}$ to $e_1^+, e_1^-, k_1^{\pm 1}$. 
Regard $V$ as a ${\mathcal U}$-module via this algebra homomorphism. 
Then the weight-space 
decomposition of this ${\mathcal U}$-module $V$ is 
$ V = \bigoplus_{i=0}^d U_{d-i}$ $(K|_{U_{d-i}} = q^{2i-d})$, 
where 
$V = \bigoplus_{i=0}^d U_i$ is 
the weight space-decomposition of the $\T$-module $V$ via $\varphi_s$. 
Since dim\,$W_{d-i} = $dim\,$W_i$ $(0 \leq i \leq d)$ 
by $(\ref{Wi W d-i})$ 
and 
$X^+ W \subseteq W,\, K\,W \subseteq W$, 
we have by Lemma \ref{lemma: X}  
$X^- W \subseteq W$, i.e., $e_1^- W \subseteq W$. 
Since $x W \subseteq W$, i.e., 
$(s e_0^+ + s^{-1} e_1^- k_1) W \subseteq W$, we obtrain 
$e_0^+ W \subseteq W$ by $e_1^- k_1 W \subseteq W$. 
Thus $W$ is $\LL'$-invariant. 
Recall that we have already shown 
that $W$ cantains the highest weight space $U_0$ of 
the $\T$-module $V$ via $\varphi_s$. 
By the following lemma, we obtain $W=V$ and the 
`if' part of Theorem \ref{thm: 2nd kind} (i) is completed. 

\begin{lemma}
\label{lemma: L'}
Assume that a multi-set 
$\{ S (\ell_i,a_i)\}_{i=1}^n$ of q-strings is in general position. 
Consider the $\LL'$-module 
$$ V = 
V(\ell) \otimes V(\ell_1,a_1) \otimes \cdots \otimes V(\ell_n,a_n) $$
and let
$$ V = \bigoplus_{i=0}^d U_i, \qquad k_0|_{U_i} =  q^{2i-d}$$
be the eigenspace decomposition of $k_0$, where 
$d=\ell+\ell_1 + \cdots + \ell_n$. 
If $W$ is an $\LL'$-submodule of $V$ and contains $U_0$, 
then $W=V$. 
\end{lemma}
Proof. 
Set 
$$ V' = V(\ell_1,a_1) \otimes \cdots \otimes V(\ell_n,a_n) $$
and let $\B$ denote the subalgebra of $\LL'$ generated by 
$e_0^+,\, e_1^+,\, k_0^{\pm 1}$. Note that $V'$ is irreducible as an 
$\LL'$-module, since it is already irreducible as a $\B$-module 
by Theorem \ref{thm: 3rd kind} (i). 
We may assume $\ell \geq 1$, since if $\ell=0$, then $V=V'$ and 
the $\LL'$-module $V$ is irreducible. 

\medskip
Let $V(\ell)=\langle v_0,v_1, \cdots, v_{\ell} \rangle$ be 
a standard basis as an $\LL'$-module: 
$e_0^+ v_i=0,\  
e_1^+ v_i = [\ell-i+1] v_{i-1}, \,
e_1^- v_i=[i+1]\, v_{i+1}, \,
k_0 v_i= q^{2i-\ell} v_i$ 
$(0 \leq i \leq \ell)$, where $v_{-1}= v_{\ell+1}=0$. Then 

$$ V = \bigoplus_{i=0}^{\ell} \langle v_i \rangle \otimes V'. $$
We show $W \supseteq \langle v_i \rangle 
\otimes V' ~~ (0 \leq i \leq \ell)$ by induction on $i$. 
For $i=0$, some element 
$$v_0 \otimes v' ~~~(V' \ni v' \neq 0)$$ 
is containd in $W$ by $W \supseteq U_0$. Since 
$e_0^+ (v_0 \otimes v') = q^{-\ell} v_0 \otimes (e_0^+ v')$,  
$e_1^+ (v_0 \otimes v')=q^{\ell} v_0 \otimes (e_1^+ v')$ and 
$k_0^{\pm 1} (v_0 \otimes v') = q^{\mp \ell} v_0 \otimes (k_0^{\pm 1} v')$,
it follows from $\B W \subseteq W$ that 
$v_0 \otimes (e_0^+ v')$, $v_0 \otimes (e_1^+ v')$, 
$v_0 \otimes (k_0^{\pm} v')$ are contained in $W$. 
Since the elements $ e_0^+,\, e_1^+,\, k_0^{\pm 1}$ generate $\B$, 
we obtain 
$$ \langle v_0 \rangle \otimes \B v' \subseteq W.$$
Since $V'$ is irreducible as a $\B$-module by 
Theorem \ref{thm: 3rd kind} (i), 
we have $\B v' = V'$ so 
$v_0 \otimes V' \subseteq W$. 
Suppose 
that $\langle v_i \rangle
\otimes V' \subseteq W$. Choose a nonzero element $v'$ from $V'$. 
Then $e_1^- (v_i \otimes v') = [i+1]\, v_{i+1} \otimes (k_1^{-1} v') 
+ v_i \otimes (e_1^- v')$. 
Since $e_1^- (v_i \otimes v')$ and $v_i \otimes (e_1^- v')$ are contained 
in $W$, we have 
$$ v_{i+1} \otimes v'' \in W, $$
where $v'' = k_1^{-1} v' \neq 0$. So $e_0^+ (v_{i+1} \otimes v''), \,
e_1^+ (v_{i+1} \otimes v''),\, k_0 (v_{i+1} \otimes v'')$ are contained 
in $W$. 
Since $e_0^+ (v_{i+1} \otimes v'') = 
q^{2i+2-\ell} v_{i+1} \otimes (e_0^+ v''), \,
e_1^+ (v_{i+1} \otimes v'') = [\ell-i]\, v_i \otimes v'' + 
q^{\ell-2i-2} v_{i+1} \otimes ( e_1^+ v''),\, 
 k_0 (v_{i+1} \otimes v'') = q^{2i+2-\ell} 
v_{i+1} \otimes (k_0 v'')$, it follows from 
$v_i \otimes v'' \in W$ that 
$$ v_{i+1} \otimes (e_0^+ v''),\, v_{i+1} \otimes (e_1^+ v''), \,
v_{i+1} \otimes (k_0 v'') $$
are all contained in $W$. So 
$$ \langle v_{i+1}\rangle \otimes \B v'' \subseteq W.$$
Since $\B v'' = V'$, we have 
$\langle v_{i+1}\rangle \otimes V' \subseteq W$. 
This completes the proof of Lemma \ref{lemma: L'}. 
\hfill $\Box$

\subsection {Proof of part (iii)} 

The part (iii) of 
Theorem \ref{thm: 1st kind}, 
Theorem \ref{thm: 2nd kind}, 
Theorem \ref{thm: 3rd kind} 
follows from the part (i) 
together with 
Theorem $\ref{thm: sigma}'$, 
Theorem \ref{thm: product formula},  
and some combinatorial observations as in   
Lemma~\ref{lemma: q-string}; we prove 
Lemma \ref{lemma: q-string} 
at the end of this subsection separately. 
Let $s$ and $d$ be a nonzero scalar and a positive integer respectively, 
chosen arbitrarily. 
We are given a polynomial $P(\lambda)$ in ${\mathcal P}_d^s$, i.e., 
$P(\lambda)$ is a monic polynomial of degree $d$ such that 
$P(\varepsilon s^{-2} + \varepsilon^* s^2) \neq 0$. 
We want to construct an irreducible $\T$-module $V$ via 
$\varphi_s$ such that the Drinfel'd polynomial $P_V(\lambda)$ 
coincides with $P(\lambda)$. 
Let $\lambda_1, \lambda_2, \cdots, \lambda_d$ denote the roots of 
$P(\lambda)$, allowing repetition.

\medskip
If $(\varepsilon,\varepsilon^*)=(1,1)$, let $\Omega_i$ denote 
the set of solutions of 
$$ \lambda_i + \zeta + \zeta^{-1} = 0 $$
for $\zeta$. We understand that $\Omega_i$ is a multi-set if 
$\lambda_i= \pm 2$. So $|\Omega_i|=2$ $(1 \leq i \leq d)$.
Set 
$$ \Omega = \bigcup_{i=1}^d\, \Omega_i $$
as a multi-set. Then $|\Omega|=2d$ as a multi-set. 
By Lemma \ref{lemma: q-string}, 
there exists a multi-set 
${\{ S(\ell_i,a_i) \}}_{i=1}^n$ of $q$-strings strongly in 
general position such that 
$$ \Omega = \bigcup_{i=1}^n \,
\bigl(S(\ell_i,a_i) \cup S(\ell_i,a_i^{-1})\bigr) $$
as multi-sets. Since $|S(\ell_i,a_i)|=\ell_i$, we have 
$d=\ell_1 + \cdots + \ell_n$. The $\T$-module 
$$ V = V(\ell_1,a_1) \otimes \cdots \otimes V(\ell_n,a_n) $$
via $\varphi_s$ has Drinfel'd polynomial 
$$ P_V (\lambda) = \prod_{i=1}^n P_{V(\ell_i,a_i)} (\lambda) $$
by Theorem \ref{thm: product formula}, where 
$$ P_{V(\ell_i,a_i)} (\lambda) = \prod_{\zeta \in S(\ell_i,a_i)} 
(\lambda + \zeta + \zeta^{-1}). $$
Thus $P_V (\lambda) = P (\lambda)$. Since 
$P(s^{-2} + s^2) \neq 0$, 
we have $-s^2 \notin S(\ell_i,a_i) \cup S(\ell_i,a_i^{-1})$ 
$(1 \leq i \leq n)$. 
So by Theorem \ref{thm: 1st kind} (i), 
the $\T$-module $V$ via $\varphi_s$ is irreducible. 

\medskip
If $(\varepsilon,\varepsilon^*)=(1,0)$, set 
$$ \Omega = 
\{ - \lambda_i ~|~ \lambda_i \neq 0,\, 1 \leq i \leq d \} $$ 
as a multi-set. 
We may assume that 
$\Omega = 
\{ - \lambda_{i}~|~ \ell+1 \leq i \leq d \}$
and 
$\lambda_1 = \cdots = \lambda_{\ell} = 0$, 
allowing $\ell=0$ . 
It is well-known and easy to show that there exists a multi-set 
${\{ S(\ell_i,a_i) \}}_{i=1}^n$ of $q$-strings 
in general position such that 
$$ \Omega =\, \bigcup_{i=1}^n S(\ell_i,a_i) $$
as multi-sets. Since $|S(\ell_i,a_i)|=\ell_i$, 
we have $d-\ell=\ell_1 + \cdots + \ell_n$. 
The $\T$-module 
$$ V = V(\ell) \otimes 
V(\ell_1,a_1) \otimes \cdots \otimes V(\ell_n,a_n) $$
via $\varphi_s$ has Drinfel'd polynomial 
$$ P_V (\lambda) = \lambda^{\ell} 
\prod_{i=1}^n P_{V(\ell_i,a_i)} (\lambda) $$
by Theorem \ref{thm: product formula}, where 
$$ P_{V(\ell_i,a_i)} (\lambda) = 
\prod_{c \in S(\ell_i,a_i)} (\lambda + c).$$
Thus $P_V (\lambda) = P (\lambda)$. 
Since $P(s^{-2}) \neq 0$, we have 
$-s^{-2} \notin S(\ell_i,a_i)$ $(1 \leq i \leq n)$. 
So by Theorem \ref{thm: 2nd kind} (i), 
the $\T$-module $V$ via $\varphi_s$ is irreducible. 

\medskip
If $(\varepsilon,\varepsilon^*)=(0,0)$, set
$$ \Omega = \{ - \lambda_1, \lambda_2, \cdots, - \lambda_d \}, $$
as a multi-set. 
Since $P(\varepsilon s^{-2} + \varepsilon^* s^2)= P(0) \neq 0$, 
we have $\lambda_i \neq 0$  $(1 \leq i \leq n)$. 
There exists a multi-set of $q$-strings 
${\{ S(\ell_i,a_i) \}}_{i=1}^n$ in general position such that 
$$ \Omega = \bigcup_{i=1}^n S(\ell_i,a_i) $$
as multi-set.  Since  $S|(\ell_i,a_i)|=\ell_i$, we have 
$d=\ell_1 + \cdots + \ell_n$. The $\T$-module 
$$ V = V(\ell_1,a_1) \otimes \cdots \otimes V(\ell_n,a_n) $$
via $\varphi_s$ has Drinfel'd polynomial 
$$ P_V (\lambda) = \prod_{i=1}^n P_{V(\ell_i,a_i)} (\lambda)$$
by Theorem \ref{thm: product formula}, where 
$$ P_{V(\ell_i,a_i)} (\lambda) = 
\prod_{c \in S(\ell_i,a_i)} (\lambda + c). $$
Thus $P_V (\lambda) = P (\lambda)$. 
By Theorem \ref{thm: 3rd kind} (i), the $\T$-module 
$V$ via $\varphi_s$ is irreducible. This completes the proof of 
the part (iii) for 
Theorem \ref{thm: 1st kind}, Theorem \ref{thm: 2nd kind}, 
Theorem \ref{thm: 3rd kind}.

\bigskip 
\noindent
{\bf Proof of Lemma \ref{lemma: q-string}}. 
We proceed by induction on $|\Omega|$, where $|\Omega|$ denotes the 
number of elements in $\Omega$, counting the multiplicities. 
Recall the partial ordering $(\ref{partial ordering})$
on $\C ~\backslash~ \{ 0 \}$ introduced in Section 7.1: 
$$a \leq b ~\Leftrightarrow~ 
b = a\, q^{2i} ~\mbox{for some integer}~ i \geq 0.$$
Choose a maximal element c in $\Omega$ 
with respect to this partial ordering. 
Note that $c^{-1}$ is minimal in $\Omega$. Set 
$$ \Omega' = \Omega ~\backslash~ \{ c, c^{-1} \} $$
as multi-sets  of nonzero scalars. 
Then by induction, there exists a multi-set 
${\{ S(\ell_i',a_i')\}}_{i=1}^{n'}$ of $q$-strings strongly in general 
position such that 
$$ \Omega' = 
\bigcup_{i=1}^{n'} \, \bigl( S(\ell_i',a_i') 
\cup  S(\ell_i',{a_i'}^{-1}) \bigr) $$ 
as multi-sets of nonzero scalars. 
Moreover such a multi-set  
${\{ S(\ell_i', a_i') \}}_{i=1}^{n'}$ of $q$-strings 
is uniquely determined by $\Omega'$ up to equivalence. 
Observe that the union $S(\ell_i', a_i') \cup \{ c \}$ 
(resp. $S(\ell_i',a_i') \cup \{ c^{-1})$ ) 
as a multi-set of nonzero scalars 
is a $q$-string if and only if 
$c = a_i' q^{\ell_i' +1}$ (resp. $c^{-1} = a_i' q^{-\ell_i' -1})$, 
in which case 
$$
S(\ell_i',a_i') \cup \{c\} = S(\ell_i' +1, a_i' q)~~~
\bigl( {\rm resp.}~ 
S(\ell_i', a_i') \cup \{c^{-1}\} = S (\ell_i' +1, a_i' q^{-1}) \bigr).
$$

\medskip
If there exist $i$'s such that either 
$c=a_i' q^{\ell_i'+1}$ or $c^{-1}=a_i' q^{-\ell_i'-1}$, 
choose one among such $i$'s that has the largest $\ell_i'$. 
By rearranging the ordering of the $q$-strings, 
we may assume $i=n'$. 
By replacing $a_{n'}'$ by $a_{n'}'^{-1}$ 
if $c^{-1}=a_{n'}' q^{-\ell_{n'}'-1}$, 
we may assume $c=a'_{n'} q^{\ell_{n'}'+1}$. 
Thus $c=a_{n'}' q^{\ell_{n'}'+1}$ is 
maximal in $\Omega$ and if $c=a_i' q^{\ell_i'+1}$ or 
$c^{-1}=a_i' q^{-\ell_i'-1}$ holds for some $i$, 
then $\ell_n' \geq \ell_i'$. In this case, define $q$-strings  
$S(\ell_i,a_i)$  $(1 \leq i \leq n')$ by 
\begin{eqnarray}
\label{S(ell i, a i) case 1}
S(\ell_i,a_i) &=& S(\ell_i',a_i') ~~~~~ (1 \leq i \leq n'-1), \\
\label{S(ell n, a n) case 1}
S(\ell_{n'},a_{n'}) &=& S(\ell_{n'}'+1,a_{n'}'q). 
\end{eqnarray}
Then the multi-set ${\{ S(\ell_i,a_i)\}}_{i=1}^{n'} $ of 
$q$-strings is strongly in general position and 
$$ \Omega = \bigcup_{i=1}^{n'}\, 
\bigl( S(\ell_i,a_i) \cup S(\ell_i,a_i^{-1}) \bigr)$$
as multi-set of nonzero scalars. 

\medskip
If there exist no $i$'s such that either $c=a_i' q^{\ell_i'+1}$ or 
$c^{-1}=a_i' q^{-\ell_i'-1}$, then 
define $q$-strings $S(\ell_i,a_i)$  $(1 \leq i \leq n'+1)$ by 
\begin{eqnarray}
\label{S(ell i, a i) case 2}
S(\ell_i,a_i) &=& S(\ell_i',a_i') ~~~~~ (1 \leq i \leq n'), \\
\label{S(ell n, a n) case 2}
S(\ell_{n'+1},a_{n'+1}) &=& S(1,c). 
\end{eqnarray}
Then the multi-set ${\{ S(\ell_i,a_i)\}}_{i=1}^{n'+1} $ of 
$q$-strings is strongly in general position and 
$$ \Omega =  \bigcup_{i=1}^{n'+1}\, 
\bigl( S(\ell_i,a_i) \cup S(\ell_i,a_i^{-1}) \bigr) $$
as multi-sets of nonzero scalars.
In any case, there exists a desired multi-set of $q$-strings.

\medskip
Next we show the 
uniqueness of such a multi-set of $q$-strings up to equivalence. 
Let 
$\{S(m_i,b_i)\}_{i=1}^n $ 
be a multi-set of $q$-strings 
strongly in general position such that 
$$ \Omega =  \bigcup_{i=1}^n\, 
\bigl( S(m_i,b_i) \cup S(m_i,b_i^{-1}) \bigr). $$
Then the maximal element $c$ of $\Omega$, 
which was chosen in the course of the construction 
of a desired multi-sets of $q$-strings,  
belongs to either $S(m_i,b_i)$ or 
$S(m_i,b_i^{-1})$ for some $i$. 
Among such $i$'s, choose one that has the smallest 
$m_i$. We many assume $i=n$ and 
$c \in S(m_n,b_n)$ by rearranging the ordering of 
the $q$-strings and replacing $b_n$ by $b_n^{-1}$ if necessary. 
Thus $b_n q^{m_n-1}$ is the maximal element $c$ 
and $m_n \leq m_i$ holds if $c \in S(m_i,b_i)$ or 
$c \in S(m_i,b_i^{-1})$, i.e., $c=b_i q^{m_i-1}$ or 
$c^{-1}=b_i q^{-m_i+1}$. 

\medskip
If $m_n \geq 2,$ then the multi-set 
$$ {\{ S(m_i,b_i)\}}_{i=1}^{n-1} \cup \{ S(m_n -1 ,b_n q^{-1})\}$$
of $q$-strings is strongly in general position and 
covers the multi-set 
$\Omega' = \Omega \,\backslash \, \{ c,c^{-1}\}$ 
of nonzero scalars as the union of 
$S(m_i,b_i),\,S(m_i,b_i^{-1})$ $(1 \leq i \leq n-1)$ and 
$S(m_n -1,b_n q^{-1}),\,S(m_n -1,b_n^{-1} q)$. 
Such a multi-set of $q$-strings is unique up to 
equivalence by induction. 
So the multi-set 
${\{ S(m_i,b_i)\}}_{i=1}^{n-1} \cup \{ S(m_n -1 ,b_n q^{-1})\}$ 
of $q$-strings 
is equivalent to 
$\{S(\ell'_i,a'_i)\}_{i=1}^{n'}$, the one which was chosen 
in the course of the construction of a desired multi-sets of $q$-strings.  
Observe that $c=(b_n q^{-1}) q^{(m_n-1)+1}$ and 
$m_n -1 \geq m_i$ if $c=b_i q^{m_i+1}$ or $c^{-1}=b_i q^{-m_i-1}$ 
for some $i$ $(1 \leq i \leq n-1)$, since $S(m_n, b_n)$ includes either 
$S(m_i, b_i)$ or $S(m_i, b_i^{-1})$ for such an $i$. 
Thus we have $n=n'$ and we may assume 
\begin{eqnarray*}
& & S(m_i,b_i) = S(\ell_i',a_i') ~~~~~ (1 \leq i \leq n-1), \\
& & S(m_n-1,b_n q^{-1}) = S(\ell_n',a_n'). 
\end{eqnarray*}
By $(\ref{S(ell i, a i) case 1})$, 
$(\ref{S(ell n, a n) case 1})$,  
the multi-set $\{S(m_i,b_i)\}_{i=1}^n $ of $q$-strings is 
equivalent to the one we constructed 
by means of $\{S(\ell'_i,a'_i)\}_{i=1}^{n'} $. 

\medskip
If $m_n=1$, then $S(m_n,b_n)=\{ c \}$ and the multi-set 
$$ {\{ S(m_i,b_i) \}}_{i=1}^{n-1} $$
of $q$-strings is strongly in general position and 
covers 
the multi-set $\Omega' = \Omega \,\backslash \, \{ c,c^{-1}\}$ 
as a union of $S(m_i,b_i),\, S(m_i,b_i^{-1})$ $(1 \leq i \leq n-1)$. 
Such a multi-set of $q$-strings is unique up to equivalence. 
Observe that there exist no $i$'s $(1 \leq i \leq n-1)$ such that 
either $c=b_i q^{m_i+1}$ or $c^{-1}=b_i q^{-m_i-1}$, 
since otherwise $S(m_n,b_n) \cup S(m_i,b_i)$ or 
$S(m_n,b_n^{-1}) \cup S(m_i,b_i)$ would be a $q$-string for such 
an $i$. 
Thus
we have $n'=n-1$ and we may assume 
$$ S(m_i,b_i) = S(\ell_i',a_i') ~~(1 \leq i \leq n-1).$$
By $(\ref{S(ell i, a i) case 2})$, 
$(\ref{S(ell n, a n) case 2})$, 
the multi-set $\{S(m_i,b_i)\}_{i=1}^n $ of $q$-strings is 
equivalent to the one we constructed 
by means of $\{S(\ell'_i,a'_i)\}_{i=1}^{n'} $. 
This completes the proof of Lemma \ref{lemma: q-string}. 
\hfill $\Box$

\medskip
\noindent
Tatsuro Ito \hfil\break
Division of Mathematical and Physical Sciences \hfil\break
Kanazawa University \hfil\break
Kakuma-machi, Kanazawa 920--1192, Japan \hfil\break
E-mail: tatsuro@kenroku.kanazawa-u.ac.jp  \hfil\break

\medskip
\noindent
Paul Terwilliger \hfil\break
Department of Mathematics\hfil\break
University of Wisconsin-Madison \hfil\break
Van Vleck Hall \hfil\break
480 Lincoln drive  \hfil\break
Madison, WI 53706-1388, USA \hfil\break
E-mail: terwilli@math.wisc.edu \hfil\break

\end{document}